\documentclass[leqno,12pt]{article}
\usepackage{amsmath, amsthm}
\usepackage{amsfonts}
\usepackage{amssymb}
\usepackage{graphicx}
\usepackage{color}
\usepackage{hyperref}

\newtheorem{theorem}{Theorem}[section]
\newtheorem{corollary}[theorem]{Corollary}
\newtheorem{lemma}[theorem]{Lemma}
\newtheorem{proposition}{Proposition}[section]

\definecolor{blue}{rgb}{0.05, .05, .65}

\newcommand{\re}{\mathbb{R}}
\newcommand{\ren}{\mathbb{R}^N}
\newcommand{\ve}{\varepsilon}

\newcommand{\dy}{\,{\rm d}y}

\numberwithin{equation}{section}
\def\qed{\,\unskip\kern 6pt \penalty 500
\raise -2pt\hbox{\vrule \vbox to8pt{\hrule width 6pt
\vfill\hrule}\vrule}\par}

\definecolor{ao}{rgb}{0.0, 0.0, 1.0}

\newcommand{\nc}{\normalcolor}

\definecolor{darkblue}{rgb}{0.05, .05, .65}
\definecolor{darkgreen}{rgb}{0.1, .65, .1}
\definecolor{darkred}{rgb}{0.8,0,0}
\title{\bf The  fractional p-Laplacian evolution equation in $\ren$ in the sublinear case
}
\author{ \sc \large J. L. V\'azquez\footnote{Univ. Aut\'onoma de Madrid, Spain} \footnote{ Real Academia de Ciencias, Spain.}}

\date{ }
\begin{document}
\maketitle

\begin{abstract} We consider the  natural  time-dependent fractional $p$-Laplacian equation posed in the whole Euclidean space, with parameter $1<p<2$ and  fractional exponent $s\in (0,1)$.  Rather standard theory shows that the  Cauchy Problem for data in the Lebesgue $L^q$ spaces is well posed, and the solutions form a family of  non-expansive semigroups with regularity and other interesting properties. The superlinear case $p>2$ has been dealt with in a recent paper.

We study here the ``fast'' regime $1<p<2$ which is more complex. As  main results, we construct the self-similar fundamental solution for every mass value $M$ and any $p$ in the subrange $p_c=2N/(N+s)<p<2$, and we show that this is the precise range  where they can exist. We also prove that general finite-mass solutions converge  towards the fundamental solution having the same mass, and convergence holds in all $L^q$ spaces. Fine bounds in the form of global Harnack inequalities are obtained.

Another main topic of the paper is the study of solutions having strong singularities. We  find a type of singular solution called Very Singular Solution that exists for $p_c<p<p_1$, where $p_1$ is a new critical number that we introduce, $p_1\in (p_c,2)$. We extend this type of singular solutions to the ``very fast range'' $1<p<p_c$. They represent examples of weak solutions having finite-time extinction in that lower $p$ range. We briefly examine the situation in the limit case $p=p_c$. Finally, very singular solutions are related to fractional elliptic problems of nonlinear eigenvalue form.
\end{abstract}
	
\

\vskip 1cm
\newpage

\small
\tableofcontents
\normalsize

\newpage

%
%

\section{Introduction. The problem}
\normalsize

This is a companion to paper \cite{VPLP2020} dealing with the natural fractional $p$-Laplacian evolution equation posed in the whole Euclidean space with fractional exponent $s\in (0,1)$. While the previous paper treated  the superlinear parameter case $p>2$, here we consider the sublinear case $1<p<2$ \ in the same setting. We point out that even if the very general facts of the theory bear a similar flavor, the detailed analysis of some important aspects leads to quite different results and methods. Thus, the long time behaviour, which is a main concern of both papers, offers some similarities along with remarkable novelties as $p$ goes down. And we find new types of solutions, like the very singular solutions (VSS) which play an important role in the theory.

For the reader's convenience we begin our present study by a brief review of the basic facts explained in \cite{VPLP2020}. The nonlocal energy functional
\begin{equation}\label{Jsp1}
{\mathcal J}_{p,s}(u)= \frac1{p}\int_{\ren }\int_{\ren } \frac{|u(x)-u(y)|^p}{|x-y|^{N+sp}}\,dxdy
\end{equation}
is a power-like functional with nonlocal kernel of the $s$-Laplacian type. It has attracted a great deal of attention in recent years. It is just the $p$-power of the Gagliardo seminorm, used in the definition of the $W^{s,p}$ spaces, \cite{AdFou, DiNPV}. We may consider it  for   exponents $0<s<1$ and $1<p<\infty$ in dimensions $N\ge 1$. Its  subdifferential 
${\mathcal L}_{s,p}$  is the nonlinear operator defined  a.e. in its domain by the formula
\begin{equation}\label{frplap.op}
\displaystyle \qquad {\mathcal L}_{s,p}(u):= P.V.\int_{\ren}\frac{\Phi(u(x,t)-u(y,t))}{|x-y|^{N+sp}}\,dy\,,
\end{equation}
where $\Phi(z)=|z|^{p-2}z,$ and $P.V.$ means principal value. It is a usually called the $s$-fractional $p$-Laplacian operator. It is  well-known from general theory that ${\mathcal L}_{s,p}$ is a maximal monotone operator in $L^2(\ren)$ with dense domain.
As in \cite{VPLP2020}, we will study the corresponding gradient flow, i.e., the homogeneous evolution equation
\begin{equation}\label{frplap.eq}
\partial_t u + {\mathcal L}_{s,p} u=0\,,
\end{equation}
posed in the Euclidean space $x\in \ren$, $ N\ge 1$, for $t>0$. We refer to it as the {\sl fractional $p$-Laplacian evolution equation}, FPLEE for short.  Motivation and related equations for this model can be seen in the  \cite{VPLP2020} and its references. See also \cite{MazRT},\cite{Puhst} and \cite{Vaz16}, where the equation is posed on a bounded domain.
In this paper we always take $p<2.$  Note that the case $p=2$ corresponds to the fractional linear heat equation, $u_t+(-\Delta)^{s} u=0,$ that has a well-established theory, \cite{BPSV, Blumenthal-Getoor, BSV17, Kwas2017, Va18}.

We supplement the equation with an initial datum \begin{equation}\label{init.data}
\lim_{t\to 0} u(x,t)= u_0(x),
\end{equation}
where a most common choice is $u_0\in L^2(\ren)$. However, the theory shows that Equation \eqref{frplap.eq} generates a continuous nonlinear semigroup in any $L^q(\ren)$ space, $1\le q<\infty$, in fact a nonexpansive semigroup in any of those spaces. As in the companion paper, we  concentrate  on data $u_0\in L^1(\ren)$, which leads to the class of finite-mass solutions and  produces a specially rich theory. In this paper we still consider all fractional exponents $0<s<1$, but the asymptotic theory of finite-mass solutions forces us to further restrict the range $1<p<2$. It happens that the asymptotic behaviour offers a reasonably unified aspect only if
\begin{equation}\label{pc.range}
p_c=\frac{2N}{N+s}<p<2,
\end{equation}
loosely speaking, when $p$ is  close to 2. Note that the gap $2-p_c=2s/(N+s)$ shrinks to 0 as $s\to 0$. This limited range will be considered in the bulk of the paper since its theory offers a certain unity with the theory for $p\ge 2$ and is quite rich.

We will be specially interested in taking a Dirac delta as initial datum. In that case the solution is called a {\sl fundamental solution \rm (FS, also called a {\sl source-type solution}). For $p>2$ such solutions were constructed in  \cite{VPLP2020}. We will show that the range \eqref{pc.range} is optimal for existence, since there are no fundamental solutions in the range  $1<p\le p_c$,  as we will explain later in the paper. The restricted range will be further divided into $p_c<p<p_1$ and $p_1<p<2$ because  the detailed analysis shows strong differences in the  fine behaviour. The new critical exponent $p_1$ is characterized below.

The reader is here reminded that the strong qualitative and quantitative separation between the two exponent ranges, $p>2$ and $p<2$ is a key feature of the non fractional $p$-Laplacian equation, i.e., \ $u_t=\nabla\cdot(|\nabla u|^{p-2}\nabla u)$, that we will call in the sequel the {\sl standard} $p$-Laplacian equation, SPLE, for the sake of clarity.
Therefore, it is no wonder we find traces of it here in the fractional setting. We recall that in the standard equation  the range $p>2$ is called the \sl slow gradient-diffusion case \rm (with finite speed of propagation and free boundaries),  while the range $1<p<2$ is called the \sl fast gradient-diffusion case \rm (with infinite speed of propagation), cf. \cite{DiBeBk}, \cite[Section 11]{VazSmooth}. Such denominations do not apply here in a strict sense, but many features have a kind of correspondence here, while other features are subject to remarkable changes, as we will see. We will study the agreement of our results with the ones for the standard equation when $s\to 1$, and with the linear ones when $p\to 2$.

The diffusion range with lower values of $p$, $1<p<p_c$, called ``very fast diffusion'' in the standard case, has very peculiar properties. We will refer to it  only  at the end of the paper to establish some remarkable properties and  comment on its main differences \nc with the case $p>p_c$. Attention is paid to the singular solutions we call Very Singular Solutions, VSS. Much work remains to be done otherwise.


\subsection{Outline of the paper and main results}

 We focus on Problem \eqref{frplap.eq}-\eqref{init.data},  posed in $\ren$  with $p$ in the range \eqref{pc.range}. It is known that this Cauchy problem is well-posed in all $L^q(\ren)$ spaces,  $1\le q< \infty$.
This parallels what is known in the case of bounded domains or the case $p\ge 2$. In view of such works, we first review the main facts of the theory in Section \ref{sec.basic}. In particular, we  define the class of continuous strong solutions that correspond to $L^2$ and $L^1$  initial data and derive its main properties. To note, we prove that bounded strong semigroup solutions are H\"older continuous using the elliptic theory of \cite{IannMS2016}, cf. Theorem \ref{them.reg.new}. See full details below.

\medskip

\noindent \bf Fundamental solutions. \rm As already explained in \cite{VPLP2020}, a most interesting question is that of finding the FS, i.e., the solution such that
\begin{equation}
\lim_{t\to 0} \int_{\ren} u(x,t)\varphi(x)\,dx=M\varphi(0),
\end{equation}
for every smooth and  compactly supported test function $\varphi\ge 0$, and some $M>0$.  This is our first main contribution to the question.

\begin{theorem}\label{thm.ssfs} Let $p_c<p<2.$ For every given mass $M>0$ there exists a unique self-similar solution of Problem \eqref{frplap.eq}-\eqref{init.data}. It has the form
\begin{equation}\label{eq.ssf1}
U(x,t;M)=M^{sp\beta}t^{-\alpha}F(M^{(2-p)\beta} |x|\,t^{-\beta})\,,
\end{equation}
with self-similarity exponents
\begin{equation}\label{eq.sse1}
\alpha=\beta N, \quad \beta=\frac{1}{sp-N(2-p)}\,.
\end{equation}
The profile $F(r)=F_{s,p}(r)$  is a continuous, positive, radially symmetric (i.e., $r=|x|\,t^{-\beta}$), and decreasing function  such that $F(r) \to 0$ with a certain power decay rate that depends on the value of $p$.
\end{theorem}
 Let us point out that $p_c$ is precisely the exponent for which the rates $\alpha$ and $\beta$ diverge, so the representation we offer has to fail below the chosen $p$ range.  We see that all fundamental solutions with $M>0$ are obtained from the one with unit mass, $M=1$, by a simple scaling rule, $F_M(r)=M^{sp\beta}F(M^{2-p)\beta}r)$. The scaling group is a basic tool of our theory, see Subsection \ref{ssec.scaling}. Fundamental solutions with $M<0$ are obtained by just reversing the sign of the solution for mass $|M|$. In the limit $M=0$ the fundamental solution becomes the null function.

 It is interesting to recall what happens in the case of the standard $p$-Laplacian evolution ($s=1$). Then the fundamental solution is the well-known Barenblatt self-similar solution that has the form of Theorem \ref{thm.ssfs} with a profile $F_{1,p}(r)$, given by the explicit  formula
\begin{equation}\label{sol.Bar}
F_{1,p}(r)=\left(C+ kr^{\frac{p}{p-1}}\right)^{\frac{p-1}{2-p}} \qquad \mbox{\rm for} \quad \frac{2N}{N+1}<p<2,
\end{equation}
cf. \cite{VazSmooth}, Formula 11.8. Here, $C>0$ is a free constant and $k=k(s,p,N)>0$ is an explicit constant. On the contrary, in the fractional case $0<s<1$ the work to establish existence of the self-similar profile is highly nontrivial. As a consequence, the result is proved only in Section \ref{sec.fs}.  Important previous steps are the barrier constructions of Sections \ref{sec.barr1} and \ref{sec.barr2},  followed by other tools explained in Sections \ref{sec.masscon} to \ref{sec.diff}. On the other hand, the proof of uniqueness of fundamental solutions is simpler in both cases: for the standard case see \cite{KV88}, for our fractional case see Section \ref{sec.fs} below. The study of the convergence of our solutions in the limit cases of the parameters is continued in Section \ref{sec.limitcases}.

\medskip

\noindent \bf Profile decay rates. \rm A marked novelty with respect to \cite{VPLP2020} happens when describing the actual decay rate of profile $F_{s,p}(r)$ as $r\to\infty$, so-called  \sl tail decay\rm, an important issue in itself. We show that it involves the existence of a new critical exponent $p_1$ that solves the algebraic expression $sp(p-1)=N(2-p)$ and lies in the interval $(p_c,2)$. Solving the quadratic equation in $p$ we have
\begin{equation}\label{eq.p1}
p_1(N,s)=\frac{\sqrt{(\lambda-1)^2+8\lambda}-\lambda+1}{2}
=\frac{4\lambda}{\sqrt{(\lambda-1)^2+8\lambda}+\lambda-1}\,,
\end{equation}
where $\lambda=N/s>1$. It ranges from  $p_1(1)=\sqrt2 $ to $p_1(+\infty)=2$. For comparison, $p_c(\lambda)=2\lambda/(\lambda+1)$ ranges from 1 to 2. The reader may check that $p_c<p_1< 2$.

 \begin{theorem}\label{thm.beh} Let $F$ be the self-similar profile of the fundamental solution. If $p_c<p<p_1$ we have the decay rate
\begin{equation}\label{eq.decay1}
 F(r)\sim C_0(N,s,p) r^{-sp/(2-p)}\,,
\end{equation}
 while for $p_1<p<2$ the decay rate is given by
\begin{equation}\label{eq.decay2}
 F(r)\approx C r^{-(N+sp)}\,.
\end{equation}
 In the critical case $p_1$ both exponents coincide and there appears a logarithmic factor as correction, given in Theorem \ref{thm.barrcomp_1}.
  \end{theorem}

\begin{figure}[h!]
\centerline{\includegraphics[width=4in]{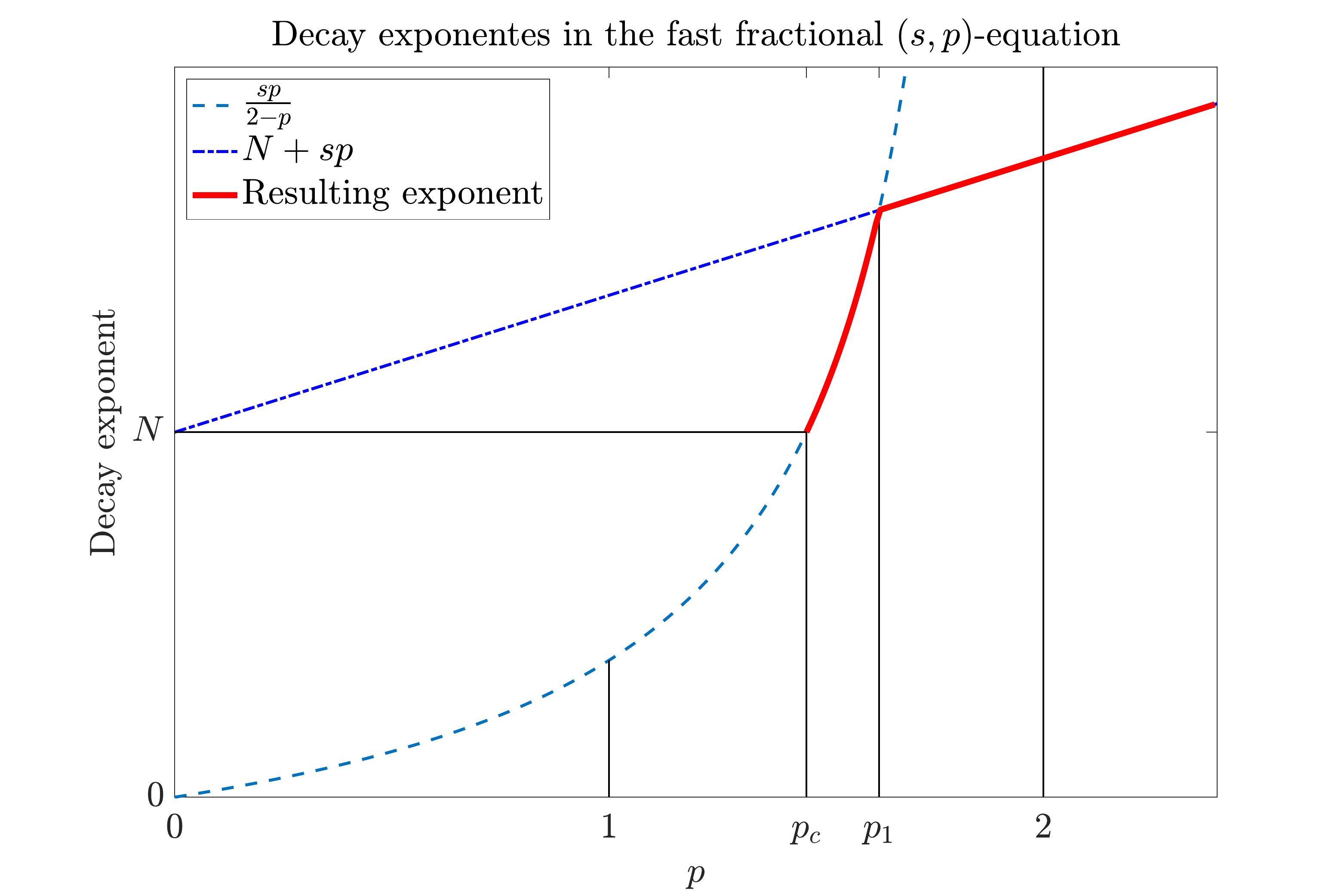}}
\caption{Spatial decay represented against $p$ for $s=0.75$}
\label{fig:Exponents}
\end{figure}

The sign $\approx$ means equivalence up to a constant factor that depends only on the parameters $N,s$ and $p$.
The stricter sign $\sim$  means that the ratio tends to 1.

Thus, there are two possible different decay rates at the tail of the self-similar profile. It is interesting to compare them with the standard $p$-Laplacian evolution equation, where the explicit Barenblatt profiles \eqref{sol.Bar} only exhibit a tail behaviour of type \eqref{eq.decay1} in the whole range $p_c<p<2$. By analogy, we propose to call this behaviour a fast diffusion decay.

The new behaviour of type \eqref{eq.decay2} in Theorem \ref{thm.beh} is due to the influence of the nonlocal (long-range) fractional kernel in Definition \eqref{frplap.op}, hence we propose to call it typical fractional diffusion decay.  We will show after a delicate analysis that this fractional effect becomes dominant in the upper subrange $p_1<p<2$, and it connects with the behaviour for $p\ge 2$ described in paper \cite{VPLP2020}. With respect to the  \eqref{eq.decay1}-type behaviour, we observe a fattening of the tail in the upper subrange $p_1<p<2$.

We ask the reader to check that both decay formulas imply that $F$ is integrable in $\ren$, and also that the integrability is formally lost for the  exponent $p=p_c$, precisely the value that marks the limit of the range where our constructions are possible. The tail behaviour is settled in Sections \eqref{sec.vss1}, \ref{sec.pos.lr}, and \ref{sec.pos.tail.ur}, see in particular formula \eqref{lim.profile.gr} and Corollary \ref{cor.ob.fs}.\nc

\medskip

\noindent \bf Asymptotics. \rm The fundamental solution is the key to the  long-time behaviour of our problem with general initial data, since it represents,  in Barenblatt's words, the \sl intermediate asymptotics, cf. \rm \cite{Barbk96}. Here is the  asymptotic result we obtain.

\begin{theorem}\label{thm.ab1} Let $p_c<p<2.$ Let $u$ be a solution of Problem \eqref{frplap.eq}-\eqref{init.data} with initial data $u_0\in L^1(\ren)$ of integral $M$,  and let $U_M$ be the fundamental solution with the same mass. Then,
\begin{equation}\label{lim.ab.L1}
\lim_{t\to\infty} \|u(t)-U_M(t)\|_1=0\,.
\end{equation}
We also have the sup-norm estimate
\begin{equation}\label{lim.ab.infty}
\lim_{t\to\infty} t^{\alpha}\|u(t)-U_M(t)\|_q=0\,,
\end{equation}
and the $L^q$-estimate
\begin{equation}\label{lim.ab.q}
\lim_{t\to\infty} t^{\alpha(q-1)/q}\|u(t)-U_M(t)\|_q=0\,.
\end{equation}
for all $1<q<\infty$.\rm
\end{theorem}
The theorem is proved in Section \ref{sec.ab}. There is no restriction on the sign of the solution. The basic result is \eqref{lim.ab.L1},  we can easily obtain rates in all $L^q$ spaces, $1<q<\infty$,  by interpolation, see for instance examples in \cite{Vaz17}. The sup-norm estimate is more delicate. Of course, for $M=0$ we just say that  $ \|u(t)\|_1$ goes to zero.

It is interesting to interpret Theorem \ref{thm.ab1} in terms of the rescaled variables defined in Subsection \ref{ssec.ssv}, formula \eqref{eq.rescflow}. Indeed, we may rephrase the result as saying that the rescaled solution $v(y,\tau)$ converges to the equilibrium state $F_M(y)$ of the flow equation \eqref{eq.resc} in all $L^q$-norms, $1\le q\le \infty$. In other words, $F_M$ attracts along this flow all finite-mass solutions with the same mass. The corresponding results for the standard $p$-Laplacian, were proved in \cite{KV88, Vascppme}.

 Numerical computations 
confirm the above theoretical statements, see Figures 2 and 3 at the end of Section \ref{sec.fs}, where an indication of the employed numerical methods is given.\nc

\medskip

\noindent \bf Very singular solutions and nonlinear elliptic problem. \rm
Section \ref{sec.vss1} contains the study of so-called Very Singular Solutions, a kind of singular solutions starting from a very strong singularity at the origin of space that stays for all times. They look similar to the ones of the standard $p$-Laplacian case, namely,
\begin{equation}\label{eq.vss}
U_\infty(x,t)=\,C_\infty  \,(t+T)^{1/(2-p)} { |x|^{-sp/(2-p)} }
\end{equation}
with $T\ge 0$ arbitrary,  $C_\infty$ a universal constant, $C_\infty(N,s,p)$. It is to be noted that in the fractional case these singular solutions exist only in the lower range $p_c<p<p_1$, see Theorem \eqref{thm.vss}. This range reduction needs for an explanation that is again given as the consequence of the long-range effect of the fractional diffusion that turns out to be incompatible with the existence of the VSS.

In Theorem \ref{cor.eigen} we derive as a consequence of the VSS the solution of the singular nonlinear elliptic problem
\begin{equation}\label{eq.eigen1}
\mathcal L_{s,p} F(y)+ \lambda F(y)=0 \qquad \mbox{for } \ y\ne 0, \quad \lambda>0,
\end{equation}
with $F>0$ and $F(0)=+\infty$. It works for $p_c<p<p_1$. This function $F$ has the same spatial form as \eqref{eq.vss},
$F(|x|)=A\,|x|^{-sp/(2-p)}$ with $A$ depending on $\lambda$ and $N,s,p$.

We point out that these results  provide rare examples of explicit solutions for the fractional $p$-Laplacian. Even if they have a simple separate-variables, power-like form, the existence is not obvious, and our proof is based on the previous existence of the fundamental solution with suitable bounds. We still do not any formula to calculate the constant $C_\infty$.

\medskip

\noindent \bf Global Harnack type estimates. \rm
Section \ref{sec.gharn} contains the study of two-sided global estimates for nonnegative solutions with compactly supported and bounded data, where the fundamental solution plays a key role, see in particular Corollary \eqref{two-sided.lr},  Theorem \ref{thm.GH}, and  Corollary \eqref{twosided.ur}. They are known as global Harnack inequalities, though they depend on some information on the initial data. This delicate section culminates the qualitative and quantitative information
about the  equation in the range $p_c<p<2$. The fine results are technically based on the quantitative analysis of positivity of Sections \ref{sec.pos.lr} and \ref{sec.pos.tail.ur}.

\medskip

The next two sections contain additional information.
Section \ref{sec.limitcases} examines the limit cases $s\to 1$ and $p\to 2$. Section \ref{sec.comp.two} is devoted to the comparison of our results with the ones known for the so-called fractional porous medium equation, FMPE: $\partial_tu+(-\Delta)^{s}(u^m)=0$, that was studied in \cite{VazBar2014},

    \medskip

\noindent \bf Very fast diffusion and nonlinear elliptic equations.} Section \ref{sec.vfd} deals with the lowest diffusion range, $1<p\le p_c$. It is not a complete study in any way, since this is a rich territory yet to be explored. We concentrate on the existence of VSS and its consequences, since the results are nontrivial and shed light on the theory developed in this paper for $p>p_c$.  We add a brief study  on explicit solutions for the elliptic equation in the limit case $p=p_c$. The following elliptic classification emerges.

\begin{theorem}\label{thm.ell} Let $F(x)=r^{-sp/(2-p)}$ with $r=|x|$, \ $0<s<1$, and $1<p<2$. Three situations arise:
\noindent (i) for $1<p<p_c$, $F$ is a weak solution of the nonlinear eigenvalue problem
$$
\mathcal L_{s,p} F(x) =c F(x), \quad  \ x\in \ren,
$$
for some $c>0$. It is a classical solution for $r>0$.

\noindent (ii) for $p=p_c$  we have $F=r^{-N}$, which is not integrable at the origin, and we get
$$
\mathcal L_{s,p} F =0, \quad \forall r>0\,.
$$
\noindent  (iii) for $p_c<p<p_1$ there is a constant $c>0$ such that
$$
\mathcal L_{s,p} F +c F=0 \quad \forall r>0\,.
$$
Of course, $c=c(s,p,N)$. Since $p>p_c$ we have $sp/(2-p)>N$, so that $F$ is not integrable at the origin.
\end{theorem}

 A final section contains  a long list of comments and open problems. \nc

\medskip

\noindent {\sc Notations}. We sometimes write a function $u(x,t)$ as $u(t)$ or $u$  when one some of the variables can be safely understood. We use the notation $u_+=\max\{u,0\}$. The letters $\alpha$ and $\beta$ will be fixed at the values given in the self-similar formula \eqref{eq.sse1}. We call universal constant a constant with respect to the variables $(x,t)$ that depends only on $s,p$ and $N$. We also use the symbol $\|u\|_q$ as shortened notation for the norm of $u$ in the $L^q$-space over the corresponding domain when no confusion is to be feared. We denote the duality product in $L^q\times L^{q'}$, with $q$ and $q'$ dual exponents, by $\langle \cdot , \cdot\rangle$. For a function $u(x)\ge 0$ we call mass or total mass the integral $\int_{\ren} u(x)\,dx$, either finite or infinite.
For signed functions that integral does not coincide with the $L^1$ norm, so the use of the term is only justified by analogy.  As already said, the sign $\approx$ means equivalence up to a constant positive factor, while $\sim$ means equivalence with limit 1. We use sometimes the term ``fast diffusion'' for the FPLEE in the range $1<p<2$, and ``very fast diffusion'' for $1<p<p_c$ by analogy with the standard PLE. This is partly justified by the self-similarity analysis since the spread rate of the space profiles is $O(t^\beta)$ and $\beta$ increases as $p$ decreases. If we take as reference linear case $p=2$, where $\beta=1/2s$, then $p<2$ is faster than linear. In view of the results, we can say that $p_c<p<2$ covers the ``good fast diffusion'' range. This name was already coined in the non-fractional setting because of its appealing qualitative and quantitative properties.

%
%
\section{Basic theory}\label{sec.basic}

\smallskip

We establish well-posedness of Problem \eqref{frplap.eq}-\eqref{init.data} in different functional spaces, starting by the consideration of the equation as a gradient flow in $L^2(\ren)$. We obtain unique strong solutions $u\in C([0,\infty]:L^2(\ren))$
 that decay in time as expected. Much of the theory is common to all cases $1< p<\infty$, so we rely on what was said in \cite{VPLP2020} for $p>2$, but some details do differ, like boundedness and continuity. \nc We present an account of  the main qualitative and quantitative properties. Some of the results of the section are new, in particular Subsection \ref{ssec.reg}.

\medskip

\subsection{Existence and uniqueness}

We can solve the evolution problem for equation \eqref{frplap.eq} with initial data  $u_0\in L^2(\ren)$  by using the fact that the equation is the gradient flow of a maximal monotone operator associated to the  convex functional \eqref{Jsp1}, see for instance \cite{MazRT, Puhst, Vaz16}. The domain of that operator is
$$
D_2({\mathcal L}_{s,p})= \{\phi\in L^2(\ren): \ {\mathcal J}_{s,p}(u)<\infty, \ {\mathcal L}_{s,p}u\in L^2(\ren)\}.
$$
Well known theory implies that for every initial $u_0\in L^2(\ren)$ there is a unique strong solution $u\in C([0,\infty): L^2(\ren))$, that we may call the semigroup solution. Strong solution means that $u_t$ and $ {\mathcal L}_{s,p}u\in L^2(\ren)$ for every $t>0$, and the equation is satisfied a.e in $x$ for every $t>0$. The semigroup is denoted as $S_t(u_0)=u(t)$, where $u(t)$ is the solution emanating from $u_0$ at time $0$. Typical a priori estimates for gradient flows follow, cf.  \cite{BrBk73, Ko67}.  The next results are part of the standard theory:
\begin{equation}
\frac12\frac{d}{dt}\|u(t)\|_2^2= -\langle  {\mathcal L}_{s,p}u(t),  u(t)\rangle=-p\,
{\mathcal J}(u(t)),
\end{equation}
where ${\mathcal J}={\mathcal J}_{p,s}$ as in the introduction. Also,
\begin{equation}
\frac{d}{dt}{\mathcal J}(u(t)) =\langle  {\mathcal L}_{s,p}u(t), u_t(t)\rangle=-\|u_t(t)\|^2_2,
\end{equation}
where integrals and norms are taken over $\ren$. It follows that both $\|u(t)\|_2$ and ${\mathcal J}(u(t))$ are decreasing in time, and we get the easy estimate ${\mathcal J}(u(t)\le \|u_0\|_2^2/2pt$ for every $t>0$.
See other properties below.

Moreover, for given $p>1$ (the index of the operator) and every $1\le q\le \infty$, the $L^q$ norm of the solution is non-increasing in time.  We can extend the original set of solutions to form a continuous semigroup of contractions in $L^q(\ren)$ for every $1\le q< \infty$: for every $u_0\in L^q(\ren)$ there is a unique strong solution such that $u\in C([0,\infty): L^q(\ren))$.  The class of solutions can be called the {\sl $L^q$ semigroup} for equation \eqref{frplap.eq} posed in $\ren$. These $q$-semigroups coincide on their common domain. The Maximum Principle applies, and more precisely $T$-contractivity holds in the sense that for two solutions $u_1, u_2$ and any $q\ge 1$ we have
\begin{equation}
\|(u_1(t)-u_2(t))_+\|_q \le \| (u_1(0)-u_2(0))_+\|_q.
\end{equation}
This implies that we have an ordered semigroup for every $q$ and $p$.  An operator with these properties in all $L^q$ spaces is called completely accretive, see \cite{BeCr91}. We can also obtain the solutions by Implicit Time Discretization, cf. the classical references \cite{CL71, Ev78}. The word mild solutions is used in that context, but mild and strong solutions coincide by uniqueness. The operator is also accretive in $L^\infty$ and this allows to generate a semigroup of contractions in $C_0(\ren)$, the set of continuous functions that go to zero at infinity.

Semigroup solutions satisfy the definition of weak solution. An $L^q$ solution $u(x,t)$ is a weak solution if for every smooth test function $\varphi(x,t)$ such that $\varphi(x,t)=0$ for $ t\ge T$ and also for $|x|\ge R$ we have
\begin{equation}\label{weak.sol.p}
\begin{array}{c}
\displaystyle \frac12\iint \frac{\Phi(u(x,t)-u(y,t))\,(\varphi(x,t)-\varphi(y,t))}{|x-y|^{N+sp}}\,dxdydt =\\
\displaystyle \iint u(x,t)\,\partial_t\varphi(x,t)\,dxdt + \int u_0(x)\varphi(x,0)\,dx\,,
\end{array}
\end{equation}
with integration for $x,y\in\ren$, $0\le t\le T$. The set of test functions can be extended by density arguments, thus for an data in $L^2$ with $\mathcal J u\infty$ we may use $\varphi\in W^{s,p}(\ren)\cap L^\infty(\ren)$.

This part of the theory can be done for solutions with two signs, but we will often reduce ourselves in the sequel to nonnegative data and solutions. Splitting the data into positive and negative parts most of the estimates apply to signed solutions. To be precise,
for a signed initial function $u_0$ we may consider its positive part, $u_{0,+}$
and its negative part $u_{0,-}= -u_0+ u_{0,+}=-\max\{-u_0, 0\}$. Then, both $u_{0,+}$ and $u_{0,-}$ are nonnegative and $-u_{0,-}\le u_0\le u_{0,+}$. It follows from the comparison property of the $L^q$ semigroups that
$$
-S_t(u_{0,-})\le S_t(u_0)\le S_t(u_{0,+}).
$$
Therefore, we may reduce many of the estimates to the case of nonnegative solutions.

Finally, let us note that, since the operator acts on a function only by involving differences of the value of the function at two places, all the theory is valid after adding a constant to any admissible function. It follows that the above theory works for functions in the spaces $X_{q,c}=L^q(\ren)+c$ (not a linear space, but an affine space). Moreover, a solution $u$ in the space  $X_{q,c}$ with $c>0$ will become a supersolution for the equation posed in $X_{q,0}=L^q(\ren)$. More precisely,  if $u_{01}(x)\le u_{02}(x)+c$  with $u_{0i}\in L^q(\ren)$, then
 $$
 S_t(u_{01})\le  S_t(u_{02}) +c.
 $$
\subsection{Scaling}\label{ssec.scaling}

In our study we will use the fact that the equation admits a scaling group that conserves the set of solutions. Thus, if $u$ is a weak or strong solution of the equation, then we obtain a two-parameter family of solutions of the same type,
$$
\widehat u(x,t)= A u(Bx,Ct), \quad A, B, C>0,
$$
on the condition that $A^{2-p}C=B^{sp}$. Special choices are  $u_k={  \mathcal T}_k u$, given by
\begin{equation}\label{scal.trn1}
{\mathcal T}_k u(x,t)= k^{N} u(kx, k^{N(p-2)+sp}t),\quad k>0\,.
\end{equation}
This transformation preserves mass.  It can be combined with a second one that keeps invariant the space variable:
\begin{equation}\label{scal.trn2}
{\widehat{\mathcal  T}}_M u(x,t)= M u(x, M^{p-2}t), \quad M>0\,.
\end{equation}
This  changes mass and preserves space. It is often used to reduce the calculations to solutions with unit mass, $M=1$. A third option that we will use is
\begin{equation}\label{scal.trn3}
{\overline{\mathcal  T}}_h u(x,t)= h^{sp/(2-p)} u(h x, t), \quad h>0,
\end{equation}
This one preserves time (i.e., it is of elliptic type). It also changes mass according to the formula
$$
\|{\overline{\mathcal  T}}_h u(t)\|_1= h^{1/(\beta(2-p))}\|u(t)\|_1.
$$

Let us point out that the set of solutions of the equation is invariant under a number of isometric transformations, like: change of sign: $u(x,t)$ into $-u(x,t)$, rotations and translations in the space variable, and translations in time. They will also be used in the sequel.

\subsection{A priori bounds}\label{apriori}

\noindent $\bullet$ In the sublinear range $1<p<2$ \ the operator is homogeneous of degree $d=p-1<1$ in the sense that ${\mathcal L}_{s,p}(\lambda\,u)=\lambda^{p-1} {\mathcal L}_{s,p}u$. Using the general results by  B{\'e}nilan-Crandall \cite{BC81b}  for homogeneous operators in Banach spaces, we can prove the a priori bound
\begin{equation}
(2-p)t u_t< u,
\end{equation}
which holds for all nonnegative solutions, in principle in the sense of distributions. This a priori bound is quite universal, independent of the individual solution. It is based on  the scaling properties and comparison (that hold for all our semigroup solutions).  Therefore, we have almost  monotonicity in time if $u\ge0$. In particular, if a strong solution is positive at a certain point $x_0$ at $t=t_0$, then for all previous times $u(x_0,t)>0$. This is called backward positivity (for nonnegative solutions), a property that has been extensively used in studies of the fast diffusion equation.

\noindent $\bullet$ On the other hand, paper \cite{BC81b} also implies the estimate for all $p\ne 2$ and  Lebesgue exponent $q>1$ we have
\begin{equation}\label{form.brh}
\|u_t\|_q\le \frac 2{(2-p)t}\|u_0\|_q
\end{equation}
for every $1< q\le \infty$. For $q=1$ it is formulated as different quotients.
Formula {\eqref{form.brh} is valid for all signed solutions. For $q=\infty$ we apply the formula for
solutions in $X=L^1(\ren)\cap L^\infty(\ren)$.

\subsection{Energy estimates}

\noindent $\bullet$ {\bf Quadratic estimate.} As we have seen before, for solutions with data in $L^2(\ren)$ and times $0\le t_1<t_2$ we  have the identity
\begin{equation}\label{en.est}
\int_{\ren} u^2(x,t_1)dx-\int_{\ren} u^2(x,t_2)\,dx=2\int_{t_1}^{t_2}
\int_{{\ren}}\int_{{\ren}} |u(x,t)-u(y,t)|^p\,d\mu(x,y)dt\,,
\end{equation}
\noindent where ${d\mu(x,y)}=|x-y|^{-(N+sp)}dxdy $. In the sequel we omit the domain of integration of most space integrals when it is $\ren$ and the time interval when it can easily understood from the context.

We point out that this estimate shows that solutions with $L^2\cap L^p$ data belong automatically to the space $L^p(0,T: W^{s,p}(\ren))$.

\medskip

\noindent $\bullet$ {\bf The $q$-power  estimate.} Arguing in the same way, for solutions with data in $L^q(\ren)$  with $q>1$ and times $0\le t_1<t_2$, we have for {\sl nonnegative solutions}
\begin{equation}\label{en.est.q}
\begin{array}{c}
\displaystyle \int u^q(x,t_1)dx-\int u^q(x,t_2)\,dx=\\ [6pt]
\displaystyle q\iiint |u(x,t)-u(y,t)|^{p-2}\langle (u(x)-u(y)), (u^{q-1}(x,t)-u^{q-1}(y,t))\rangle \,d\mu(x,y)dt,
\end{array}
\end{equation}
with integration in the same sets as before. We will use the inequality
\begin{equation}
(a-b)^{p-1}(a^{q-1}-b^{q-1})\ge C(p,q) \,|\,a^{(p+q-2)/p}-b^{(p+q-2)/p}\,|^p
\end{equation}
which is valid for all $a>b>0$ and $p,q>1$.  This inequality is also true when $b\ge a>0$ by symmetry.
We get the new inequality
\begin{equation}\label{en.est.q}
\begin{array}{c}
\displaystyle C(p,q)\int \iint |u(x,t)^{(p+q-2)/p}-u(y,t)^{(p+q-2)/p}|^{p} \,d\mu(x,y)dt\\ [6pt]
\displaystyle \le \int u^q(x,t_1)dx-\int u^q(x,t_2)\,dx,
\end{array}
\end{equation}
which applies the solutions of the $L^q$ semigroup, $q>1$. This gives a precise estimate of the dissipation of the $L^q$ norm along the flow.

\noindent {\sl Case of signed solutions.} The above results hold after performing some careful adaptations. Thus, the difference we want to control is
\begin{equation}
\displaystyle \int |u|^q(x,t_1)dx-\int |u|^q(x,t_2)\,dx
\end{equation}
and the right-hand side of \eqref{en.est.q} is the same as there if we adopt the notation $u^{q-1}$ to mean the odd power $|u|^{q-2}u$ (a common convention in the PDE literature when it works, it should be used with care). In the next line we use the notation $a^{p-1},b^{p-1},$ to mean $|a|^{p-2}a, |b|^{p-2}b$ and the same happens with the power $(p+q-2)/p$. The equality when
$a>0>b=-b_1$ becomes
$$
(a+b_1)^{p-1}(a^{q-1}+ b_1^{q-1})\ge C(p,q) \,|\,a^{r}+b_1^{r}\,|^p,\qquad r=\frac{p+q-2}{p},
$$
for $a,b_1>0$, which is also true. We thus arrive again at estimate \eqref{en.est.q}
which is also true in this case with the present notations.


\subsection{Difference estimates}\label{ss.diff.est}

It is well known that the semigroup is contractive in all $L^q$ norms, $1\le q\le \infty$. At some moments we would like to know how the norms of the difference of two solutions decrease in time. Such decrease is called {\sl dissipation}. We present here the easiest case, decrease in $L^2$ norm.

\noindent {\bf $L^2$ dissipation.} For solutions with data in $L^2(\ren)$ and times $0\le t_1<t_2$, we  have the identity for the difference of two solutions $u=u_1-u_2$
\begin{equation}\label{diff.est.2}
\begin{array}{c}
\displaystyle\int_{\ren} u^2(x,t_1)\,dx-\int_{\ren} u^2(x,t_2)\,dx=2\int_{t_1}^{t_2}\int_{{\ren}} u u_t\,dxdt\\=
\displaystyle 2\int_{t_1}^{t_2}\int_{{\ren}}\int_{{\ren}}
\displaystyle \left(|u_1(x,t)-u_1(y,t)|^{p-1}(u_1(x,t)-u_1(y,t)) \right. \\[10pt]
 \displaystyle \left.-|u_2(x,t)-u_2(y,t)|^{p-1}(u_2(x,t)-u_2(y,t))\right)\\[8pt]
  \displaystyle (u_1(x,t)-u_2(x,t)-u_1(y,t)+u_2(y,t))\,d\mu(x,y)dt\,,
\end{array}
\end{equation}
\noindent where ${d\mu(x,y)}=|x-y|^{-(N+sp)}dxdy $ as before. Putting $a=u_1(x,t)-u_1(y,t)$ and $b=u_2(x,t)-u_2(y,t) $ and using the numerical inequality as before, we bound below the last integral by
$$
C(p) \iiint  |(u_1(x,t)-u_1(y,t))^{p/2}-(u_2(x,t)-u_2(y,t))^{p/2}|^2\,d\mu(x,y)dt\\.
$$
This is an estimate of the $L^2$ dissipation of the difference $u=u_i-u_2$.

Later on, we  will need the expression of the $L^1$ dissipation in the study of the asymptotic behaviour, but we will postpone it until conservation of mass is proved.

\subsection{Regularity of bounded solutions}\label{ssec.reg}

 We prove the following regularity result for bounded and integrable semigroup solutions.

\begin{theorem}\label{them.reg.new}  Let $u(x,t)$  be the solution of equation \eqref{frplap.eq} with initial data  $u_0\in L^\infty(\ren)\cap L^1(\ren)$.  Then $u$ is uniformly bounded in $\ren\times \re_+$. It is also H\"older continuous in space and Lipschitz continuous in time uniformly for all $t\ge \tau>0$. Precise estimates are given in the proof.
\end{theorem}

\noindent {\sl Proof.} Lipschitz continuity in time for $u(x,t) $ holds for all $t>0$. It comes from the derivative bound $\|u_t\|_\infty\le C\|u_0\|_\infty/t.$  Therefore, for all $t\ge \tau>0$  we  also know that
$\mathcal L_{s,p} u(t)=-u_t$ holds for every fixed time and is uniformly bounded.

The remaining part of the proof consists of showing the local Hölder continuity of the solution with respect to the space variable. This  follows from the elliptic study performed in Section 5 of paper \cite{IannMS2016}, that is valid for any $p>1$. More precisely, we copy below their Corollary 5.5 which reads.

\noindent {\bf Corollary\cite{IannMS2016}}\label{corlocalpha}
{\sl There exist universal constants $C>0$ and $\alpha\in (0,1)$ with the following property: for all $u\in \widetilde{W}^{s,p}(B_{2R_0}(x_0))\cap L^\infty(B_{2R_0}(x_0))$
such that $|{\mathcal L}_{s,p}u|\le K$ weakly in $B_{2R_0}(x_0)$,
\begin{equation}
\label{tesicorloc}
[u]_{C^\alpha(B_{R_0}(x_0))}\le C\big[(K R_0^{ps})^\frac{1}{p-1}+Q(u; x_0, 2R_0)\big]R_0^{-\alpha}.
\end{equation}
}

Here, $x_0\in\ren$ and $R_0>0$. The behaviour for large $x$ is included in the tail estimate:
$$
\mbox{Tail}(u;x,R) =\left( R^{ ps}\int_{B^c_R(x)} \frac{|u(y)|^{ p-1}} {|x - y|^{ N+ps}}
\, dy
\right)^{1/(p-1)}.
$$
It is clear that boundedness of $u$ implies a finite tail. $Q$ is defined as follows:
$$
Q(u;x_0 ,R) = \|u\|_{ L^\infty (B_{ R (x_0 )})} + \mbox{Tail}(u;x_0 ,R).
$$
We conclude that the semigroup solution of the theorem $u(x,t)$ is Hölder continuous with respect to the space variable, uniformly for all $t\ge \tau>0$. \qed

\medskip

\noindent {\bf Remark.}  In case $u_0$ is continuous, i.e., belonging to the space $C_0 (\ren)$ endowed with the sup norm,  we recall that for all exponents $p > 1$,  the solution semigroup is contractive in that space, so solutions with initial data in that class will be continuous in space for every positive time as a consequence of the inequality
$$
\|u(x,t)-u(x+h,t)|_{L^\infty(\ren)}\le \|u_0(x)-u_0(x+h)\|_{L^\infty(\ren)}\le \omega(h)
$$
for a certain modulus of continuity $\omega$.

If $u_0$ is also integrable, continuity of $u$ in time for $t=0$ comes from the uniform continuity in $x$ and  the continuous $L^1$ dependence in time (via a simple triangular argument). If the condition $u_0\in L^1(\ren)$ is not assumed, then we may use   a barrier argument (using vertical displacement in time).  We will not enter into the details of this last comment that is not needed in the sequel.\nc

\subsection{Positivity of nonnegative solutions}\label{ssec.pos}

We  get the following backward positivity result for nonnegative solutions of equation \eqref{frplap.eq}: if at a point $P=(x_0,t_0)\in\ren\times \re_+$ we know that $u(x_0,t_0)$ is strictly positive, then the weighted monotonicity in $t$ of the function $u(x,t)t^{1/(2-p)}$ guarantees that $u(x_0,t)>0 $ for all $0<t<t_0$. Positivity onwards in time will be proved  to be true for $p>p_c$.  It is not true in general for $p<p_c$ because of the phenomenon of extinction in finite time that we will see in Section \ref{sec.vfd}.

Finally, this is a classical argument for positivity that applies to all classical solutions (in the sense that $u_t$ and $ {\mathcal L}_{s,p}u$ exist and are continuous everywhere). At every point $(x_0,t_0)$ where  a  solution reaches the minimum value $u=0$ and
$ ({\mathcal L}_{s,p}u)(x_0,t_0)$ exists, then it must be strictly negative according to the definition formula for the operator. On the other hand, if $u_t$ exists it must to zero. From this contradiction we conclude that a.e. $u(x,t)$ must be positive. By the already proved conservation of positivity, for any $0<t<t_0$ we have  $u(x,t_0)>0$ for a.\,e. $x\in \ren$.  For global positivity in the good fast range see Section \ref{sec.globalpos}. Quantitative positivity statements are contained in Sections \ref{sec.pos.tail.ur} and \ref{sec.gharn}.

\subsection{Comparison via symmetries. Almost radiality}

The Aleksandrov symmetry principle \cite{Alek} has found wide application in elliptic and parabolic linear and nonlinear problems. An explanation of its use for the Porous Medium Equation is given in  \cite{Vazpme07}, pages 209--211, where previous references are mentioned. In the parabolic case it says that whenever an initial datum can be compared with its reflection with respect to a space hyperplane, say $\Pi$, so that they are ordered, and the equation is invariant under symmetries, then the same space comparison applies to the solution at any positive time $t>0$.

The result has been applied to elliptic and parabolic equations of Porous Medium Type involving the fractional Laplacian in \cite{VazBar2014}, section 15. The argument of that reference can be applied in the present setting. We leave the verification to the reader. The standard consequence we want to derive is the following

\begin{proposition}\label{Alek} Solutions of our Cauchy Problem having compactly supported data in a ball $B_R(0)$ are radially decreasing  in space  along an outgoing a cone of directions $K_\theta(x)$ for all $|x|\ge 2R$ and some $\theta=\theta(x)$. Moreover, whenever $|x|>2R$ and $|x'|>|x|+2R$, then we have $u(x,t)\ge u(x',t)$
for all $t>0$.\end{proposition}

Here, $K_\theta(x)$ is the cone with vertex $x$, axis directed along the line $\overline{0x}$ and aperture angle $\theta$.

\subsection{Boundedness for positive times}\label{sec.smooth1}

From this moment on and but for the last section of the paper, we work in the range $p_c<p<2$, sometimes called the \sl good fast range\rm.
An important result valid for many nonlinear diffusion problems with homogeneous operators is the so-called $L^1$-$L^\infty$ smoothing effect. In the present case we have

\begin{theorem}[Smoothing effect] \label{L1-Linfty} Let $p>p_c$. For every solution with initial data $u_0\in L^1(\ren)$ we have
 \begin{equation}\label{smoot.effct}
 |u(x,t)|\le C_b(N,p,s) \|u_0\|_1^{\gamma}\,t^{-\alpha}\,,
 \end{equation}
 with exponents  $\alpha=N\beta$, $\gamma=sp\beta$ and $\beta=1/(N(p-2)+sp)$.
  \end{theorem}

This is also true for $p\ge 2$ as proved and used in \cite{VPLP2020}. The exponents are  given by the scaling rules (dimensional analysis). The result has been recently proved by Bonforte and Salort \cite{BS2020}, Theorem 5.3, where an explicit estimate for the constant $C_b(N,p,s)$ is given.
Note that this formula has to be invariant under the scaling transformations of Subsection \ref{ssec.scaling}.

 We may now combine Theorems \ref{them.reg.new} and \ref{L1-Linfty}  to conclude that $L^1$ solutions become H\"older continuous for all positive times with a given H\"older exponent that only depends on $p,s$ and $N$.

\subsection{On the fundamental solutions}\label{ssec.bddness}

The existence and properties of the fundamental solution of Problem \eqref{frplap.eq}-\eqref{init.data} are a main concern of this paper. We expect the FS  to be unique, positive and self-similar for any given mass $M>0$. Self-similar solutions have the form
\begin{equation*}\label{eq.ssf1b}
U(x,t;M)=t^{-\alpha}F(x\,t^{-\beta};M)\,.
\end{equation*}
(more precisely, this is called direct self-similarity). Substituting this formula into equation  \eqref{frplap.eq}, we see that time is eliminated as a factor in the resulting transformed equation on the condition that: $\alpha+1=(p-1)\alpha+\beta sp$. We also want integrable solutions that will enjoy the mass conservation property, which implies $\alpha=N\beta$. Imposing both conditions, we get
\begin{equation*}\label{eq.sse}
\alpha=\frac{N}{sp-N(2-p)}, \quad \beta=\frac{\alpha}{N}=\frac{1}{sp-N(2-p)}\,,
\end{equation*}
as announced in the Introduction. It is precisely the requirement of positivity of $\alpha$ and $\beta$ what forces the choice $p>p_c$ in the study of self-similarity.

The profile function $F(y;M)$ must satisfy the nonlinear stationary fractional equation
\begin{equation}\label{eq.ssp}
{\mathcal L}_{s,p} F= \beta \,\nabla\cdot (yF)\,.
\end{equation}
See a similar computation for the Porous Medium Equation in \cite{Vazpme07}, page 63. Using rescaling ${\widehat {\mathcal  T}}_M $, we can reduce the calculation of the profile to mass 1 by the formula
$$
F(y;M)=M^{sp\beta}F(M^{-(p-2)\beta}y;1).
$$
In view of past experience with $p\ge 2$, we will look for  $F$ to be radially symmetric, monotone in $r=|y|$, and positive everywhere with a certain behaviour as $|y|\to \infty$.

We have proved that all solutions with $L^1$ data at one time will be uniformly bounded } $p>p_c$. Thus,  $F$ must  be bounded. Moreover, bounded solutions have a bounded $u_t$ for all later times. By Corollary 5.5 of \cite{IannMS2016} $F$ must be H\"older continuous. In the case of the fundamental solution, this means that $r^{1-N}(r^NF(r))'$ is bounded, hence by \eqref{eq.ssp} $rF'$ is bounded,  and $F$ is regular for all $r>0$.
By monotonicity the limit $F(0+)=\lim_{r\to 0}F(r)$ exists and is finite, therefore $F(x)$ is a continuous function in $\ren$.

The self-similar fundamental solution must take a Dirac mass as initial data, at least in the sense of initial trace, i.e., $u(x,t)\to M\delta(x)$ as $t\to 0$ in a weak sense. It  will be invariant under the scaling group  ${{\mathcal  T}}_k $ of Subsection \ref{ssec.scaling}. All of this will be proved in this paper. The detailed statement is contained in Theorems \ref{thm.exfs} and \ref{thm.exfs2}, and whole proofs follow there.


\subsection{Self-similar variables}\label{ssec.ssv}
Often in the sequel, it will be convenient to pass to self-similar variables. This is done by zooming the original solution according to the self-similar exponents \eqref{eq.sse}.
More precisely, the change uses by the formulas
\begin{equation}\label{eq.rescflow}
u(x,t)=(t+a)^{-N\beta}v(y,\tau) \quad y=x\,(t+a)^{-\beta}, \quad \tau=\log(t+a)\\,
\end{equation}
with $\beta=(N(p-2)+sp)^{-1}$, and any $a>0$, we mostly use $a=1$. Here, $\tau$ is called the {\sl new time} or logarithmic time. The formulas imply that  $v(y,\tau)$ is a solution of the corresponding PDE:
\begin{equation}\label{eq.resc}
\partial_\tau v + {\mathcal L}_{s,p} v -\beta \nabla\cdot(y\, v)=0\,.
\end{equation}
This transformation is usually called \sl continuous-in-time rescaling \rm to mark the difference with the transformation with fixed parameter \eqref{scal.trn1}.

Note that the rescaled equation \eqref{eq.resc} does not change with the time-shift $a$, but the initial value in the new time does,
$\tau_0=\log(a)$. This has to be taken into account since it may cause confusion. Thus, when  $a=0$ we get $\tau_0=-\infty$ and the $v$ equation is defined for $\tau\in \re$. This choice is good for the self-similar solutions. In any case, the mass of the $v$ solution at new time $\tau\ge \tau_0$ equals that of the $u$ at the corresponding time $t\ge 0$.

Sometimes $\tau$ is defined as $\tau=\log((t+a)/a)$ without change in the equation. It is just a displacement in the new time, but it is important to take it into account in the computations.

\noindent {\bf Denomination.} For convenience we sometimes refer in the sequel to the solutions of the rescaled equation \eqref{eq.resc} as $v$-solutions, while the original ones are $u$-solutions.

\medskip

%
%

\section{Barrier construction and tail behaviour I}\label{sec.barr1}

 A main step in the paper is to construct an upper barrier $\widehat u$ for the solutions  of the  Cauchy problem with suitable data. The barrier will be needed in the proof of existence of the fundamental solution, obtained as limit of approximations having the same mass as the initial Dirac delta.  We will only need to consider nonnegative data and solutions. Besides, we may perform the construction using bounded radial functions with compact support as initial data and then use some  comparison argument  to  eliminate the restrictions of radial symmetry and compact support. The barrier we seek will be radially symmetric, decreasing in $|x|$ and will have  behaviour $\widehat u(x,t)=O(|x|^{-\gamma})$ for very large $|x|$, with some $\gamma>N$ to be integrable at infinity.

\subsection{\bf Barrier in the lower range} \label{ssec.barr.low}

 At this moment we find the first evidence of the two $p$ subranges with different behaviour.  We begin with a case that is quite new with respect to the analysis done in \cite{VPLP2020} for $p=2$. It is rather inspired in the analysis of the fractional porous medium done in \cite{VazBar2014}. We will choose a very specific candidate to be the upper barrier.
\begin{equation}\label{barrier2}
{\widehat v}(y)=C_1 \, r^{-sp/(2-p)}, \quad r=|y|\,,
\end{equation}
where use the rescaled $v$-variable.   We notice that this barrier  is stationary in time, $\widehat v(r)$, a very good property when we will come to large time asymptotics.   Note that ${\widehat v}(y)$ is integrable at infinity if $sp/(2-p)<N$, i.e., if $p>p_c $. The value of $p_c$ was already mentioned in the Introduction, formula \eqref{pc.range}. We will see that the barrier will not work in the whole range $p_c<p<2$, only in the lower part of that range.

\begin{proposition}\label{prop.barr1} Assume that $p_c<p<p_1$. Then, there is a constant $C_1^* =C_1^*(p,s,N) >0$ such that for $C_1\ge C_1^*$ the (stationary) supersolution condition
\begin{equation}\label{ineq.super}
E(\widehat v):= {\mathcal L}_{s,p} {\widehat v}-\beta \, r^{1-N}(r^N{\widehat v})_r\ge 0
\end{equation}
holds everywhere in $D=\ren\setminus\{0\}$, i.\,e., for all $r>0$.  The value of $p_1$ is the solution of the algebraic equation $N+sp=sp/(2-p)$ in the interval in $1<p<2$ announced in \eqref{eq.p1}, so that $p_c<p_1<2$.
\end{proposition}

\noindent{\sl Proof.}  We analyze \eqref{ineq.super}. To begin with, the term
\begin{equation}
-\beta \, r^{1-N}(r^N {\widehat v})_r=\frac1{2-p}{\widehat v}= \frac{C_1}{2-p}r^{- sp/(2-p)}>0
\end{equation}
 has a good sign. On the other hand, ${\mathcal L}_{s,p} \widehat v$  may be negative. But we do not care
as long as it is finite by virtue of the following argument.

\noindent (ii) First, we need to make sure that the singularity of ${\widehat v}$ at $r=0$ does not disturb the definition of ${\mathcal L}_{s,p} {\widehat v}(r)$ as a finite function for  $r>0$, i.\,e., we need  \  ${\widehat v}^{p-1}\in L^1(B_1(0))$. This happens if
\begin{equation}
 \frac{ps(p-1)}{2-p}<N, \quad \mbox{i.\,e. when \ } \ h(p):= \frac{p(p-1)}{2-p}<\frac{N}{s}.
\end{equation}
 Since $h(1)=0$ and $h(2)=\infty$, there is a solution $p_1$ of the equation $h(p_1)=N/s$
 for $1<p_1<2$. Moreover, for $p=p_c$ we have
$$
 \frac{p(p-1)}{2-p}=\frac{2N(N-s)(N+s)}{2s\,(N+s)^2 }=\frac{N(N-s)}{s\,(N+s)}< \frac{N}{s},
 $$
hence the solution $p_1$ lies between $p_c$ and $1$. Such a solution $p_1=p_1(N,s)$ is unique because $h$ is monotone in the interval $(1,2)$.
We can also check that the equation expressing equality of the decay exponents in Theorem \ref{thm.beh},
$$
N+sp= \frac{sp}{2-p},
$$
has exactly the same solution $p=p_1$.

\noindent (iii) Admitting the exponent restriction $p\in (p_c,p_1)$, we may continue. We will calculate ${\mathcal L}_{s,p} {\widehat v}$ at $r_0=1$. If $y_0$ is any point with $|y_0|=1$ we have
$$
{\mathcal L}_{s,p} {\widehat v}(y_0)= C_1^{p-1}\int_{\ren}   \frac{\widehat v(y_0)^{p-1}-\widehat v(y)^{p-1}}{|y_0-y|^{N+sp}} \,dy\,.
$$
The part of this integral performed in the exterior domain $\{r>1\}$ is positive, so this integral, or a part of it, can be disregarded for the purpose of proving inequality \eqref{ineq.super}. Besides, the integral is convergent for $y$ near $0$,  since ${\widehat v}^{p-1}\in L^1(B_1(0))$ by the previous point and the denominator approaches 1 in that region. Moreover, the integral is negative but finite in the annulus $\ve<r<1$, but for a neighbourhood of the point $y=y_0$ (located at the border surface $r=1$). Therefore, we still have to examine the integral in a small ball $B_\rho(y_0)$ centered at $y_0$. Since $\widehat  v$ is a $C^2$ function without critical points in $B_0$  the integral converges by the results of \cite{KKL}, Section 3. This  takes into account the cancelations of differences at points located symmetrically w.r.to $y_0$.  Another proof is given at \cite{DTGCV}, Lemma A.2. The integrability rate near $y=y_0$ only depends on the $C^2$ norm of the function at the point and a nonzero lower bound for $|\widehat v(r)'$|  in the neighbourhood of $r=1$. \rm

Summing up, there is a finite constant $k(N,s,p)$ such that
\begin{equation}\label{def.k}
{\mathcal L}_{s,p} {\widehat v}(1)=-k(N,s,p)\,C_1^{p-1}.
\end{equation}
The sign of $k$ will be important. We insert the minus sign because this is what happens when $s=1$ and the calculation is explicit. See more below.

\noindent (iv) In order to calculate the expression $E(\widehat v)$ ar $r\ne 1$ we use a scaling transformation that leaves the expression invariant. The transformation is defined for $h>0$ by
$$
v_h(y,\tau)={\overline{\mathcal  T}}_h(v)(y,\tau):=h^{\gamma} v(h y, \tau).
$$
Working out the details, we see that in order to leave equation \eqref{eq.resc} invariant, the correct value of the scaling exponent is $\gamma=sp/(2-p)$, and then we get
$$
E(v_h)(y,\tau)=h^{\gamma}E(v)(hy,\tau).
$$
Next, we point out that $\overline{\mathcal  T}_h$ leaves our choice function $\widehat v$ invariant, ${\overline{\mathcal  T}}_h(\widehat v)=\widehat v$. Applying these facts to a point with $|y|=r$ and choosing $h=1/r$, we get
$$
{\mathcal L}_{s,p} {\widehat v}(r)=r^{-\gamma}{\mathcal L}_{s,p} {\widehat v}(1)=-kC_1^{p-1}\,r^{-\gamma},
$$
and finally,
\begin{equation}\label{ineq.super.2}
E(\widehat v)(r)=r^{-\gamma}E(\widehat v)(1)= \left(\frac{C_1}{2-p}-kC_1^{p-1}\right)\,r^{-\gamma}, \qquad
\gamma=sp/(2-p).
\end{equation}

(v) The end of proof consists of analyzing that inequality.  In the case where $k>0$ there exists an optimal constant $C_1^*$ such that equality holds:
$$
E(C_1^* r^{-sp/(2-p)})=0, \qquad C_1^*=((2-p)k(N,s,p))^{1/(2-p)}\,,
$$
while $E(C_1 r^{-sp/(2-p)})>0$  if $C_1>C_1^*$ and $r>0$. We will confirm that possibility in Section \eqref{sec.vss1},  see Theorem \ref{cor.eigen}. We will also have in that case $E(C_1 r^{-sp/(2-p)})<0$  if $0<C_1<C_1^*$ and $r>0$.

In case $k\le 0$ the expression $E(\widehat v)(r)$ is positive for every $C_1>0$ and the supersolution has been obtained in all cases $C_1>0.$  We will later rule out this possibility.  \qed

\medskip

\noindent {\bf Remark.} Note that the mass of this special solution is infinite. However, the decay $\widehat u(x,t)=O(|x|^{-sp/(2-p)})$ guarantees that $\widehat u(\cdot,t)$ is integrable at infinity precisely for $p>p_c$.

\medskip

We can use that result to prove the comparison we are aiming at.

\begin{theorem}\label{thm.barr.lr} Let $p_c<p<p_1$. Let $v$ be a solution of the renormalized equation with $v_0\in L^1(\ren)$ and $v_0(y)\le G(y;C_1):=C_1 \,|y|^{-sp/(2-p)} $ for $|y|>0$, where $C_1>0$ is large enough as in the previous construction, $C_1\ge C_1^*$. Then, for every $\tau>0$ we have
\begin{equation}
v(y,\tau)\le G(y;C_1).
\end{equation}
In terms of $u$, this means that if $u_0(x)\le C_1|x|^{-sp/(2-p)} $, then
\begin{equation}\label{decay.u1}
u(x,t)\le \widehat u(x,t+1):=C_1\,\,\frac{(t+1)^{1/(2-p)}} { |x|^{sp/(2-p)} }\,.
\end{equation}%
    \end{theorem}

The proof of the comparison theorem follows the lines of Theorem 3.2, (iv) of \cite{VPLP2020}.  By slightly changing the proof,  we may assume that $u_0(x)\le A \,|x|^{-sp/(2-p)} $ for $|x|\ge 0$, with $A>C_1^*$ and conclude that
\begin{equation}\label{decay.u1.b}
u(x,t)\le C_1^*\,\frac{\,(t+T)^{1/(2-p)}} { |x|^{sp/(2-p)} }\,,
\end{equation}%
if $T\ge (A/C_1^*)^{2-p}$.

\medskip

\noindent  {\bf Tail behaviour. } The result implies that in this $p$-range the spatial decay of the  class of solutions under consideration is at least $O(|x|^{-sp/(2-p)} )$ for every fixed time. We will prove below that such a rate is exact for every data in the subclass of nonnegative and nontrivial initial data.

On the other hand, the reader will see that  expression \eqref{decay.u1} decays in $|x|$ but increases in time. This, in principle  surprising fact, already happens for the standard fast $p$-Laplacian diffusion ($s=1$) as  the Barenblatt solutions show. Using \eqref{eq.ssf1} and \eqref{sol.Bar}, we find an expression of the form $O(|x|^{p/(2-p)}t^{1/(2-p)})$, in agreement with \eqref{decay.u1} with $s=1$. The upper barrier \eqref{decay.u1.b} decays with $|x|$ in the far field limit but it builds up with time. This is not a defect of the barrier but a very specific property of fast diffusion flows, already described in the PME case in \cite{Vascppme}.

\medskip

\noindent {\bf Other remarks.}  1) In view of this theorem we say that when $C_1>0$  is  large enough, then $G(y;C_1)$ is a {\sl global supersolution} for the $v$ equation, and    $\widehat u(x,t;C_1)$ is {\sl global supersolution} for the $u$ equation.  Note that $G(y;C_1)$ is not integrable at the origin since $sp/(2-p)>N$ for $p>p_c$. This implies that some extra care must be taken in the above computations, but in the end it causes no problem.

\noindent 2)  The global supersolution constructed in this section will appear again, subject to closer scrutiny in the study of the Very Singular Solutions in Section \ref{sec.vss1}.  \nc

%
%
\section{Barrier construction and tail behaviour II}\label{sec.barr2}
We go on to consider the remaining exponent interval, $p_1\le p<2$. We will divide the study into the open upper interval  $p_1< p<2$ as the main part, plus the new critical case $p=p_1$ as a border case.

\subsection{\bf Barrier construction for the upper range}\label{ssec.barr1}

 Here we want to copy the method and result that we have successfully used for $p>2$ in \cite{VPLP2020}. Again, we will use the rescaled solution  and Equation \eqref{eq.resc} introduced in Subsection \ref{ssec.ssv}. Translating previous a priori bounds for the original equation into the present rescaled version, we see that all the rescaled solutions are bounded $v(r,t) \le A((t+a)/t)^{N\beta}$. If $u_0$ is bounded above by say $L>0$, then we get an $L^\infty$ bound for $v$ of the form $v(y,\tau)\le L(t+a)^{N\beta}$. We also get a bound of the form $v(r,t)\le B r^{-N}$, as a consequence of finite mass  $M>0$, radially symmetry and monotonicity in $|x|$; this decay at infinity  is uniform in time.  $A$ and $B$ depend on the mass of the solution, $A= cM^{2p\beta},$ $B=cM$. In this section  $c$ denotes a positive constant that can be computed as a function of $N$, $s$ and $p$. It can vary from line to line.

 The upper barrier  we consider in self-similar variables   will be stationary in time, $\widehat v(r)$.  We again use the notation $r=|y|>0$ when we work with self-similar variables. The construction of the supersolution takes a different form from the previous section. We forget the origin and concentrate on finding a new upper estimate for large $r$ so that we get an integrable barrier in a region $r\ge R_1\ge 0$.
  The barrier will have the form of an inverse power in the far field region. We want to be able  to compare a given solution $v(r,t)$ with $\widehat v(r)$ in an outer domain. Therefore, we have to restrict the solution concept into the solution of a Dirichlet problem in a time interval.

 To be precise, the barrier will be defined by different expressions in two regions: For $|y|\le 1$ we take
${\widehat v}(y;C_1)=C_1,$ \ while for  $r=|y|>1$ \ it has the form ${\widehat v}(y;C_1)=C_1 r^{-(N+\gamma)},$ hence
\begin{equation}\label{barrier1}
\widehat v(y;C_1)=C_1\,\min\{r^{-(N+\gamma)},1\},
\end{equation}
with a suitable $\gamma>0$. We will make the choice $\gamma=sp$ which will turn out to be the best choice. Note the difference with the exponent choice in the previous section.
A key step in proving that this choice produces a barrier is the following supersolution result in an exterior domain.

\begin{proposition} If $p\in (p_1,2)$ and $C_1$ is large enough, the function \ ${\widehat v}(y)$ defined by
\eqref{barrier1}   with $\gamma=sp$ satisfies the supersolution inequality
\begin{equation}\label{ineq.super1}
E(\widehat v):={\mathcal L}_{s,p} {\widehat v}-\beta \, r^{1-N}(r^N{\widehat v})_r\ge 0
\end{equation}
for every $r>2$.
\end{proposition}

\noindent {\sl Proof.} We see that
\begin{equation}
-\beta \, r^{1-N}(r^N{\widehat v} )_r=\beta sp \, C_1r^{-N-sp}>0  \,,
\end{equation}
which has a good sign. On the other hand, ${\mathcal L}_{s,p} \widehat v(r)$  may be negative. We have to  estimate the contribution of different regions against the previous bound. In the sequel we fix a point $y=y_0$ with $r=r_0>2$ and operate in different subregions to evaluate ${\mathcal L}_{s,p} \widehat  v(r_0).$

(i) A first term comes from the influence of the inner core $B_1(0)=\{r: r<1\}$ where we have  ${\widehat v} \approx C_1$. When we evaluate the contribution from this region to the integral ${\mathcal L}_{s,p} \widehat v$ on the exterior region of concern for us, $r\ge 2$, we get the quantity
$$
I_1:= \left.{\mathcal L}_{s,p} \widehat  v\,\right]_1\sim -c C_1^{p-1}\, r_0^{-N-sp}\,.
$$
This implies a first condition that we must impose: $cC_1^{p-1}r^{-N-sp}\le \ve C_1r^{-N-\gamma}$. It holds if  $c \le \ve C_1^{2-p}$ with $\ve$ small, say $\ve<1/4$.

(ii)  The contribution of the region $\{r>r_0\}$ need not be counted (either the whole region or a part of it, see below)  since the integrand is positive, as seen in the formula of the operator. So we have to calculate the contribution from the annulus $\{1<r<r_0\}$. We split it into several pieces.

(iii) For $D_2=\{1<r<r_0/2\}$ we get at $r_0$ the negative contribution  $I_2:=
\left.{\mathcal L}_{s,p} \widehat  v(r_0)\,\right]_{D_2 }$, i.e., with integral extended to $D_2$. We have:
$$
|I_2| \le c
\int_{1}^{r_0/2}  C_1^{p-1} r^{-(N+sp)(p-1)}r_0^{-(N+sp)}r^{N-1}dr=
$$
$$
=c  C_1^{p-1} r_0^{-(N+ps)}\int_{1}^{r_0/2}  r^{N(2-p) -sp(p-1)}r^{-1}dr
$$
Since $N(2-p)<sp(p-1)$ precisely for $p_1<p<2$, the integral converges as $r\to\infty, $ and  we can estimate this contribution as
$$
|I_2|\le c C_1^{p-1} r_0^{-(N+ps)}.
$$
Hence, we need a second condition: $c\le \ve C_1^{2-p}$. The coincidence of the limit of the admissible $p$-range with the end of the $p$-range obtained in the previous section is one of the lucky moments of the present paper.

(iv) A similar but simpler calculation can be done for the contribution of the region
$D_2'=\{y: r_0/2\le |y|\le r_0, \ |y-y_0|\ge r_0/2\}$.

(v)  We are left with the part of the ball $B_{r_0/2}(y_0)$ contained in $B_{r_0}(0)$. However, it is easier to add the rest of the ball and consider the region $D_3=B_0= B_{r_0/2}(y_0)$ (we recall that the exterior of $B_{r_0}(0)$ can be counted or not at will). We have
$$
I_3:=\left.{\mathcal L}_{s,p} \widehat  v(r_0)\,\right]_3 =
\int_{B_{r_0/2}(0)}\frac{\Phi_p\left(\widehat v(y_0)-\widehat v(y_0+z)\right)}{|z|^{N+sp}}\,dz\,,
$$
where $z=y-y_0$, $|y_0|=r_0$. Since $\widehat  v$ is a $C^2$ function without critical points in $B_0$ the integral converges by the results of \cite{KKL}, Section 3, or  \cite{DTGCV}, Lemma A.2, even if there is a singularity of the integrand at $r=r_0$. We must evaluate how this integral   changes with variable $r_0$. This is easily done by rescaling, taking two choices, $r_0''= \lambda r_0'$, and computing the change of the integral $I_3$ when passing from $r_0'$ to $r_0''$, that happens to be
a factor of $\lambda$.  It easily follows that for $r_0''>r_0'\ge 4$:
$$
I_3(r_0'')= \lambda^{-\gamma'}I_3(r_0')\,  \qquad   \mbox{with \ } \gamma'=(N+sp)(p-1)+sp.
$$
Therefore, after putting $r_0=4$ and $r_0''=r_0$, we get for large $r_0$: $|I_3(r_0)|= cC_1^{p-1}\, r_0^{\gamma'}$. We now check that
$$
\gamma'-(N+sp)= (N+sp)(p-2)+sp=sp(p-1)+N(p-2)
$$
is a positive  quantity  precisely for $p>p_1$. \nc Therefore, the calculation enters our scheme if we impose the condition\
$$
c\le \ve C_1^{2-p}r_0^{\gamma'}.
$$
This last condition is weaker than the other two. The result follows. \qed

We need a technical quantitative lemma.

\begin{lemma}\label{thm.Linf}
Let $v$ be a solution of the renormalized equation with nonnegative data $u_0=v_0$ such that
$\|v_0\|_\infty\le L$ and $\|v_0\|_1\le M$. There exists a constant $K(L;M)>0$ such that
\begin{equation}
\|v(y,\tau)\|_\infty\le K  \mbox{ for } \ y \in\ren,  \tau>0\,.
\end{equation}
\end{lemma}

\noindent {\sl Proof.}  Let us pick some $\tau_0>0$. Starting from initial mass $M>0$, from the smoothing effect \eqref{L1-Linfty} and the scaling transformation \eqref{eq.rescflow} (we put $\tau=\log (t+1)$), we know that
$$
v(y,\tau)= (t+1)^{\alpha}u(x,t)\le C_0 M^{sp\beta}((t+1)/t)^{N\beta}= C_0 M^{sp\beta}(1-e^{-\tau})^{-N\beta},
$$
where $C_0$ is universal. We have $\|v(\tau)\|_\infty\le K$ for all $\tau\ge \tau_0$ if
\begin{equation}\label{barr.cond1}
C_0 M^{sp\beta}\le K (1-e^{-\tau_0})^{N\beta}.
\end{equation}
On the other hand, for $0\leq\tau<\tau_0 $ we argue as follows: from $v_0(y)\le L$  we get $u_0(x)\le L$, so $u(x,t)\le L$, therefore
$$
\|v(\tau)\|_\infty \le L (t+1)^{\alpha}= L \,e^{\alpha\tau}.
$$
We now impose
\begin{equation}\label{barr.cond3}
 K=\min\{L \,e^{\alpha\tau_0}, C_0 M^{sp\beta}(1-e^{-\tau_0})^{-N\beta}\},
\end{equation}
we get $\|v(y,\tau)\|_\infty\le K$  for every $\tau>0$. \qed

This is a useful consequence:

\begin{corollary}\label{thm.Linf.strong} Under the additional assumption that
\begin{equation}\label{barr.cond2}
C_0 M^{sp\beta}< L,
\end{equation}
there is a $\tau_0>0$ such that $\|v(\tau)\|_\infty \le L $ \ for every $\tau \ge \tau_0$. The precise condition is
\begin{equation}\label{barr.cond2.tau}
C_0 M^{sp\beta}\le L (1-e^{-\tau_0})^{\alpha}\,,
\end{equation}
 \end{corollary}

We have taken $a=1$ for convenience since then $\tau_0=0$, $x=y$  and $v(y,0)=u_0(x)$.
 The same formula holds with $(t+a)$ instead of  $(t+1)$ but then $C$ changes. We can now derive the  main comparison result.

\begin{theorem}[Barrier comparison] \label{thm.barrcomp}
Let $u$ be a solution with nonnegative data $u_0$ such that
$\|u_0\|_\infty\le L$ and $\|u_0\|_1\le M$. There exists a constant $C(L;M)>0$ such that whenever
\begin{equation}
u_0(y)\le C_1 \,r^{-(N+sp)} \quad \mbox{ for } \ r>2\,.
\end{equation}
with  $C_1\ge C(L;M)$, then the corresponding $v$-solution satisfies
\begin{equation}\label{barrcomp.v}
v(y,\tau)\le \widehat v(y;C_1) \qquad \forall y\in \ren, \tau>0\,.
\end{equation}
In other words,
\begin{equation}
u(x,t)\le\widehat u(x,t;C_1):= (1+t)^{-\alpha}\widehat v(|x|\,(1+t)^{-\beta};C_1) \qquad \forall x\in \ren, t>0.
\end{equation}
    \end{theorem}

\noindent {\sl Proof.} (i)  Let us pick some $\tau_0>0$ and apply the lemma to conclude that
we get $\|v(y,\tau)\|_\infty\le K(L,M)$  for every $\tau>0$. This gives a comparison between  $v(y,\tau) $ with $\widehat v(y;C_1)$  in the inner  cylinder \ $Q_0=B_2(0)\times (0,\infty)$ on the condition that $C_1\ge K(L,M)2^{N+sp}:=C(L,M)$.

(ii) We still  have to compare both functions in the exterior cylinder  $Q= \{y: |y|=r> 2\}\times (0,\infty)$.
Now, we have already proved that $\widehat v(y)$ is a pointwise supersolution for the equation in $Q$. The standard comparison theorem for the fractional evolution equation  applies to our $p$-Laplacian  equation, see similar details in \cite{VPLP2020}. Hence, in order to obtain comparison in $Q$, only the initial data and data in the complement \ $Q_0$ have to be checked  The latter has been established in (i), the former depends on the assumption on the initial data  $v_0(y)\le C_1 \,r^{-(N+sp)}$ for $r>2$. We conclude that \eqref{barrcomp.v} holds. \qed

\medskip

\noindent {\bf Spatial decay of solutions.} The barrier can be used to find a rate of space decay of the solutions which is uniform for bounded mass, bounded initial sup, and controlled initial tail.  In fact, under the conditions of Theorem \ref{thm.barrcomp}
it follows that in the outer region \ $|x|\ge c(t+1)^{\beta}$ \ we have
 \begin{equation}\label{decay.upur}
   u(x,t)\le C |x|^{-(N+sp)}(t+1)^{sp\beta}\,.
\end{equation}
Again, we see that the upper barrier decays with $|x|$ in the far field limit but it builds up with time.

\subsection{Barrier for the critical exponent}\label{sec.barr.cexpo}

We examine the remaining case $p=p_1$, where $N(2-p)=sp(p-1)$. Formally, the estimates from both regions we have studied coincide in suggesting a decay rate $\widehat v(y)=O(r^{-sp/(2-p)})=O(r^{-(N+sp)})$ for $p=p_1$. Normally, in such critical cases some type of correction is needed. For a recent example see \cite{BFV18e}, formula 3.4. \nc

\begin{proposition}\label{thm.barr.cc} The proof of the previous subsection holds in the critical case $p=p_1$ if we take  the following expression for the barrier in the outer region $r>2$:
\begin{equation}\label{form.barr.cclog}
\widehat v(y)=C_1r^{-(N+sp_1)}\log(r)^\gamma, \qquad  \gamma=\frac1{2-p_1},
\end{equation}
and $C_1$ is large enough.
\end{proposition}

\noindent {\sl Proof.} The contribution from the first-order term in \eqref{ineq.super1} is now
\begin{equation*}
-\beta \, r^{1-N}(r^N{\widehat v} )_r\approx \beta sp \, C_1r^{-N-sp}\log(r)^\gamma>0  \,,
\end{equation*}
which has a good sign, we call it the principal term. Using the same regions and notations as before, the contribution from $r\ge 1$ amounts to
$$
I_1\sim -c r^{-N-sp}\,.
$$
which is a lower value than the principal estimate. The contribution $I_2$ is
$$
|I_2| \le
\int_{1}^{r_0/2} c  C_1^{p-1} r^{-(N+sp)(p-1)}\log(r)^{\gamma (p-1)}r_0^{-(N+sp)}r^{N-1}dr\,.
$$
Since $N(2-p)=sp(p-1)$ precisely for $p=p_1$,  we estimate this contribution as
$$
I_2\approx c C_1^{p-1} r_0^{-(N+sp)} \int_{1}^{r_0/2} \log(r)^{\gamma (p-1)}r^{-1}dr= c  C_1^{p-1} r_0^{-(N+2s)}\log(r_0)^{1+\gamma (p-1)},
$$
which will be a lower value with respect to the principal term whenever $1+\gamma (p-1)<\gamma$, i.e. for $\gamma(2-p)>1$.
For the limit case $\gamma(2-p)=1$ we need the condition that $C_1$ must be large enough.

 Finally, for $I_3$ we write
$$
|I_3|\le c\int_{B_{r_0/2}(x_0)} \frac{(\widehat v(r_0)-\widehat v(r))^{p-1}} {|z|^{N+sp}}\,dz,
$$
with $z=x-x_0$. This integral is convergent at the singularity as before. Taking into account that for $r\gg 1$ we have $\widehat v'(r)\sim -(N+sp)\widehat v(r)/r$, we get
$$
|I_3|\le cC_1^{p-1}   r_0^{-(N+sp)(p-1)-sp}\log(r_0)^{\gamma(p-1)}
$$
which for $p=p_1$ gives
$$
|I_3|\le c  C_1^{p-1} r_0^{-(N+sp_1)}\log(r_0)^{\gamma(p-1)},
$$%
so the rate of decay of this term is faster than the principal term, and we are done proving that \eqref{form.barr.cclog} gives a good barrier if $\gamma\ge 1/(2-p_1)$\nc. \qed

The conclusions of Theorem \eqref{thm.barrcomp} read as follows for $p=p_1$.

\begin{theorem}[Barrier comparison] \label{thm.barrcomp_1}
Let $u$ be a solution with nonnegative data $u_0$ such that
$\|u_0\|_\infty\le L$ and $\|u_0\|_1\le M$. There exists a constant $C(L;M)>0$ such that whenever
\begin{equation}
u_0(y)\le C_1 \,r^{-(N+sp_1)} \log (r)^\gamma \quad \mbox{ for } \ r>2\,.
\end{equation}
with  $C_1\ge C(L;M)$, then the corresponding $v$-solution satisfies
\begin{equation}\label{barrcomp.v}
v(y,\tau)\le \widehat v(y;C_1) \qquad \forall y\in \ren, \tau>0\,.
\end{equation}
In other words,
\begin{equation}
u(x,t)\le\widehat u(x,t;C_1):= (1+t)^{-\alpha}\widehat v(|x|\,(1+t)^{-\beta};C_1) \qquad \forall x\in \ren, t>0.
\end{equation}
    \end{theorem}

We will go back to the issue of upper bounds for this case in Section \eqref{sec.p_1}.

\medskip

\section{Mass conservation}\label{sec.masscon}

We now proceed with the mass analysis. The main result is the conservation of the total mass for the Cauchy problem posed in the whole space with nonnegative data.

\begin{theorem}\label{thm.mc} Let $p>p_c$. \nc Let $u(x,t)$ be the semigroup solution of Problem \eqref{frplap.eq}, \eqref{init.data}, with $u_0\in L^1(\ren)$, $u_0\ge 0$. Then for every $t>0$ we have
\begin{equation}
\int_{\ren} u(x,t)\,dx=\int_{\ren} u_0(x)\,dx.
\end{equation}
\end{theorem}

Before we proceed with the proof we make a reduction:   We may always assume that $u_0\in L^1(\ren)\cap L^\infty(\ren)$ and compactly supported. If mass conservation is proved under these assumptions, then it follows for all data $u_0\in L^1(\ren)$ by the semigroup contraction property.

We recall that the $L^1$ mass is not conserved  for Problem \eqref{frplap.eq} with positive and integrable data when the exponent lies in the range $1<p<p_c$ because of the appearance of solutions that vanish in finite time, see  \cite{BS2020} and Section \ref{sec.vfd} below. This phenomenon also happens for the standard porous medium and $p$-Laplacian equation, see for instance the original source \cite{BC81}, and a long account on this problem in \cite{VazSmooth}.\nc

The proof of the theorem is divided into several cases in order to graduate the difficulties. Note under our assumptions $0<sp<2$.

\subsection{First case: $N=1<sp$.}\label{ssec.4.1} Here the mass calculation is quite straightforward. We do
a direct calculation for the tested mass. Taking a smooth and compactly supported test function $\varphi(x)\ge0 $, we have for $t_2>t_1>0$:
\begin{equation}\label{mass.calc}
\left\{\begin{array}{l}
\displaystyle \left|\int u(t_1)\varphi\,dx-u(t_2)\varphi\,dx\right|\le \iiint
\left|\frac{\Phi(u(y,t)-u(x,t))(\varphi(y)-\varphi(x)}{|x-y|^{N+sp}}\right|\,dydxdt\\[10pt]
\le \displaystyle \left(\iiint  |u(y,t)-u(x,t)|^p\,d\mu(x,y)dt\right)^{\frac{p-1}{p}}
\left(\iiint |\varphi(y)-\varphi(x)|^p\,d\mu(x,y)dt\right)^{\frac{1}{p}}\,,
\end{array}\right.
\end{equation}
with space integrals over $\ren$ (here, $N=1$) and time integrals over $[t_1. t_2]$. Use now the sequence of test functions $\varphi_n(x)=\varphi(x/n)$ where $\varphi(x)$ is a cutoff function which equals 1 for $|x|\le 2$ and zero for $|x|\ge 3$. Then,
$$
\int \int |\varphi_n(y)-\varphi_n(x)|^p\,d\mu(x,y)=
n^{N-sp}\int \int |\varphi(y)-\varphi(x)|^p\,d\mu(x,y)
$$
and this tends to zero as $n\to\infty$. Using \eqref{en.est} we conclude that the triple integral involving $u$ is also bounded in terms
of $\|u(\cdot,t_1)\|_2^2$, which is bounded independently of $t_1$. Therefore, taking the limit as $n\to\infty$ so that $\varphi_n(x)\to 1$ everywhere, we get
$$
\int u(x,t_1)\,dx=\int u(x,t_2)\,dx,
$$
hence the mass is conserved for all positive times for data in $L^2\cap L^1$. The statement of the theorem needs to let $t_1\to 0$, but this can be done thanks to the continuity of solution of the $L^1$ semigroup as a curve in $L^1(\ren)$.

The limit case $1=sp$ also works by revising the integrals, but we get no rate.

\subsection{Case $N\ge sp$}\label{ssec.4.2}
Since $sp<2$, in dimensions two or more we always have $N>sp$. We can also have $N=1> sp$. The proof of mass conservation in this case does not seem to be easy and we offer a proof with some delicate step. Actually, in order to obtain the mass conservation in this case  we need to use a uniform estimate of the  decrease of the solutions in space so that they help in estimating the convergence of the integral. This will be done by using the barrier estimates that we have obtained in previous sections.

\noindent {$\bullet$} We go back to the first line of \eqref{mass.calc}. The proof relies on some calculations with the multiple integral in that line. We also have to consider different regions. We first deal with exterior region $A_n=\{(x,y): |x|,|y|\ge n\}$,  where recalling \eqref{mass.calc} we have
\begin{equation*}\label{mass.calc2}
\begin{array}{c}
 \displaystyle I(A_n):= \int_{t_1}^{t_2}\iint_{A_n} \frac{|\Phi(u(y,t)-u(x,t))|\,|\varphi_n(y)-\varphi_n(x)|}{|x-y|^{N+sp}}\,dydx\,dt\\[10pt]
\le \displaystyle \left(\iiint  |u(y,t)-u(x,t)|^{p}\,d\mu(x,y)\,dt \right)^{\frac{p-1}{p}}
\left(\iiint |\varphi_n(y)-\varphi_n(x)|^{p}\,d\mu(x,y)dt \right)^{\frac{1}{p}}
\end{array}\end{equation*}
which we write as $I=I_1. I_2$. In the rest of the calculation we omit the reference to the limits that is hopefully understood.

We already know that $I_2 \le C_p\,n^{(N-sp)/p}(t_2-t_1)$. On the other hand,
we want to compare $I_1$ with the dissipation $D_\ve$ of the $L^r$ norm, for $r=1+\ve$. We recall that  according to \eqref{en.est.q}
\begin{equation*}\label{mass.calc3}
\displaystyle D_\ve= \iiint |u(y,t)-u(x,t)|^{p-1}\,|u^\ve(y,t)-u^\ve(x,t)|  \,d\mu dt \le C(\ve,p)\int |u|^{1+\ve}(x,t_1)\,dx\,,
\end{equation*}
which is bounded by $  C(\ve,p,u_0)$. Next, we use the elementary equivalence: for all $\ve\in (0,1)$ and all $a,b>0$ we have
$$
|a^\ve-b^\ve| \ge C(\ve) (a-b) (a+b)^{\ve-1}\,.
$$
It follows that
$$
 \displaystyle D_\ve \ge C_\ve \iiint |(u(y,t)-u(x,t)|^{p}\,(|u|^\ve(y,t)+|u|^\ve(x,t))^{\ve-1}  \,d\mu dt\,.
$$
After comparing the formulas, we conclude that
$$
I_1^{p/(p-1)}\le C\,D_\ve\|2u\|_\infty^{1-\ve}\,.
$$
In view of the value of $u$ in the region $A_n$, $u\approx n^{-\gamma}$, we have $I(A_n)\le C n^{-\sigma}$ with
$$
\sigma =\frac1{p}\left(\gamma(p-1)(1-\ve)-(N-sp)\right).
$$
Let us examine the cases. In Case 1 for $p_c<p<p_1 $, we have $\gamma=sp/(2-p)$, hence
$$
p\sigma =\frac{sp(p-1)}{2-p}+sp-N-\ve\frac{sp(p-1)}{2-p}\approx \frac{sp}{2-p}-N,
$$
which is positive for $p>p_c$. In Case 2 $p_1<p<2$ we have $\gamma=N+2s$, hence
$$
p\sigma = N(p-2)+ sp^2-\ve (N+sp)(p-1)\approx sp+N(p-2)+ s(p^2-1)
$$
which is positive for $\ve$ small, in both cases we get the vanishing in the limit  $n\to\infty$ of this term that contributes to the conservation of mass.  The same happens in the critical case $p=p_1$.  Note that the argument holds for all $p> p_c$ and $0<s<1$.

\noindent $\bullet$  We still have to make the analysis in the other regions. In the inner region $B_n=\{(x,y): |x|,|y|\le 2n\}$ we get
$\varphi_n(x)-\varphi_n(y)=0$, hence the contribution to the integral \eqref{mass.calc} is zero. It remains to consider the cross regions
$C_n=\{(x,y): |x|\ge 2n ,|y|\le n\}$ and $D_n=\{(x,y): |x|\le n ,|y|\ge 2n\}$. Both are similar so we will look only at $C_n$. The idea is that we have an extra estimate: $|x-y|>n$ so that
\begin{equation*}
 \begin{array}{c}
 \displaystyle I(C_n)\le n^{-(N+sp)}\int_{t_1}^{t_2}\iint_{C_n}|\Phi(u(y,t)-u(x,t))|\,|\varphi_n(y)-\varphi_n(x)|\,dydx\,dt\\[10pt]
 \displaystyle  \le Cn^{-(N+sp)}(\int_{|y\le n} dy)( \int dt \int_{n\le |x|\le 3n} |u(x,t)|^{p-1}dx)\,.
 \end{array}
\end{equation*}
Note that for $1<p<2$ we have $0<p-1<1$. We now argue as follows:
$$
\int_{n\le |x|\le 3n} |u(x,t)|^{p-1}dx\le n^{N(2-p)}\left(\int_{n\le |x|\le 3n} |u(x,t)|\,dx\right)^{p-1}
$$
Therefore, using the $L^1$ contraction of the semigroup we get
$$
 \displaystyle I(C_n)\le n^{-(N+sp)} n^N (t_2-t_1)\|u_0\|_1^{p-1}= n^{N(2-p)-sp} n^N (t_2-t_1)\|u_0\|_1.
 $$
 We have $N(2-p)-sp<0$ if $ p>2N/(N+s)$, which is true. Therefore, $I(C_n)$
tends to zero as $n\to \infty$ with a power rate. Same  for $I(D_n)$. This concludes the proof. Note that these regions overlap but that is no problem. \qed

\medskip

\noindent {\bf Signed data.} Theorem \ref{thm.mc} holds also for signed data and solutions. However, the denomination mass for the integral over $\ren$ is physically justified only when $u\ge 0$. For signed solutions the theorem talks about conservation of the whole space integral. The above proof has be  reviewed. Subsection \ref{ssec.4.1} needs no change. As for Subsection \ref{ssec.4.2}, the elementary equivalence has to be written for all $a,b\in\re$
$$
|a^\ve-b^\ve| \ge C(\ve) |a-b| (|a|+|b|)^{\ve-1}\,.
$$

\noindent {\bf Remark.} As said before, the law of mass conservation does not hold for exponents $1<p<p_c$, the most extreme counterexample  being given by solutions that vanish identically after a finite time, see Section \ref{sec.vfd}.


\section{Positivity of nonnegative solutions}\label{sec.globalpos}

We prove in this section that nonnegative and nontrivial solutions are actually positive everywhere if $p>p_c$, a range where mass conservation holds. We give first a more precise quantitative statement (but only local positivity) since it will be needed as a technical tool. Then, we state and prove a general qualitative result.

Let us recall that in the limit case $s = 1$, with $1<p < 2$ fixed, we get the standard $p$-Laplacian equation, where  everywhere positivity is true for all nonnegative solutions due to the property of infinite propagation of the fast $p$-Laplacian equation (unless there is complete extinction in finite time, something that may happen only for very fast diffusion). The difficulty in finding convenient explicit lower barriers has forced us to introduce new ideas to tackle the present fractional equation.

\subsection{A quantitative positivity lemma}

As a consequence of mass conservation and the existence of the upper barrier, we obtain a partial positivity lemma for certain solutions of the equation.

\begin{lemma}\label{lem.pos}  Let $v$ be the solution of equation \eqref{eq.resc}   such that the initial data $v_0\ge 0$ is a nontrivial bounded function, bounded above by an upper barrier $G(y)$ as in Sections \ref{sec.barr1}, \ref{sec.barr2}, and we also assume that $\int v_0(y)\,dy=M>0$. Besides, $v_0$ is radial and radially decreasing.   Then,
there is a continuous nonnegative function $\zeta(y)$, positive in a ball of radius $r>0$, such
that for every $\tau>0$
\begin{equation}
v(y,\tau)\ge  \zeta(y) \quad \mbox{ for all } \ y\in \ren, \ \tau>0.
\end{equation}
In particular, we may take $\zeta(y)\ge c_1>0$ in $ B_{r_0}(0)$ for suitable $r_0$ and $c_1>0$, to be computed  below. The function $\zeta$ depends only on the data $G$, $M$ and $\|v_0\|_\infty$.
\end{lemma}

\noindent {\sl Proof.} (i) We know that for every $\tau>0$ the solution $v(\cdot,t)$ will be nonnegative, radial, and radially nonincreasing. 
By Sections \ref{sec.barr1}, \ref{sec.barr2}, there is the same upper barrier $G(y)$ will be on top of $v(y,\tau)$ for every $\tau$. Since $G$ is integrable at infinity, for every $\varepsilon>0$ small there is $R(\ve)>0$
such that
$$
\int_{\{|y|>R(\ve)\}} v(y,\tau)\,\dy\le
\int_{\{|y|>R(\ve)\}} G(y)\,\dy\le \ve
$$
for all $\tau>0$. Moreover, there is a radius $r_0>0$ such that
$$
\int_{\{|y|<r_0\}} v(y,\tau)\,\dy\le M/3
$$
for all $\tau>0$ since $v$ will be uniformly bounded for all $\tau$.  Therefore,
$$
\int_{\{r_0\le |y|\le R(\ve)\}} v(y,\tau)\,\dy\ge M-\ve-M/3 >M/2.
$$
We use this estimate and the fact that $v$ is monotone in $r=|y|$  to conclude that
$$
v(r_0,\tau)(R(\ve)^N-r_0^N) \ge c(n)M/2,
$$
hence $v(r,\tau)\ge c_1$ for all $r\le r_0$ and $\tau>0$, with $c_1=c(N,s,p,M,R,L)$. \qed


\subsection{Everywhere positivity}\label{ever.pos}

 We will prove positivity in the sense of locally strict positivity. This means that near a point of positivity  $(y_0,\tau_0)\in\ren\times \re_+$ there exist a small $\ve>0$ and a constant $c=c>0$ (depending on the point) such $v\ge c$ a.e. in the neighbourhood $Q_\ve=B_\ve(y_0)\times (\tau_0-\ve,\tau_0+\ve)$. Since we have proved later on continuity the concept of locally strict positivity is just strict positivity,but we have kept the argument because it might be useful in other contexts. \nc

\begin{theorem}\label{lem.genpos}
Let $v$ be the solution of equation \eqref{eq.resc}  with  nontrivial and locally integrable initial data $v_0\ge 0$. Then $v(y,\tau)$ is (locally strictly) positive everywhere.
\end{theorem}

\noindent {\sl Proof.} \bf I. \rm  We may  assume by approximation from below that the function is bounded and integrable,   $\|v_0\|_\infty\le L, $ and $\|v_0\|_1\le M $, and also that it is supported in the ball $B_R(0)$. By scaling we may also assume that $R=1$.

In the first part of the proof we will establish positivity away from the initial support.  There are a number of steps similar to the last lemma.

 (i) Arguing as before, there is a stationary upper barrier $G(y)$, depending on $L$. Therefore, for $R_1(L,M)\gg1$ large enough we have
$$
\int_{\{|y|>R_1\}} v(y,\tau)\,\dy\le
\int_{\{|y|>R_1\}} G(y)\,\dy\le M/3\,,
$$
valid for all $\tau>0$.

(ii) Due to the (possibly increasing) $L^\infty$ bound for $v$, we find that if $L$ is small enough and $\tau\le \tau_1$
$$
\int_{\{|y|\le 6\}} v(y,\tau)\,\dy\le M/3.
$$
It follows that the mass in the complement of these two regions (an annulus) will be equal or larger than $M/3$.

(iii) Next step is to use the available monotonicity. We apply Aleksandrov's Principle (see Proposition \ref{Alek}) on the $u$-solution to prove that for $t\ge 0$ and $|x|\ge k>1$ the solution $u(x,t)$ is monotone in a cone of outgoing directions $K_\theta(x)$, and the aperture $0<\theta<\pi/2$ grows with increasing $k>1$. The same property happens to the solution $v(y,\tau)$ since the space variable shrinks with time but the aperture of the cones is kept under the scaling. Let us consider the set of cones with vertex in the surface $|y|=4$. We only need to find a finite number of those cones to cover the whole exterior domain $|y|\ge 6$, let $j>4$ be a sufficient number. Let $y_1,y_2,\cdots, y_j$ be the vertex points and let $\widehat K_\theta(y_i)$ be $K_\theta(y_i)\cap\{y: |y|\le R_1\}$. Then the $y$-cone containing the largest mass, with vertex at say $y=y_1$, will contain a mass larger or equal than $M/3j$. Since $y_1$ is point of maximum for $v$ in that cone,  we estimate the value
$$
v(y_1,\tau)\le cM/j|\widehat K_\theta(y_1)|\sim c_1M / R_1^N=c_2.
$$
We now use the second assertion of Aleksandrov's principle to assert that
$$
v(y,\tau)\ge v(y_1,\tau)\ge c_2
$$
for  $y$ such that $1<|y|<2$ and $\tau\ge \tau_1$. In terms of $u$, this means the occurrence of a set of positivity
that covers the whole region $\{(x,t): \ t\ge t_2, \ (1+t)^\beta\le |x|\le 2(1+t)^\beta\}$ for $t\ge t_2\ge t_1$.

 (iv) Due to the monotonicity in time of the function $u(x,t)t^{1/(2-p)}$ (see Subsection \ref{apriori}) we conclude that $u$ is locally strictly positive in the cylinder
$$
Q= A\times (0,t_2), \qquad A=\{ x:  (1+t_2)^\beta\le |x|\le 2(1+t_2)^\beta\}.
$$
Adding sets with different values of $t_2\ge t_1$, we get the conclusion that $u$ is locally strictly positive in the outer set
\begin{equation}\label{pos.set2}
Q^*=\{ (x,t): |x|\ge  (1+t_1)^\beta, \ 0<t<t_1\}.
\end{equation}
This completes the step of instantaneous creation of a positivity outer region.

\medskip

\noindent \bf II. \rm We need a further argument to cover the possible hole $H=\{|x|\le  (1+t_1)^\beta \}$ left by formula \eqref{pos.set2} near the original support. We do it in three steps, i.e., constructing a small solution, displacing it in space far enough, and then comparing with $u$ for later times.

(i) We first do a special case of the previous computation where $\hat u_0(x)=\chi(B_{1}(0))$. The argument leads to a time $\hat t_1>0$ (depending only on $N,s,p$) such that  $\hat u(x,t)$ is locally strictly positive at all points outside the ball \ $B_{\hat R}(0)$, \ $\hat R=(1+\hat t_1)^\beta$, \ for times $0<t<\hat t_1$. We need to stretch that set and we do it by defining a rescaled initial function
\begin{equation}\label{eps.c}
\hat u_{\ve,c}(x,0)= c \chi(B_{\ve}(0)).
\end{equation}
According to the scaling transformation of Subsection \ref{ssec.scaling} we get a solution
$$
\hat u_{\ve,c}(x,t)=c\,\hat u(\ve^{-1} x, c^{p-2}\ve^{-sp} t)
$$
We impose the condition $\ve^{sp}c^{2-p}=\delta$ small enough. In this way the positivity set of $\hat u_{\ve,c}$ is a scaling of the one for $\hat u$, i.e,.
$$
\hat Q^*_{\ve,c}=\{ x: |x|\ge \ve (1+\delta \hat t_1\,)^\beta\}\times (0,\delta \hat t_1).
$$

(ii) We now shift the origin of coordinates to  a point $x_0$ with $|x_0|> r_1^*= 2(1+t_1)^\beta$, and we consider the last solution with data shifted in space by the amount $x_0$, and we also take the time $t_1$ of Part I  as the origin of times:
\begin{equation}\label{eps.c.2}
\overline u(x,t+t_1)=\hat u_{\ve,c}(x-x_0,t).
\end{equation}
The corresponding positivity set for $\overline u$ is then
$$
\overline Q^*=\{ x: |x-x_0|\ge \ve (1+\delta \hat t_1\,)^\beta\}\times (t_1,t_1 +\delta \hat t_1).
$$
This set covers $H$ for $t_1<t<t_1 +\delta \hat t_1$ as long as $|x_0|\ge 2\ve (1+\delta \hat t_1\,)^\beta.$

(iii) To finish the proof, we compare $u(x,t)$ with $\overline u(x,t))$ at $t=t_1$ and conclude that for $t>0$
$$
u(x,t+t_1)\ge \overline u(x,t)).
$$
This happens if $\ve$ and $c$ in \eqref{eps.c.2} are small enough (depending on the point $(x_0,t_1)$. We obtain positivity of $u$ in $B_{r_1*}(0)$ for $t=t_1+\ve^{sp}\hat t_1$. By time monotonicity, it holds for all previous times. \qed


\noindent {\bf Remark.} By comparison the result is true even if the initial support is not bounded, and the solution is merely
a function in a Lebesgue space $L^q(\ren)$, $1\le q\le \infty$. But notice that the compact support condition plays a big role in the proof.


\section{$L^1$ dissipation for differences}\label{sec.diff}

In subsequent sections we will need to estimate the dissipation of the difference $u=u_1-u_2$ in the framework of the $L^1$ semigroup.
This is a very delicate result, typical of the fractional operator. We multiply the equation by $\phi=s_+(u_1-u_2)$, where $s_+$ denotes the sign-plus or Heaviside function, and then integrate in space and time. We get in the usual way, with $u=u_1-u_2$, $u_+=\max\{u,0\}$,
\begin{equation}\label{diff.est.1}
\begin{array}{c}
\displaystyle \int u_+(x,t_1)\,dx-\int u_+(x,t_2)\,dx=\int_{t_1}^{t_2}\int s_+(u)u_t\,dx\\ [6pt]
\displaystyle  =\int_{t_1}^{t_2}dt\int ({\mathcal L}_{s,p}  u_1-{\mathcal L}_{s,p}  u_2)\,s_+(u_1-u_2)\,dx =\\ [6pt]
\displaystyle \int_{t_1}^{t_2}dt \iint \left\{|u_1(x,t)-u_1(y,t)|^{p-2}(u_1(x,t)-u_1(y,t)) \right. \\[10pt]
 \displaystyle \left.-|u_2(x,t)-u_2(y,t)|^{p-2}(u_2(x,t)-u_2(y,t))\right\}\\[8pt]
\big( s_+(u(x,t))-s_+(u(y,t))\big)\,d\mu(x,y)\,.
\end{array}
\end{equation}
We recall that $s_+(u(x,t))=1$ only when $u_1(x,t)>u_2(x,t)$, and
$s_+(u(y,t))=0$ only when $u_1(y,t)<u_2(y,t)$. If we call the last factor in the above display
$$
I=s_+(u_1(x,t)-u_2(x,t))-s_+(u_1(y,t)-u_2(y,t))\,,
$$
we see that $I=1$ if $u_1(x,t)>u_2(x,t)$ and $u_1(y,t)\le u_2(y,t)$. Therefore, on that set
$$
u_1(x,t)-u_1(y,t)>u_2(x,t)-u_2(y,t).
$$
In that case we examine the other factor,
$$
F= |u_1(x,t)-u_1(y,t)|^{p-2}(u_1(x,t)-u_1(y,t)) -|u_2(x,t)-u_2(y,t)|^{p-2}(u_2(x,t)-u_2(y,t))\,,
$$
and conclude that it is positive. The whole right-hand integrand is positive.

In the same way, $I=-1$ if $s_+(u(x,t))=0$ and $s_+(u(y,t)=1$ i.e., only when $u_1(x,t)\le u_2(x,t)$ and $u_1(y,t)>u_2(y,t)$. Then,
$u_1(x,t)-u_1(y,t)<u_2(x,t)-u_2(y,t)$ and $F<0$. The whole right-hand integrand is again positive. We conclude that

\begin{proposition} In the above situation we have the following dissipation estimate:
\begin{equation}\label{diff.est.1b}
\begin{array}{c}
\displaystyle \int (u_1-u_2)_+(x,t_1)\,dx-\int (u_1-u_2)_+(x,t_2)\,dx\\ [8pt]
\ge \displaystyle\iiint_D \left| |u_1(x:y,t)|^{p-2}u_1(x:y,t) -|u_2(x:y,t)|^{p-2}u_2(x:y,t)\right|\,d\mu \,dt.
\end{array}
\end{equation}
where $u_1(x:y,t)=u_1(x,t)-u_1(y,t)$, $u_2(x:y,t)=u_2(x,t)-u_2(y,t)$, and $D\subset \re^2$ is the domain where
$$
\{ u(x,t)>0, \  u(y,t)\le 0\}\cup\{  u(x,t)\le 0, \  u(y,t)>0\}\\,
$$
that includes the whole domain where $u(x,t)\,u(y,t)<0$. There is no dissipation on the set where $u(x,t)\,u(y,t)>0$.
\end{proposition}


\section{Existence of the fundamental solution}\label{sec.fs}

This section deals only with nonnegative solutions. The first result is

\begin{theorem}\label{thm.exfs} Let $2>p>p_c$. For any value of the mass $M>0$ there exists a fundamental solution $U(x,t)$ of Problem \eqref{frplap.eq}-\eqref{init.data} having the following properties:

 (i) It is a nonnegative strong solution of the equation in all $L^q$ spaces, $q\ge 1$, for $t\ge t_0>0$.

 (ii) $U(\cdot,t)$ is  radially symmetric and decreasing in the space variable for every $t>0$.

 (iii) It decays in space as predicted by the barriers, $u(t)=O(|x|^{-N-sp})$ or $u(t)=O(|x|^{-sp/(2-p)})$ depending on the $p$ range, with a logarithmic correction for $p=p_1$.

 (iv) It decays in time $O(t^{-\alpha})$ uniformly in $x$.
\end{theorem}

The result is completely similar to the case $p>2$ treated in the companion paper \cite{VPLP2020}, it is given here only for reference. The proof is the same,  and the will be omitted since  will prove a stronger result in the next subsection in great detail. See the remarks that are made  in \cite{VPLP2020}. In particular, uniqueness of such a solution is not proved.

We are going to improve this result to complete the proof of Theorem  \eqref{thm.ssfs} by proving existence of  a unique self-similar fundamental solution, which is a more specific object.

\begin{theorem}\label{thm.exfs2} If $p_c<p<2$ there is a fundamental solution of Problem \eqref{frplap.eq}-\eqref{init.data}  with the properties of Theorem \ref{thm.exfs} that is also self-similar. Moreover, the self-similar fundamental solution is unique.
The profile $F_M(r)$ is a nonnegative and radial continuous function that is nonincreasing along the radius, is positive everywhere and satisfies with decay rate prescribed in (iii) of previous theorem. Moreover, for $p_c<p<p_1$ we have the precise bound
\begin{equation}\label{decay.Fbarrier.lr}
F_M(y)\le C_1^*\,|y|^{-sp/(2-p)} \,,
\end{equation}%
where  $C_1^*>0$ is the constant introduced in Theorem \eqref{thm.barr.lr}, so the estimate is independent of $M$.
\end{theorem}

We will prove further below that the asymptotic estimate is exact, i.e., $F_M$ goes to zero at infinity with the exact rate $F_M(r) \sim C_1^* r^{-sp/(2-p)}$ for $p<p_1$, while \
$$
F_M(r)\approx r^{-(N+sp)}\quad \mbox{ if } p>p_1,
$$
but this has to wait for the study of lower bounds, done in Sections \ref{sec.pos.lr} and \ref{sec.pos.tail.ur}.


\subsection{Proof of uniqueness of the self-similar profile}

This proof is valid for every $p_c<p<2$, and even for $p\ge 2$ as shown in \cite{VPLP2020}.
We know that any self-similar profile $F$ is bounded, radially symmetric and non increasing.
 We know that $0\le F\le C$, that $F$ decreases like a power of $r$, either $F\le Cr^{-(N+sp)}$ or $F\le Cr^{-sp(/2-p)}$, with a log correction for $p=p_1$.
We prove regularity for the profile by using the regularity of the equation. We recall that
$U_t(x,1)=-\beta \nabla \cdot(xF) $ is bounded, so that $F$ is a $C^1$ function for $r>0$.

 The main step is to use mass difference analysis, since this is a strict Lyapunov functional,
hence we arrive at a contradiction when two self-similar profiles meet. This is an argument
taken from the book \cite{Vazpme07}. It goes as follows: We take two profiles $F_1$ and $F_2$ and assume the same mass $\int F_1\, dx=\int F_2\, dx=1$. If $F_1$ is not $F_2$ they must intersect and then $\int (F_1-F_2)_+dx=C$ is not zero and $\int (F_2-F_1)_+dx=C$ too. By self-similarity,  $\int (U_1-U_2)_+dx=C$ must be constant in time. But we have proved in Section \ref{sec.diff} that whenever $C>0$ at one time, the integral  must be a decreasing quantity in time. \qed


\subsection{Existence of the fundamental self-similar profile}\label{ssec.exselfsim.ur}

We do a separate analysis in the different exponent subranges. It will be obtained by a method that in a first step  proves existence of a periodic $v$ solution.

\medskip

 \noindent {\bf I. Case  $p_1<p<2$}. We prove here the existence of SS FS in the upper  sublinear range. The absence of a global supersolution in the barrier construction of Section \ref{sec.barr2} forces us develop a delicate fixed point argument
that uses ideas of our work with Feo and Volzone \cite{FVV2020}, and previously in \cite{VPLP2020}. It is as follows.

(i)  Let $X = L^1(\mathbb{R}^N)$. We consider a subset $\mathcal K$ of $ X$ defined as follows:
$ \mathcal K=\mathcal K(M,L)$ is the set of all $\phi\in L_{+}^{1}(\ren)\cap L^{\infty}(\re^{N})$  such that \

(a) $\int \phi(y)\,dy = M,$ for some $M>0$, \

(b) $\,\phi\,$ is radially symmetric and nonincreasing w.r.to $r=|x|$, \

(c) $\phi$  is bounded above by a fixed admissible barrier function $G$ as in Theorem \ref{thm.barrcomp}, and

(d) $\phi$  is uniformly bounded above by a constant $L>0$.

\medskip

\noindent It is easy to see that  $\mathcal K (M,L)$ is a non-empty, convex, closed and bounded subset with respect to the norm of the Banach space $X$.

(ii) Next, we prove the existence of periodic orbits. For all $\phi\in \mathcal K(M,L)$ we consider the solution $v(y,h)$ to equation \eqref{frplap.eq} starting at $\tau = 0$ with data $v(y, 0) =\phi(y)$.
 For $h>0$ we  consider the semigroup map $S_h : X \rightarrow X$ defined by $\mathcal{S}_h(\phi) =v(y, h)$.

Fix $h>0$ and let us check that  $S_h(\mathcal K(M,L))\subset \mathcal K(M,L)$.
Conditions (a), (b), and (c) are a consequence of known properties, based on Theorem \ref{thm.barrcomp}. We now impose the smallness condition \eqref{barr.cond2} on $M$ with respect to $L$. Then, there is an $h>0$ such that
\begin{equation}
C_0 M^{sp\beta}\le L (1-e^{-h})^{\alpha}
\end{equation}
and we get the estimate \ $\|v(\tau)\|_\infty \le L $ \ for every $\tau \ge h$. This settles condition
(d) for $\phi_1=S_h \phi$.

(iii) The relative compactness of $S_h(\mathcal K(M,L))$ comes from known regularity theory. We may even assert uniform local positivity using Lemma \ref{lem.pos} since the mass is constant in time.
It now follows from the Schauder Fixed Point Theorem,  \cite{EvansPDE},
that there exists at least fixed point $\phi_h \in \mathcal K$, \textit{i.\,e.,} $ S_h(\phi_h) = \phi_h$. The fixed point is in $\mathcal K$, so it is not trivial because of the conservation of mass and the $L^\infty$ bound.
Iterating the equality we get periodicity for the orbit $v_h(y, \tau)$ starting at $\tau = 0$
\begin{equation}\label{fn.periodic}
v_h(y,\tau+ kh) = v_h(y,\tau )\quad  \forall \tau > 0,
\end{equation}
valid for all integers $k\geq1$. It is not a trivial orbit, $v_h\not\equiv 0$.

We will abandon here line of work proposed in \cite{VPLP2020} consisting in passing to the limit in the family $v_h$ as $h\to 0$ to obtain a stationary solution. This has the problem that  according to the smallness condition the proof does not work unless we lower the mass of the solutions, and in this way may land on the trivial stationary solution $\widehat{v}\equiv 0$.\nc

\medskip

\noindent (iv) Here is the new ingredient. \sl We claim that any periodic solution like  $v_h$ defined in \eqref{fn.periodic} must be stationary in time. \rm The proof follows the lines of the uniqueness proof of previous subsection. Thus, if $v_1$ is a periodic solution that is not stationary, then $v_2(y,\tau)=v_1(y,\tau+ c)$ must be different from $v_1$ for some $c>0$, and both have the same mass. With notations as above we consider the functional
 \begin{equation*}
 J[v_1,v_2](\tau)=\int (v_1(x,\tau)-v_2(x,\tau))_+\,dx.
 \end{equation*}
 By known  $L^1$ accretivity of the operator,  this is a Lyapunov functional, i.e., it is nonnegative and nonincreasing in time.
 By the periodicity of $v_1$ and $v_2$, this functional must be periodic in time. Combining those properties we conclude that it is constant.
 We have to decide whether it is a positive constant or zero. In the latter  case, we arrive at a contradiction with the assumption that the solutions are different and we are done.
 To eliminate the possibility that it is  a positive constant, we point out that two different radial solutions with the same mass  must intersect. At this moment we apply the dissipation results of Section \ref{sec.diff} that imply that the functional must be decreasing in time. This contradicts what was already proven about periodicity. \qed

\medskip

(v) From this moment on, we set $F(y)=S_1(\phi_1)(y)$, for some fixed point $\phi_1$. Going back to the original variables, it means that the
corresponding function
\begin{equation*}
\widehat{u}(x,t)=(t+1)^{-\alpha} F(x\,(t+1)^{-\beta})
\end{equation*}
is a self-similar solution of equation \eqref{frplap.eq}. We know that the mass of $F$ is the same $M$.

(vi) Initial data. It is easy to see that the family $\rho_t(x)=\widehat{u}(x,t)$ is a nonnegative mollifying sequence as $t\to0$
since the functions are nonnegative, continuous and the integral outside of  a small ball goes to 0 as $t\to 0$ (use the upper barrier to prove it). Therefore, since there is conservation of mass in the sequence, by standard theory
\begin{equation}\label{conv delta}
\lim_{t\rightarrow 0^+}
\int_{\mathbb{R}^N} \widehat{u}(x,t)\varphi(x)\,dx=M\,\varphi(0)
\end{equation}
for all $\varphi \in C_0^\infty(\mathbb{R}^N).$

\medskip

(vii) Local positivity: we know from the proof that $\widehat{v}(y)\ge \zeta(y)$, which is positive in the ball $B_0(r_0) $ of radius $r_0>0$. This means positivity for $\widehat{u}(x,t)$ in sets of the form $|x|\ge r_0^\beta$. Positivity everywhere follows from Theorem \ref{lem.genpos}. This argument applies for $p_c<p<2$.

\medskip

(viii) We know from the proof that $M_1>0$. If $M_1=M$, we are done. If $M_1\neq M$ we apply the mass changing scaling transformation \eqref{scal.trn2}. In order to consider a negative mass, $M<0$, the fundamental solution is obtained by just putting
$U_M(x,t)=-U_{-M}(x,t)$.\qed

\medskip

\noindent {\bf II. Case $p=p_1$.}  The proof works with minor modifications for $p=p_1$. In the statement only the space decay of point (iii) changes into  $u(t)=O(|x|^{-N-sp}\log(|x|)^{\gamma}$ with $\gamma=1/(2-p_1)$.

\medskip

\noindent {\bf III. Case $p_c<p<p_1$.}

The main point is to prove the different upper estimate
$$
F_M(y)\le C_1^*\,|y|^{-sp/(2-p)} \,,
$$
where  $C_1^*>0$ is the constant introduced in Theorem \eqref{thm.barr.lr}, so the upper bound will not depend on $M$.
This is an important fact that will lead to the existence of the Very Singular Solution in Section \ref{sec.vss1}.
Note that the right-hand side of \eqref{decay.Fbarrier.lr} is invariant under the elliptic scaling transformation ${\overline{\mathcal  T}}_h $, which allows to change the mass of the self-similar solutions $F_M$. This explains why the tail estimate \sl does not depend \rm  on $M$. It is called a universal upper bound.

\noindent {\sl Review of the proof of existence.} We recall that the definition of $\mathcal K\subset X=L^1(\ren)$ involves the upper bound that we call $G(y;C_1)$ and is now taken from the supersolution construction of Theorem \ref{thm.barr.lr}. The rest of the proof is the same as before. \qed

The fixed point idea was already used in \cite{VPLP2020} but the end of proof is different.

{\subsection{Computed graphics}

Figures 2 and 3  show  the profiles of the self-similar fundamental solutions in the two ranges of $s$ and $p$. They were computed by numerically integrating the evolution equation starting with smooth initial data with compact support.

\begin{figure}[h!]
\includegraphics[width=\textwidth]{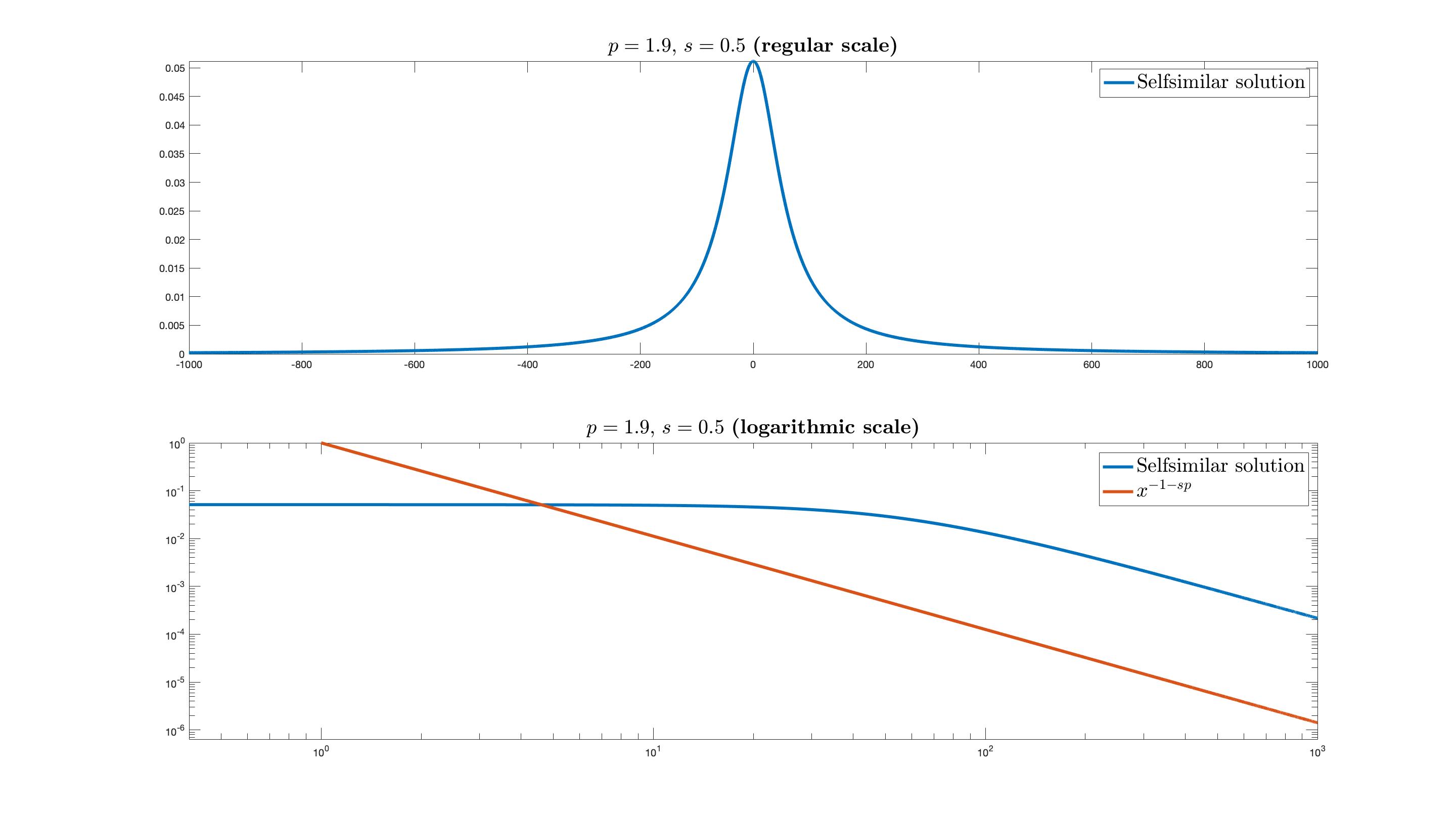}\caption{Self-similar fundamental solution for $s=0.5$ and $p=1.9>p_1$ }
\label{fig1:Self-sim graphics}
\end{figure}

\begin{figure}[h!]
\includegraphics[width=\textwidth]{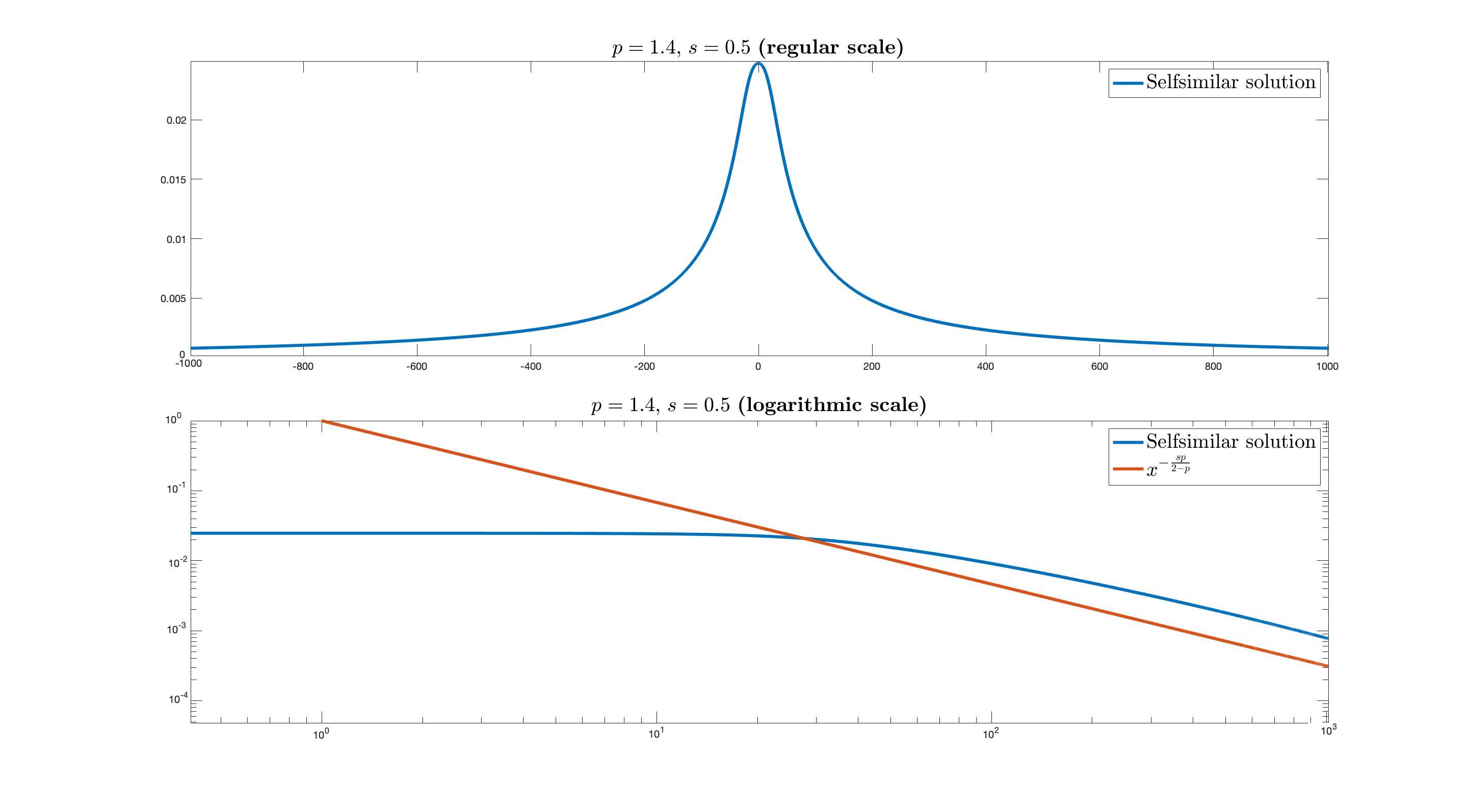}
\caption{Self-similar fundamental solution for $s=0.5$ and $p=1.4<p_1$.
 The second graphic in logarithmic scale shows clearly the decay with exponent $-sp/(2-p)$.
 }
\label{fig2:Self-sim graphics}
\end{figure}

 The shown computed examples correspond to dimension $N=1$ with $s=.5$, so that $p_c=4/3=1.33$ and $p_1=(\sqrt{15}-1)/2=1.436.$ The first figure is an example of approximation to self-similarity in the upper range $p_1<p<2$. We choose $p=1.9$. The second part of the plot displays the same function using the logarithmic scale. In this way the predicted decay is clearly shown, as the asymptotic line with slope $-(N+sp)\approx -1.95$. Figure 3 is an example of approximation to self-similarity in the lower range $p_c<p<p_1$. We use $N=1$, $s=.5$ and and $p=1.4$. 

The series of numerical experiments was performed by  F. del Teso using an explicit Euler finite difference scheme. He provides further details of the process as follows. The numerical discretization of the fractional $p$-Laplacian includes:  1) the corresponding weights of a quadrature for the singular integral, which are taken form \cite{dTEnJa18}, and 2) a Lipschitz regularization of the nonlinearity $\Phi_p$ (when $p<2$) to make the scheme stable and monotone. This is a well known trick in the context of explicit schemes for fast diffusion equations (see \cite{dTEnJa19}). The strong attraction properties of the asymptotic profiles, typical of nonlinear diffusion models, seems responsible for the fast convergence of the numerical method towards self-similarity. A rigorous analytical study of the numerical treatment of this equation has delicate points and is still in progress. \nc

\section{Asymptotic Behaviour}\label{sec.ab}

We establish here the asymptotic behaviour of finite mass solutions as $t\to\infty$, as stated in Theorem \ref{thm.ab1}.
We may assume that $M>0$ since the case $M<0$ can be reduced to positive mass by changing the sign of the solution.
The solutions are not necessarily nonnegative. We make a comment on $M=0$ below.

\subsection{Proof of the $L^1$ convergence} By scaling we may also assume that $M=1$.  The proof   relies on the previous results plus the existence of a strict Lyapunov functional,  that happens to be
\begin{equation}
J(u_1,u_2;t):=\int (u_1(x,t)-u_2(x,t))_+\,dx
\end{equation}
where $u_1$ and $u_2$ are two solutions with finite mass.

\begin{lemma}
Let $u_1$ and $u_2$ are two solutions with finite mass. Then, $J(u_1,u_2;t)$ is strictly decreasing in time
unless the solutions are ordered.
\end{lemma}

\noindent {\sl Proof. } By previous analysis, Section \ref{sec.diff},
we know that
\begin{equation}
\frac{d}{dt}J(u_1,u_2;t)=-\iint_D \left| |u_1(x:y,t)|^{p-2}u_1(x:y,t) -
|u_2(x:y,t)|^{p-2}u_2(x:y,t)\right|\,d\mu(x,y)\,,
\end{equation}
with notation as in \eqref{diff.est.1b}. In particular, the set $D\subset \re^{2N}$ contains the points where
$$
(u_1(x,t)-u_2(x,t))\,(u_1(y,t)-u_2(y,t))<0.
$$
Now, in order to $dJ/dt$ to vanish at a time $t_0>0$ we need $u_1(x:y,t)=u_2(x:y,t)|$ on $D$, i.\,e., $u_1(x,t)- u_2(x,t)=u_1(y,t)- u_2(y,t)$. But this is incompatible with the definition $D$, so $D$ must be empty, hence $u_1$ and $u_2$ must be ordered at time $t$. This implies that they have the same property for $t>t_0$. \qed

\medskip

\noindent {\sl  Proof of Theorem \ref{thm.ab1} continued.} (i) Let us assume that $u_0$ is bounded and compactly supported. It is convenient to  consider the $v$ version of both solutions, namely $v_1$ and $V_M=F_M$. Notice that $V_M$ is stationary and $v(y,\tau)$ is bounded above by a barrier function $G(y)$ as used above. We can show that $v(y,\tau+ n_k)$ converges strongly in $L^1(\ren)$, along a subsequence $n_k\to \infty$, towards a new solution $w_1$ of the $v$-equation. Under our assumptions,  $w_1$ must be a fundamental solution but maybe not self-similar. \nc

We know from the Lemma that $J(v_1, V_M;t)$ is strictly decreasing in time, unless $w_1(t)=V_M$ for all large $t$, in which case we are done. If this is not the case, we continue as follows. By monotonicity there is a limit
$$
\lim_{t\to\infty} J(u_1, U_M;t)=\lim_{\tau\to\infty} J(v_1, V_M;\tau)=C\ge 0\\.
$$
We want to prove that $C=0$, which implies our result. If the limit is not zero, we consider the evolution of the new solution $w_1$ together with $V_M$. We have
$$
 J(w_1, V_M;t_0)=\lim_{\tau \to \infty}J(v_1, V_M; t_0+\tau)=C\,,
$$
i.e., is constant for all $t_0>0$, which means that $w_1=V_M$  by equality of mass \normalcolor and the lemma. By uniqueness of the limit, we get convergence along the whole half line $t>0$ instead of a sequence of times.

(ii) For general data $u_0\in L^1(\ren)$, $M>0$, we use approximation. For $M<0$ just change signs.
Finally, in the case $M=0$ we just bound our solution from above and below by solutions of mess $\ve$ and $-\ve$ respectively, apply the Theorem and pass to the limit $\ve\to 0$.

\subsection{Proof of  convergence in $L^q$ norm}

When $1<q<\infty$, the $L^q$-convergence with rate to the fundamental solution, formula \eqref{lim.ab.q}, comes from interpolation between the same convergence in $L^1(\ren)$ plus the uniform boundedness of the $v$ orbit in $L^\infty(\ren)$ for large $\tau$, see Theorem
\eqref{L1-Linfty}. We leave the easy details to the reader.

\subsection{Proof of  convergence in $L^\infty$ norm}\label{ssec.asymp.Linfty}

\begin{theorem}\label{thm.ab.infty} Let $p_c<p<2.$ Let $u$ be a solution of Problem \eqref{frplap.eq}-\eqref{init.data} with initial data $u_0\in L^1(\ren)$, $u_0\ge 0$, and $u_0$ has compact support. Let $\int u_o\,dx=M>0$ and let $U_M$ be the fundamental solution with the same mass. Then,
\begin{equation}\label{lim.ab.inf}
\lim_{t\to\infty} t^{\alpha}\|u(t)-U_M(t)\|_\infty=0
\end{equation}
uniformly $x\in\ren$.
\end{theorem}

\noindent {\sl Proof.} (i) We can take $M=1$ to simplify, and also assume that $u_0$ has support inside the ball $B_1(0)$.
We begin by a simple remark. By the properties of the initial data and the properties of the self-similar solution, there exists a mass $M_u$ (maybe large) and a constant $T$ such that $u_0\le U_{M_u}(x,T)$, hence
$$
u(x,t)\le U_{M_u}(x,t+T)=(t+T)^{-\alpha}F_{M_u}(x,(t+T)^{-\beta})\,.
$$
In view of the fact that $F_{M_u}$ goes to zero at infinity, the uniform convergence \eqref{lim.ab.inf} holds on exterior sets of the form  $|x|\ge C\,t^\beta$ with $C$ large enough.

(ii) Let us argue in region  $\{ a\,t^{\beta}\le |x|\le C\,t^\beta\}$ where $F_1(r) $ is uniformly positive and bounded. By Aleksandrov's Principle, we conclude that $u(x,t)$ will have monotonicity properties along finite cones of the form $K_{\theta,R}(x_0)=K_{\theta}(x_0)\cap\{|x-x_0|\le R\}$ with vertex located any point $x_0$ away from $B_0$, say $|x_0|\ge 4$. The axis of the cone passes through the origin, there is a certain aperture angle $\theta$ that is uniform for $|x_0|\ge 4$ and the length is unbounded in the away direction, but only until a fraction of $|x_0|$ in the direction pointing to $x=0$.

Now, assume that we pass to rescaled variables $v(y,\tau)$. Similar monotonicity applies but the limitation is now
$|y_0|\ge 4/(t+1)^{\beta}$, leaving a rapidly decreasing hole near $y=0$. Suppose now that $v(y,\tau)$ does not converge to $F(y)$ for $|y|\ge  a$. There is a point $y_0$ at distance $r_0>a$ such that for an infinite sequence
$\tau_k\to\infty$ we have
$$
|v(y_0,\tau_k)-F(r_0)|\ge c>0\,.
$$
The difference can be by excess or defect. In the latter case we have
$$
v(y_0,\tau_k)\le F(r_0)-c\,.
$$
By using monotonicity along the cone in the outwards direction, we find a  finite cone $ K=K_{\theta,R}(x_0)$  such that
$$
v(y,\tau_k)\le F(r_0)-c \quad \mbox{for } \ y\in K.
$$
But when the cone is small we have by continuity $F(y)\ge F(r_0)-(c/2)$ in $K$, so that
$$
F(y)-v(y,\tau_k)\ge(c/2) \quad \mbox{for } \ y\in K.
$$
This is incompatible with the convergence in $L^1$ that has been established in Theorem \ref{thm.ab1}.

(iii) The argument for the excess case where $v(y_0,\tau_k)\ge F(r_0)+c$ is similar. Now we recall that $u$ and $v$ grow in a cone of directions pointing inwards to the origin.  \qed

(iv) The convergence near $x=0$ comes from the uniform H\"older continuity of the solutions obtained in Theorem \ref{them.reg.new} after a convenient rescaling. Indeed, for $t_1\gg 1$ we perform the rescaling
$$
{\mathcal T}_k u(x,t)= k^{\alpha} u(k^\beta x, k t),\quad k>0\,.
$$
with $k=t_1$. In this way we produce another solution $u_k ={\mathcal T}_k u$ of the equation with the same mass that takes at $t=1$ the values of the original $u$ after scaling factors. We check that this ${\mathcal T}_k u$ is very close to the Fundamental solution $F_M(x,1)$ at time 1, both for the $L^1$ norm and for the $L^\infty$ norm outside a small hole around the origin. Now we use the uniform H\"older continuity of \ref{them.reg.new} to conclude that the sup-convergence extends to the hole at $t=1$. Undoing the transformation we get the desired global sup-convergence up to the factor $t_\alpha$, as desired.

\subsection{Proof without additional assumptions.}

%

Assuming such result on  uniform  H\"older regularity for bounded solutions of equation \eqref{frplap.eq} we may follow like this to prove Theorem \ref{thm.ab.infty} without extra assumptions on $u_0$, just  nonnegative  and integrable data with positive mass $M$. We want to prove formula  \eqref{lim.ab.inf} uniformly in $\ren$. We return to the proof of the  previous step  and discover that the bounded sequence $v(y,\tau+ n_k)$ is locally relatively compact in the set of continuous functions in $\ren\times (\tau_1, \tau_2) $ thanks to the assumed results on H{\"o}lder continuity, once they are translated to the $v$-equation. Hence, it converges locally to the same limit as before, but now in uniform norm. In order to get global convergence we need to control the tails at infinity. We use the following argument: a sequence of space functions $v(\cdot, \tau)$ that is uniformly bounded near infinity in $L^1$ (thanks to the convergence to $F_M$) and is also uniformly H{\"o}lder continuous, it must also be  uniformly small in $L^\infty$. This implies that the previous uniform convergence was not only local but global in space. Using the correspondence \eqref{eq.rescflow},  we get the convergence of the $u(\cdot,t)$ with factor $t^{\alpha}$. This part of the theorem is proved.   \qed

\medskip

\noindent {\bf Comment.} We have proved that $t^\alpha u(x,t)$ converges to $F_M(y)$ uniformly in outer sets, but since both expressions go to zero as $|x|\to\infty$, the result does not say anything about the relative behaviour of the ``tails''. This is a delicate matter that we will study with great attention in upcoming sections and will only be complete with the so-called global two-sided bounds of Section \ref{sec.gharn}.

\section{Very singular solutions in the lower range}\label{sec.vss1}

 We open here a window towards a new topic that will be both a tool and an aim in itself. We consider the limits of fundamental solutions  with increasing mass in the lower exponent range $p_c<p<p_1$. We know that the set of fundamental solutions
 $$
 U_M(x,t)=t^{-\alpha}F_M(x\,t^{-\beta})
 $$
 is ordered with respect to the mass $M>0$. Moreover,  in this $p$  range there is an a priori estimate
in terms of the very singular supersolution, that we have found as a barrier in Theorem \ref{thm.barr.lr} and worked out in formula
\eqref{decay.Fbarrier.lr}
$$
F_M(y)\le C_1^*\,|y|^{-sp/(2-p)} \,,
$$
This family is ordered  with $M>0$. Passing to the monotone limit $M\to \infty$ in this expression,
we get
 $$
0<V_\infty(y):=\lim_{M\to\infty}F_M(y)\le C_1^*(N,s,p)\,|y|^{-sp/(2-p)}.
$$
 Now, if we apply the transformation $\overline{\mathcal T}_h$ to $F_M$ we find another fundamental solution profile with mass $M_h=h^{sp/(2-p)-N} M= h^{1/(2-p)\beta}M$, i.\,e.,
$$
\overline{\mathcal T}_h F_M(y)= F_{h^{1/(2-p)\beta}M}(y).
$$
Passing to the monotone limit $M\to \infty$ in the last expression, we get \ $  \mathcal T_h F_\infty(y)= F_{\infty}(y)$.
This immediately implies that
 $$
V_{\infty}(y)= C_\infty\,r^{-sp/(2-p)}
$$
for some  $ 0<C_\infty\le C_1^*.$ Note that the mass of this special {\sl limit solution} is infinite.

\begin{theorem}\label{thm.vss} Let $p_c<p<p_1$. There exists a constant $C_\infty(s,p,N) >0$ such that $ V_{\infty}(r)=C_\infty  \,r^{-sp/(2-p)}$ is a singular solution of the rescaled equation \eqref{eq.resc} with a non-integrable singularity at $r=0$. In original variables, we obtain
\begin{equation}\label{decay.vss.u}
U_\infty(x,t)=\,C_\infty  \,(t+T)^{1/(2-p)} { |x|^{-sp/(2-p)} }
\end{equation}
for any $T>0$. More precisely, this is  a classical solution of  equation \eqref{frplap.eq} for $x\ne 0$. It is called the Very Singular Solution, VSS.
\end{theorem}

\noindent{\sl Proof.} Going back to  Subsection \ref{ssec.barr.low}, we see that we cannot have the supersolution condition \eqref{ineq.super.2} for any $C_1<C_\infty$. This means that the constant $k$ in expression
$$
E(\widehat v)(r)= \left(\frac{C_1}{2-p}-kC_1^{p-1}\right)\,r^{-sp/(2-p)}
$$
must be positive, as already announced. The expression must vanish for $C_1=C_\infty$, and $\widehat v(y;C_1) $ is a solution for $C_1=C_1^*=C_\infty$.  Moreover, $E(\widehat v)(r)$ becomes positive for $C_1>C_1^*$ and negative for $0<C_1<C_1^*$.  We will use the small $C_1<C_\infty$ case  as a subsolution in the next Section. \qed

\noindent {\bf Remarks.} (1) Note $\partial_t U_\infty(x,t)>0$ for $x\ne 0$ and $t>0$, in agreement with the fact that $\mathcal L_{s,p} (r^{-sp/(2-p)})$ is always negative in our range.

\noindent (2) This type of singular solution is well known in the non-fractional setting $s=1$, and it exists for all $p_c<p<1$, cf. \cite{VazSmooth}.
The range we find here is clearly smaller than the standard case, and that needs an explanation that goes as follows: It is due to the strong influence of the nonlocal operator at long distances that dominates in the upper sublinear range. This will be precisely described below.

\noindent (3) The value of the constant for $s=1$ is given in \cite{VazSmooth}, formula 11.24, as
$$
C_\infty(1,p,N)= \left((p/(2-p))^{p-1}\beta^{-1}\right)^\frac1{2-p}\nc\,.
$$
Note that it goes to zero as $p\to p_c$.

\noindent (4)  The space decay $U_\infty(x,t)=O(|x|^{-sp/(2-p)})$ guarantees that $U_\infty(\cdot,t)$ is integrable at infinity precisely for $p>p_c$.

\noindent (5) $V_\infty$ is not a global weak solution including the origin, because $V_\infty $ is not locally integrable at $y=0$, even if $V_\infty^{p-1}\in L^1_{loc}$ . The origin is a non-removable singularity. Same comments for $U_\infty$.

\noindent (6) The same technique does not work for $p_1<p<2$ where the limit of the $U_M(x,t)$ will be shown to be infinite, see Subsection \ref{nonex.vss}.

\subsection{Nonlinear elliptic problem}
Here is an interesting elliptic consequence of the study of the VSS.

\begin{theorem}\label{cor.eigen}  Given $\lambda>0$, the function $F=A\,|y|^{-sp/(2-p)}>0$,  satisfies the nonlinear and singular elliptic  problem
\begin{equation}\label{eq.eigen}
\mathcal L_{s,p} F(y)+ \lambda F(y) =0 \qquad \mbox{for } \ y\ne 0,
\end{equation}
for the value
$$
A=c(N,s,p)\,\lambda^{-1/(2-p)}, \qquad c(N,s,p)=(2-p)^{-1/(p-2)}\,C_\infty .
$$
Note that \ $F(0)=+\infty$ with a non-removable singularity.
\end{theorem}

This in particular justifies the assertion made at the end of proof of Proposition  \ref{prop.barr1} that
$\mathcal L_{s,p}(r^{-sp/(2-p)})$ is negative for all $r>0$ if $p_c<p<p_1$.

\subsection{Asymptotic consequence for self-similar profiles}\label{ssec.pos1}
 We point out that $U_\infty(x,t)$ has a nonintegrable singularity at $x=0$, while on the contrary it is perfectly integrable at infinity. Since the fundamental solutions satisfy
 $$
 U_M(x,t)=M^{sp\beta}t^{-\alpha}F(x\,(M^{p-2}t)^{-\beta})\,,
 $$
 if we take $|x|=t=1$ and pass to the limit  in $M$ we get
 $$
 \lim_{M\to\infty } M^{sp\beta}F(M^{(2-p)\beta)})=U_\infty(1,1)=C_\infty  (N,s,p).
 $$
 In other words,
\begin{equation}\label{lim.profile.gr}
 \lim_{r\to\infty } r^{sp/(2-p)}F(r)=C_\infty  (N,s,p),
\end{equation}
 which gives the sharp behaviour of the self-similar profiles for $p_c<p<p_1$, as announced in formula \eqref{eq.decay1} of Theorem \ref{thm.beh}. By the way, the limit in \eqref{lim.profile.gr} is taken in an increasing way.

Notice also that this positive limit implies at the same time  that $F$ is strictly positive for all $r>0$ in a quantitative way, a clear manifestation of the infinite speed of propagation of the equation in this setting and range of parameters.

%
%
\section{Positivity of general solutions in the lower range}\label{sec.pos.lr}

We want to bound from below any given positive solution by a fundamental solution, maybe with small mass. We can do it after inserting a time delay.

\begin{theorem}\label{lower.par.est1.lr}
Let $0<s<1$ and $p_c<p<p_1$. Let $ u(x,t)$ be a solution of Problem  \eqref{frplap.eq} with nontrivial initial data $u_0(x)\ge 0$. Then, $u$ is positive  everywhere in $Q=\ren\times (0,\infty)$. More precisely, for every $t_1>0$ there exist a small mass $M_d$  such that for every $t\ge 2t_1>0$
\begin{equation}\label{low.est.lowrange1}
u(x,t)\ge U_{M_d}(x,t-t_1)\,.
\end{equation}
$M_d$ depends on the solution and on $t_1$. In any case $M_d\le M(u_0):=\int_{\ren} u_0\,dx$.
\end{theorem}

\noindent {\sl Proof.} (i)  Since both functions in the statement \eqref{low.est.lowrange1} are solutions we only have to prove the result for $t=2t_1$, i.e., we must have
\begin{equation}
u(x,2t_1)\ge U_{M_d}(x,t_1)
\end{equation}
everywhere in $\ren$. In order to get this, we will use a comparison argument in the time interval $(t_1,2t_1)$

\noindent (ii) We recall that the given solution is strictly positive, $u>c_1$ around the origin, say in $B_1(0)$, for some time $t=2t_1$. By monotonicity in time we get  $u>c_2>0$  for $t_1\le t\le 2t_1$ also in $B_1(0)$.
On the other hand, we have seen in Section \ref{sec.vss1} that for all $0<C<C_\infty$ the function $\widehat v(y;C)=C \, r^{-sp/(2-p)}$ is a subsolution of the $v$ equation for $r\ne 0$. Therefore,
$$
\widehat u(x,t;C)=C \, t^{1/(2-p)}|x|^{-sp/(2-p)}
$$
will be a subsolution of the $u$ equation for $x\ne 0$.

\noindent (iii) We will  use the known parabolic comparison principle for fractional equations, in the outer region, the cylinder $Q_e$ where $r>1$ and $t_1<t<2t_1$. There is no problem in comparing $u$ and $\widehat u$ inside that exterior cylinder, $Q_e$, since one of the functions is a solution, the other one a subsolution, so the fractional comparison principle says that (1) we also have to check the initial data outside ($u>0$ is enough since $\widehat u(x,0;C)=0$), and (2) we have to check the pointwise comparison in the inner cylinder with base the ball $B_1(0)$, which is more delicate, and will need a change of lower candidate since it has a singularity at $x=0$.

\noindent (iv) In order to do this last step we make the following modification.
We   cut at the level $\widehat v(y,\tau;C)=c_3$. Let the new candidate function will be
$$
\tilde v(y)=\min\{\widehat v(y;C), c_3\}.
$$
for some small $c_3>0$. We hope that this new comparison function will not create a problem. The important observation is that when we take $C$ very small, the cutting of $\widehat v(y)$  at level $c_3 $ will take place in a very small hole $B_\rho$, $\rho\ll 1$.
Indeed,
$$
\rho(C)=(C/c_3)^{(2-p)/sp}
$$
which tends to zero with $C$.

\begin{lemma} For $C$ small (depending on $c_3$) the modified $\tilde v(y)$ is still a subsolution of the $v$-equation outside the unit hole $B_1$.
 \end{lemma}

 \noindent {\sl Proof of the Lemma.} We recall that $\mathcal L_{s,p}\widehat v(y;C)=-kC^{p-1}\,r^{-(N+sp)}$. Also, for $p<p_1$ we have
$$
sp(p-1)<N(2-p), \qquad \frac{sp}{2-p}<N+sp.
$$
The key observation is that for $C$ small and $|y_0|= r_0>1$ the perturbation
$$
\mathcal L_{s,p}(\tilde v(y_0)-\widehat v(y_0;C))=
\int \frac{((\widehat v(y;C)^{p-1}-c_3^{p_1})_+}{|y-y_0|^{N+sp}}\,dy
$$
can be estimated as  less than $\ve(C)r_0^{-(N+sp)}$, where
$$
\ve(C)=\int_{B_\rho(0)} \widehat v(y;C)^{p-1}\,dy=c C^{p-1}\int_{B_\rho(0)}r^{N-1-\frac{sp(p-1)}{2-p}}dr= c C^{p-1}\rho^{\gamma}
$$
with $\gamma=N- \frac{sp(p-1)}{2-p}>0$. Recalling the value of $\rho(C)$ we  get
$$
\ve(C)=cC^{p-1}(C/c_3)^{\gamma (2-p)/sp}=cc_3^{-\mu} C_1^{p-1+\mu}
$$
The last exponent is larger than $p-1$  for $p<p_1$. Therefore,
$$
\mathcal L_{s,p}(\tilde v(y_0)  \le (-k+ c_4C^{\mu})C^{p-1}r_0^{-(N+sp)}\le -(k/2)C^{p-1}r_0^{-(N+sp)}<0.
$$
Going to the subsolution condition derived from \eqref{ineq.super.2}, we conclude that  $\tilde v(y)$ is a subsolution if $C/c_3$ is small.
\qed

\noindent (v) We now apply the fractional parabolic comparison to our delayed solution $u_1(x,t)=u(x,t+t_1)$ and to
$$
\tilde u(x,t)=t^{-\alpha}\tilde v(xt^{-\beta})=\min\{C \, t^{1/(2-p)}|x|^{-sp/(2-p)}, c_3\}.
$$
The flat region is now the ball $B_{\rho(t)}(0)$ with
$$
\rho(t)=(C/c_3)^{(2-p)/sp}t^{1/sp},
$$
which is less than 1 for a long time if $C/c_3$ is small.
In the outer domain $|x|\ge 1$ we have for all this  time $\tilde u(x,t)<c_3$, so that
$\mathcal L_{s,p}\tilde u$ is calculated at a fixed time $t$ as in the previous lemma but taking $Ct^{1/(2-p)}$ instead of $C$.
Therefore, for $t<1$ we can get for $C$ small enough
$$
\tilde u_t+\mathcal L_{s,p}\tilde u <\frac1{2-p}Ct^{(p-1)/(2-p)}|x|^{-sp/(2-p)}-(k/2)C^{p-1}t^{(p-1)/(2-p)}r_0^{-(N+sp)}<0
$$
We may now apply the comparion theorem to $u_1(x,t)$ and $\tilde u(x,t)$  and conclude that for $t=t_1$ and $|x |\ge 1$ we have
$$
u_1(x,t)=u(x,2t_1)\ge \tilde u(x,t_1)
$$
Since $u(x,2t_1)$ is bounded below by a constant in $B_1(0)$, we have succeeded in squeezing a fundamental solution with small mass and some delay everywhere below  $u(x,2t_1)$. For $t\le 2t_1$
we use plain comparison. \nc

\noindent (vi)   The fact that $M_d\le M(u_0)$ follows from  the  asymptotic results of Section \ref{sec.ab}. Just check the values near $x=0$ for large $t$. \qed

\medskip

\begin{corollary}\label{tailasym.lt} Under the above conditions, for every $\ve>0$ there exists $R>0$ such that whenever  $t\ge 2t_1>0$ and $|x|\ge R\,(t-t_1)^\beta$, then
$$
u(x,t)\ge C_\infty(1-\ve)\,t^{1/(2-p)}\,|x|^{-sp/(2-p)}(1-(t_1/t))^{1/(2-p)}.
$$
In particular, if $t\gg t_1$ we get
$$
u(x,t)\ge C_\infty(1-2\ve)\,t^{1/(2-p)}\,|x|^{-sp/(2-p)}.
$$
Note that $R$ depends on $u_0$ and $t_1$.
\end{corollary}

\medskip

\noindent {\bf Remark.} This minimum behavior of nonnegative solutions when $|x|\rightarrow \infty$ is sharp as it coincides with the rates of the VSS both in space and time.

A comparison result without a time delay is the following.

\begin{corollary} Under the above conditions we get for all $t\ge 2t_1$ and $x\in\ren$
\begin{equation}\label{low.est.mass}
u(x,t)\ge U_{M_d}(x,t)\,(1-(t_1/t))^{1/(2-p)} \,.
\end{equation}
\end{corollary}

\noindent  {\sl Proof.}  We calculate
$$
U_M(x,t-t_1)=(t-t_1)^{-\alpha}F_M(|x|\,(t-t_1)^{-\beta})
$$
We use the monotonicity formula: for all $r>0$ and $\lambda>1$ we have $F(\lambda r)\ge \lambda^{-\gamma} F(r)$ with $\gamma=sp/(2-p)$, to derive
$$
F_M(|x|\,(t-t_1)^{-\beta})\ge F_M(|x|\,t^{-\beta})\,(t/(t-t_1))^{-\gamma \beta}
$$
Therefore,
$$
U_M(x,t-t_1)\ge U_M(x,t)\,(t/(t-t_1))^{\alpha-\beta \gamma}=U_M(x,t)\,(t/(t-t_1))^{-1/(2-p)}.
$$
This and the previous theorem proves \eqref{low.est.mass}. \qed \nc

%
%
\section{Analysis of positivity in the upper fast range}\label{sec.pos.tail.ur}

The analysis of the range $p_1<p<2$ will allow us to obtain the minimum  behavior of nonnegative solutions when $|x|\rightarrow \infty$, more precisely their rate of space decay, for small times $t>0$. This will imply the precise decay rate of the profile of the fundamental solution. Our goal is to obtain a lower bound that matches the spatial behaviour of the upper barrier, as established in Section \ref{sec.barr2}.

\begin{theorem}\label{low.par.est1}
Let $0<s<1$ and $p_1<p<2$. Let $\underline u(x,t)$ be a solution of Problem \eqref{frplap.eq} with initial data $u_0(x)\ge 0$ such that $u_0(x)\ge 2$ in the ball $B_3(0)$. Then there is a time
$t_1>0$ and a constant $c>0$ such that
\begin{equation}
u(x,t)\ge c\,t\,|x|^{-(N+sp)}
\end{equation}
if $|x|\ge 2$  and $0<t<t_1$.
\end{theorem}

We will use a comparison argument based on the following construction.

\begin{lemma} There is a smooth, positive and radial function $G_1(r) $ in $\ren$ such that

(i) $G_1(r)\le 1$ everywhere, and $G_1(r)=cr^{-(N+sp)}$ for all $r>2$

(ii) ${\mathcal L}_{s,p}G_1$ is bounded and ${\mathcal L}_{s,p}G_1(r)\approx - r^{-(N+sp)}$ for all $r\ge R>2$.
\end{lemma}

\noindent {\sl Proof of the Lemma.} We define $G_1$ by  specifying it in three different regions.  For $r\le 1$ we put $G_1=1$. For $r>3/2$ we put $G_1(r)=cr^{-(N+sp)}$ as indicated, with a small constant $0<c<c_0$ that will change in the application, so we must pay attention to it. In the intermediate region we choose a smooth and radially decreasing function that matches the values at $r=3/2$ and $r=2$ with $C^2$ agreement.

It is then easy from the theory to prove that ${\mathcal L}_{s,p}G_1$ is bounded on any ball, so we only have to worry about the behaviour at infinity, more precisely for $r\gg 2$. In order to analyze that situation  we point out that, according to formula \eqref{frplap.op}, \ ${\mathcal L}_{s,p}G(x)$ is an integral with contributions from the variable $y$ in different regions.

Next, we show that the contribution from the ball $B_{3/2}(0)$ is the largest one. Indeed, we have for $|x|=r>2$
$$
-I_1(r)=\int_{B_{3/2}} \frac{1-G(x)}{|x-y|^{N+sp}}\,dy \le (1-\ve)|B_1| \ (r/2)^{-(N+sp)},
$$
that does not depend on the small parameter $c$ appearing in the tail. The other contributions depend on $c$ and can be made small with respect to $I_1(r)$ for all $r>2$  precisely if $p_1<p<2$,  see details in Section \ref{sec.barr2}. \qed

\medskip

\noindent {\sl Proof of the Theorem.} (i) We modify function $G_1$ to introduce a linear dependence on time in the outer region. We take a smooth cutoff function $\eta$ lying between 0 and 1 such that
$\eta(x)=1$ for $|x|\le 1$ and $\eta(x)=0$ for $|x|\ge 2$, and put
\begin{equation}\label{subsol1}
U(x,t)=\eta(x)\,G_1(x)+ (1-\eta(x))ct \,r^{-(N+sp)}.
\end{equation}

\noindent (ii) We want to prove that this function satisfies the subsolution condition
\begin{equation}\label{subsol2}
U_t+{\mathcal L}_{s,p} U <0
\end{equation}
in an outer region $\{r>R\}$ and for an interval of times $0<t<t_*(c)$ if $c$ is small enough. Now for $R\ge 2$ we have
$$
U_t= cr^{-(N+sp)}>0.
$$
On the other hand, the proof of the Lemma shows that in that region
$$
{\mathcal L}_{s,p} U \le -Cr^{-(N+sp)},
$$
as long as we can disregard the contributions from outside $B_{3/2}$, and this is true if $tc$ is mall enough. The conclusion \eqref{subsol2} follows.

\noindent (iii) Next, we proceed with the comparison step between $u$ and $U$ in a space-time domain of the form
$Q=\{(x,t): \ |x|\ge 2, \ 0<t<t_1\}$.  By comparison we may consider some smaller initial data $u_0$, such that $0\le u_0(x)\le 2$ and $u_0(x)=  2$ in the ball of radius 3 and $u_0$. Moreover, $u_0$ is smooth. By previous results of this paper (cf. the continuity arguments of \eqref{ssec.pos}) we know that $u(x,t)\ge  2-\ve$ in a ball of radius $2<R<3$ for all small times $0<t<t_0$.
And we know that $U(x,t)\le 1$ at all points (while $t$ is small).

Summing up, we already have the necessary inequalities for the equation inside that domain. As for initial conditions we know that $U(x,0)=0$ for all $|x|\ge 2$, while $u_0\ge 0$ everywhere. Comparison in the inner cylinder,  $|x|\le 2$ and $0<t<t_1$, has been also established.

\noindent (iv)  Now we only need to integrate by parts the difference of the two equations with multiplier $(U-u)_+$ to get the conclusion that $(u-U)_+$ must be zero a.e. in $Q$. Note that both functions belong to $L^2(\ren)\le L^\infty(\ren)$ uniformly in $t$. Since $U\le u$ in the set $\Omega=\{|x|\ge 3\}$ for $0<t<t_0$, we get for all those times
$$
\frac{d}{dt}\int_{\Omega} (U-u)_+^2\,dx= 2\int_{\Omega} (U-u)_+(U_t-u_t)\,dx=
2\int_{\ren} (U-u)_+(U_t-u_t)\,dx=I.
$$
But that integral is easily estimated
$$
I=-2\int_{\ren} ({\mathcal L}_{s,p}U-{\mathcal L}_{s,p}u)\,(U-u)_+\,dx\le 0
$$
by $T$-accretivity (better do the direct computation, see above the computation of the evolution of the $L^2$  norm of the difference of two solutions). Since  $(U-u)_+=0$ in $\Omega$ for $t=0$, we get the desired conclusion:
$$
u(x,t)\ge U(x,t)\ge ct\,r^{-(N+sp)}
 $$
 if  $r\ge 2$ and  $0<t<t_1$. \qed


\subsection{Application to the self-similar solution}

We consider the fundamental solution after a time displacement:
$$
u_1(x,t)= (t+1)^{-\alpha} F_M(|x|\,(t+1)^{-\beta}),
$$
that satisfies the assumptions of Theorem \ref{low.par.est1} if $M>0$ is large enough. We conclude that
\begin{equation}\label{asm.est.upr}
F_M(r)\ge c_1\,r^{-(N+sp)} \quad \mbox{for all large } \ r.
\end{equation}
By scaling, the same is true for $M=1$ with a different constant. Together with the upper bound from Theorem \ref{thm.barr.lr}, we get assertion \eqref{eq.decay2} of Theorem \ref{thm.beh} which was still missing to complete the proof.

For possible future reference, let us state the tail behaviour of the fundamental solution $U_M(x,t)=t^{-\alpha}F_M(|x|\,t^{-\beta})$. Let us choose $M>0$.

 \begin{corollary} \label{cor.ob.fs} On every outer region of the form $\{(x,t): \ |x|\ge C t^{\beta}, \ C>0\}$ we have constants $0<C_1<C_2$ such that
\begin{equation}\label{sigma}
C_1 \,M^{\sigma}|x|^{-(N+sp)}t^{sp\beta}\le U_M(x,t)\le C_2\, M^{\sigma\beta}|x|^{-(N+sp)}t^{sp\beta},
\end{equation}
where $\sigma=sp-(2-p)(N+sp)$, which is positive in this range.
 \end{corollary}

Actually, $\sigma= sp(p-1)-N(2-p)$ is positive precisely for $p>p_1$ and goes to $0$ when $p\to p_1$ for $s,N$ fixed. This shows good agreement with the behaviour of the lower range $p_c<p<p_1$,  where $M$ does not appear in the bounds, see \eqref{tailasym.lt} and \eqref{ghi} below.

\medskip

\noindent {\bf Remark.} Positivity estimates related to the ones in this section have been obtained for the fractional porous medium equation in \cite{VazBar2014, StanVa14, VolVa14}. Other forms of positivity estimates  were developed in \cite{BV06} for the Fast Diffusion Equation, and in the fractional case in \cite{BV14}.

\nc

\subsection{Non-existence of the very singular solution}\label{nonex.vss}

In view of estimate \eqref{asm.est.upr} we get for fixed $y\ne 0$ and all large enough $M$:
$$
F_M(y)=M^{sp\beta}F(y\,M^{(2-p)\beta})\ge c_1 |y|^{-(N+sp)} M^{\sigma\beta}\,,
$$
with $\sigma>0$. Therefore, for every $y\in \ren$
$$
\lim_{M\to \infty} F_M(y)=+\infty,
$$
so that the Very Singular Solution like in Section \ref{sec.vss1} does not exist in the exponent range $p_1<p<2$.

\subsection{Bounds for the critical exponent $p_1$}\label{sec.p_1}

We recall the upper bound from Theorem \ref{thm.barrcomp_1}

\begin{theorem}[Barrier comparison] \label{thmm}
Let $u$ be a solution with nonnegative data $u_0$ such that
$\|u_0\|_\infty\le L$ and $\|u_0\|_1\le M$. There exists a constant $C(L;M)>0$ such that whenever
\begin{equation}
u_0(y)\le C_1 \,r^{-(N+sp_1)} \log (r)\quad \mbox{ for } \ r>2\,.
\end{equation}
with  $C_1\ge C(L;M)$, then the corresponding $v$-solution satisfies
\begin{equation}\label{barrcomp.v}
v(y,\tau)\le \widehat v(y;C_1) \qquad \forall y\in \ren, \tau>0\,.
\end{equation}
In other words,
\begin{equation}
u(x,t)\le\widehat u(x,t;C_1):= (1+t)^{-\alpha}\widehat v(|x|\,(1+t)^{-\beta};C_1) \qquad \forall x\in \ren, t>0.
\end{equation}
    \end{theorem}

  This decay applies in particular to the FS, $F_1(y)\le C\,r^{-(N+sp_1)} \log (r)$ \ for  \ $r>2$. But the lower bound that we can obtain is based on the argument of
Theorem \ref{low.par.est1} and only gives an estimate of the form $F(y)\ge C\, |y|^{-(N+sp_1)}$.

  Open problem: find a sharper lower bound. \nc

\section{Two-sided global bounds. Global Harnack }\label{sec.gharn}

The uniform convergence proved in the asymptotic section implies that
$$u(x,t)/U_M(x,t)\to 1 \qquad \mbox{as } t\to \infty
$$
uniformly on sets of the form $\{|x|\le ct^{\beta}\}$. But it does not say anything about the relative error on the far away region, i.e. for the so-called tail behaviour.

We can contribute to that issue by using the positivity analysis of Sections \eqref{sec.pos.lr} and \ref{sec.pos.tail.ur} together with previous results on upper bounds. We will obtain a two-sided global estimate, assuming that the initial data are bounded, nonnegative and compactly supported. The result applies to all positive times and says that the relative quotient $u(x,t)/U_M(x,t)$ stays bounded for $t\ge \tau>0$.

This kind of two-sided bound by the fundamental solution is usually called a Global Harnack Inequality and is frequent in nonlinear diffusion problems with fast diffusion. See  applications to the fast diffusion equation in \cite{Vascppme, CrVa03, BV06}, and a very recent one in \cite{Simonov}. It is not true for equations with slow diffusion and free boundaries. There are a number of references for fractional parabolic equations like  \cite{BSV15,  BSV17}, even in the so-called slow range.  We will not mention the large literature on elliptic problems or problems in bounded domains.  

The technical details and final result are a bit different in both $p$-subranges and will be studied separately.

\subsection{Two-sided global bounds for $p_c<p<p_1$ }\label{sec.gharn1}

Once the desired lower bound was established in Section \eqref{sec.pos.lr}, we address the existence of a similar upper bound for general solutions  for the same range $p_c<p<p_1$. Actually, we need some conditions on the initial data for the result to be true.

\begin{theorem}\label{thm.GH} Let $0<s<1$, $p_c<p<p_1$, and let $u$ be the semigroup solution corresponding to initial data $u_0\in L^\infty(\ren)$, $u_0\ge 0$, $u_0\not\equiv 0$, supported in a ball of radius $R$. Then, there exist a mass $M_u>0$ and a constant $T>0$ such that
\begin{equation}
u(x,t)\le U_{M_u}(x,t+T) \qquad \mbox{ for all \ } \ x\in \ren, \ t\ge 0.
\end{equation}
The constants  $M_u$ and $T$ may depend on $u_0$. Moreover, if $M(u_0)=\int_{\ren} u_0\,dx$, then $M(u_0)\le M_u$.
\end{theorem}

\medskip

\noindent {\sl Proof}. (i) We  analyze separately the region near infinity and then the interior core.  In a region near infinity we use the comparison with the upper barrier estimate of Theorem \ref{thm.barr.lr}, formula \eqref{decay.u1.b},
\begin{equation}\label{decay.u1.b2}
u(x,t)\le C_\infty\,\frac{\,(t+T_1)^{1/(2-p)}} { |x|^{sp/(2-p)} }\,,
\end{equation}
where we use a delayed Very Singular profile. In order to insert the fundamental solution we see
for every $\ve>0$ there exists $R>0$ such that
$$
F_1(y)\ge C_\infty(1-\ve)\,r^{-sp/(2-p)}=(1-\ve)V(r)
$$
whenever $|y|=r\ge R$. Taking  $T=2T_1$, this means that for $|x|\ge  M_u^{-(p-2)\beta}R\,(t+T)^\beta$ we have
$$
{U_{M_u}(x,t+T)}\ge (1-\ve) {V(x,t+T)}.
$$
We then get
$$
\frac{U_{M_u}(x,t+T)}{V(x,t+T_1)}\ge (1-\ve)\left(\frac{t+T}{t+T_1}\right)^{1/(2-p)}\ge 1.
$$
It follows that $u(x,t)\le U_{M_u}(x,t+T)$ in the outer region $|x|\ge \widehat R= R M_u^{-(p-2)\beta}(t+T)^\beta$ and $t$ near zero.


\noindent (ii) In the interior region  we use the comparison of $u(x,t)$ and $U_{M_d}(x,t+T)$
on the inner region that follows from the condition
$$
\|u_0\|_\infty\le M_u^{sp\beta}(t+T)^{-\alpha}F_1(R)
$$
for $|x|\ge \widehat R(t)$. We need to check it at $t=0$, where we need
$$
\|u_0\|_\infty\le c M_u^{sp\beta}T^{-\alpha},
$$
that is obtained for $M_u$ large or $T$ small.

\noindent (iii)  In view of the  asymptotic results of  Section \ref{sec.ab} we have $M_d<M(u_0)<M_u$. Just check the values at $x=0$ for large $t$. \qed

The upper and lower bounds imply  a two-sided global estimate, assuming that the initial data are bounded, nonnegative and compactly supported. The result applies to all positive times and says that the relative quotient $u(x,t)/U_M(x,t)$ stays bounded for $t\ge t_0>0$.

\begin{corollary} {\bf (Global 2-sided estimate)\rm} For every solution $u(x,t) $ as in the previous theorem, there are finite constants $k_2>k_1>0$ such that the inequalities
\begin{equation}
k_1\le  \frac{ u(x,t)}{U_M(x,t)}\le k_2
\end{equation}
hold uniformly in $x\in \ren$ and $t\ge 1$.
 \end{corollary}

\noindent {\sl Proof}. The lower bound uses formula \eqref{low.est.mass} and the elementary estimate
$$
U_{\lambda M}(x,t)\ge \lambda^{sp\beta} U_{ M}(x,t)
$$
if $\lambda <1$. The inequality if reversed if $\lambda<1$. This together with
$$
U_M(x,t+T)\le U_M(x,t)\,((t+T)/t)^{1/(2-p)}
$$
proves the upper bound. \qed

Let us write in a clear way our conclusion about the very precise size of the spatial tails.

\begin{corollary}\label{two-sided.lr} For every solution $u(x,t) $ as in the previous theorem, and for every $t\ge 1$ and for every $\ve>0$ there is an $R(\ve)>0$ such that
\begin{equation}\label{ghi}
C_0-\ve< u(x,t)\,|x|^{sp/(2-p)}\,t^{-1/(2-p)}<C_0+\ve
\end{equation}
on the outer set \ $|x|\ge R(\ve) t^{\beta}$ if $t$ is large.
\end{corollary}

\subsection{Two-sided global bounds for $p_1<p<p_2$ }\label{sec.gharn2}

We are going to obtain similar upper and lower bounds in terms of fundamental solutions in the upper $p$ subrange.

\begin{theorem}\label{thm.GH} Let $u$ the semigroup solution corresponding to initial data $u_0\in L^\infty(\ren)$, $u_0\ge 0$, $u_0\ne 0$, supported in a ball of radius $R$. For every $\tau>0$ there exist constants $M_d, M_u>0$ and a delay $ c_2>0$ such that
\begin{equation}
U_{M_d}(x,t)\le u(x,t)\le U_{M_u}(x,t+c_2) \qquad \mbox{ for all \ } \ x\in \ren, \ t\ge \tau.
\end{equation}
The constants  $M_d, M_u$, and $c_2$ may depend on $\tau$. Moreover, if $M(u_0)=\int_{\ren} u_0\,dx$, then $M_d\le M(u_0)\le M_u$.
\end{theorem}

\medskip

\noindent {\sl Proof}. (i) Let us begin by the upper bound that is an easy consequence of the barrier estimate of Section \ref{sec.barr2}, in particular Theorem \ref{thm.barrcomp}, together with the lower estimate in \eqref{sigma} for the Fundamental solution. Hence, the upper bound  \eqref{decay.upur} is comparable with the self-similar profile $F_1$, hence smaller than $F_{M_u}$ for some $M_u>1$. This estimate is valid even for $\tau=0$, with $c_2=1$ and $M_u$ large enough. It holds then for every $t>0$ by comparison.

\noindent (ii) For the lower bound we need to take $\tau>0$ and use the results of Section \ref{sec.pos.tail.ur}. By comparison, translations and rescaling we my assume that $u_0$ is as in Theorem \ref{low.par.est1}. We also assume that is radially decreasing. Therefore, given a time  $\tau>0$ small enough we have the estimate \ $ u(x,\tau)\ge c\,\tau\,|x|^{-(N+sp)} $ \
for all $|x|\ge R > 2$. On the other hand,
$$
U_{M_d}(x,c_1)= c_1^{-\alpha}F_{M_d}(|x|\,c_1^{-\beta})\le C \min\{ c_1^{-\alpha}\, M_d^{sp\beta}, \,M_d^{\sigma} c_1^{sp\beta}|x|^{-(N+sp)}\},
$$
for every $|x|\ge 0$. We have used formula \eqref{sigma}. We conclude that for given $c_1>0$ there exists $M_u$ small enough such that
$$
u(x,\tau)\ge U_{M_u}(x,c_1).
$$
We may put now $\tau=c_1$. By comparison the lower bound is true for all $t\ge \tau$.

(iii)  In view of the previous asymptotic results we have $M_d<M(u_0)<M_u$. Just check the values at $x=0$ for large $t$. \qed


The result applies to all positive times and says that the relative quotient $u(x,t)/U_M(x,t)$ stays bounded  for $t\ge t_0>0$.

\begin{corollary}\label{twosided.ur} {\bf (Global 2-sided estimate)\rm} For every solution $u(x,t) $ as in the previous theorem, there are finite constants $k_1,k_2>0$ with $k_1<1<k_2$ such that the inequalities
\begin{equation}
k_1\le  \frac{ u(x,t)}{U_M(x,t)}\le k_2
\end{equation}
hold uniformly in $x\in \ren$ and $t\ge 1$.
 \end{corollary}
Let us write in a clear way our conclusion about the size of the spatial tails.

\begin{corollary} For every solution $u(x,t) $ and as in the previous theorem, and for every $t\ge 1$ there are positive constants $C_1,C_2>0$ such that
\begin{equation}\label{ghi2}
C_1< u(x,t)\,|x|^{N+sp}\,t^{-sp\beta}<C_2
\end{equation}
on the outer set  $|x|\ge C \, t^{\beta}$, $C>0$ large. Here, $C_1$ and $C_2>0$ may be quite different.
\end{corollary}
\nc

Note that the rates in formulas \eqref{ghi} and \eqref{ghi2} are quite different. They agree in the limit $p=p_1$ (where they are not really correct because of a logarithmic correction, so typical of critical cases). Moreover, on the expanding spheres $|x|=C \, t^{\beta}$ both formulas predict a time decay like $O(t^{-\alpha})$.

\section{Limit cases}\label{sec.limitcases}

In the paper we have considered  all fractional exponents in the range $0<s<1$ and nonlinear exponents $p_c<p<2$. The limit cases are interesting
as examples of continuity with the dependence on parameters. We will make here a brief sketch of relevant facts.

\medskip

\noindent $\bullet$ {\sl Limit $p\to 2$.} The limit of the $(s,p)$ semigroup as $p\to 2$ for fixed $s$ offers only minor difficulties.   Also the passage to the limit in the self-similar solutions gives the well-known profiles of the fractional linear heat equation. These profiles decay like $O(|x|^{-(N+2s)} )$ as $|x|\to\infty$, cf. \cite{Blumenthal-Getoor}, see also \cite{BSV17} and its references. The linear self-similar solutions are also asymptotic attractors, as proved in \cite{Va18}, where convergence rates are obtained. The  limit  $p\to 2$ can also be checked computationally with minor difficulty.

\nc
\noindent $\bullet$  {\sl Limit $s\to 1$.} It is well known by experts that the operator ${\mathcal L}_{s,p}$ must be conveniently renormalized by a constant factor including the factor $1-s$,  cf. \cite{BBM02, DTGCV, IshiiN, MazRT}, in order to converge to the standard $p$-Laplacian as $s\to 1$. It is then rather easy to prove that, as $s\to 1$ for fixed $p> 2$, we obtain the semigroup corresponding to the standard $p$-Laplace operator, already mentioned in the introduction.\nc

In particular, we can pass to the limit in the self-similar solutions of Theorem \ref{thm.ssfs} and check that the self-similar profile $F_{s,p}(r)$ converges to the  profile $F_{1,p}(r)$, given by the well-known Barenblatt profile \
$F_{1,p}(r)=\left(C+ kr^{\frac{p}{p-1}}\right)^{\frac{p-1}{2-p}}, $
already mentioned. The decay exponents $\alpha(s,p)$ and $\beta(s,p)$ also converge. Notice that for $0<s<1$ the self-similar profiles $F_{s,p}(r)$ are positive with tails at infinity, and the limit Barenblatt profile $F_{1,p}(r)$  has infinite tails too, but the rates obtained in the limit of the fractional decay are  different in the upper subrange $p_1\le p<2 $.\nc

Full details should be provided elsewhere. A similar study of continuous dependence with respect to parameters has been done in full detail in the case of the Fractional Porous Medium Equation in \cite{DPQRV2} (it includes fast diffusion).

\section{Comparison of two nonlinear fractional diffusion operators}\label{sec.comp.two}

We have already referred to the motivation received from the study of the standard $p$-Laplacian equation and also from the fractional linear heat equation. Both serve as limit cases of our equation and provide a number of results and tools that have been used in building this theory, as well as the theory of companion paper \cite{VPLP2020} for $p>2$. When studying fast diffusion asymptotics we need to mention the monograph \cite{VazSmooth} as a convenient previous reference.

We will devote this section to compare our paper with the theory developed in paper \cite{VazBar2014} for the so-called ``fractional porous medium equation'' (FPME)
\begin{equation}  \label{eq1.fpme}
\partial_tu+(-\Delta)^{s}(u^m)=0,  \qquad  x\in\mathbb{R}^N,\; t>0\,,
\end{equation}
that was studied with fractional exponent $s\in (0,1)$, and  $m>0$. That paper contains the proof of  existence, uniqueness and main properties of the fundamental solutions for the FPME, using the previous study made in \cite{PQRV1}. Like here, they are self-similar functions of the form \ $U(x,t)= t^{-\alpha} F(|x|\,t^{-\beta})$ \ with suitable $\alpha$ and $\beta$.
The admissible range is $m>m_c=(N-2s)/N$. The paper goes on to  prove that the asymptotic behaviour of general solutions is represented by such special solutions in the whole admissible range. Very singular solutions are also constructed in a more reduced range $m_c<m_1<1$. Since  paper \cite{VazBar2014} develops a theory with many topics in common with the present one, a comparison should be illustrative for the interested reader. As it often happens with nonlinear models, there are some results and methods that can be adapted after some effort, while other aspects demand a completely new approach.  This is a well-known fact when comparing the theories of the standard porous medium equation with the standard $p$-Laplacian equation.

In the comparison between the fractional models, a key tool in the proof of existence of fundamental solutions for the FPME  was the introduction of the so-called ``dual equation'';  much of the technical work was based on that tool that allows for very convenient maximum principles to be applied. However, such a tool is not available here, so we have been forced to follow a completely different approach to the existence of self-similar fundamental solutions, which implied developing new tools or adapting recent results of other works. A  key step was the construction of sharp upper barriers, a nontrivial task done in Sections \eqref{sec.barr1} and \eqref{sec.barr1}. This process took a certain lapse of time. The extra complication appears in many small details. Thus, the new critical exponent $p_1$ is much less evident than $m_1$.

Another marked difference between both equations lies in the difficulty of making explicit computations for the $p$-Laplacian in order to find explicit solutions. On the contrary, in the PME   \eqref{eq1.fpme} there is a linear operator $(-\Delta)^{s} $ that is easier to manage, see for instance \cite{dyda} for explicit computations with the $s$-Laplacian. This  allows to find or check explicit solutions, or sub- and supersolutions. For instance,  checking the VSS formulas is an easy task for equation   \eqref{eq1.fpme}, done in \cite{VazBar2014} and \cite{VolVa15} using Fourier analysis. In the case of the fractional $p$-Laplacian, the  VSS is a rare example of an explicit solution, and even in that case the existence is not obvious, and the value of the universal constant $C_\infty$ is unknown (to our knowledge). Actually,  we needed to establish first the fundamental solutions and derive some convenient and highly nontrivial estimates in order to prove existence of the VSS. As we have seen, the existence of these special solutions in turn implies very strong Global Harnack Inequalities in the corresponding $p$-range $p_c<p<p_1$.

We will continue the comparative analysis in the next section for lower values of $p$.

\section{On the very fast diffusion range}\label{sec.vfd}

This section contains a number of contributions on the FPLEE in the range $1<p<p_c$ that serve as a complement to the information of previous sections. We are mainly interested in a topic: the existence and properties of Very Singular Solutions. Then, we are also interested  in the nonexistence of fundamental solutions, the failure of the smoothing effect of Theorem \ref{L1-Linfty}, and the occurrence of finite mass extinction for a class of initial data. We are not aiming at a more complete theory, since appears to be rich enough, as can be inferred from the theory of the standard porous medium and $p$-Laplacian equations done in the fast diffusion range in \cite{VazSmooth}.

The basic facts of the theory contained in the first 5 subsections of Section \ref{sec.basic} are still valid with some small changes.  Thus, for every $u_0\in L^q(\ren)$, $1\le q< \infty$ there exists a unique semigroup solution with the expected basic properties. Only positivity for all times is not ensured, because solutions exist that vanish identically after a finite time.  Moreover, we do not have the $L^1$-$L^\infty$ effect of  Subsection \ref{sec.smooth1}, and the fundamental solutions of Subsection \ref{ssec.bddness} do not exist. The self-similar variables of Subsection \ref{ssec.ssv} will be replaced by a different type of self-similarity called inverse-time self-similarity or self-similarity of second type, a well-known type in the specialized literature, \cite{Barbk96, BerrHoll, GaKing2002, GaPel1997, King93, VazSmooth}.

\noindent\subsection{VSS  and nonlinear elliptic eigenvalue problem}

By analogy with Section \ref{sec.vss1} we construct a singular solution that we will still call  the VSS (Very Singular Solution), even if it not so singular in this range of $p$. The spatial part of this solution satisfies a nonlinear elliptic eigenvalue problem, that we state first.

\begin{theorem}\label{cor.eigen.vf}   Let $1<p<p_c$. Given $\lambda>0$, the function $F=A\,|y|^{-sp/(2-p)}>0$,  satisfies the nonlinear and singular elliptic  problem
\begin{equation}\label{eq.eigen.vf}
\mathcal L_{s,p} F(y)= \lambda F(y)  \qquad \mbox{for } \ y\ne 0,
\end{equation}
if the following condition holds: \
$ A^{2-p}\lambda=k(s,p,N).$ \
Note that \ $F(0)=+\infty$ with an integrable singularity.  The solution is not integrable at infinity (of the space variable) since \ $sp/(2-p)<N$ precisely for $p<p_c$. On the other hand, in this range $F$  is a weak solution in the whole space: for every smooth test function $\varphi(x)$ with compact support  we have
\begin{equation}
\begin{array}{c}
\displaystyle \iint \frac{\Phi_p(F(x)-F(y))\,(\varphi(x)-\varphi(y))}{|x-y|^{N+sp}}\,dxdy =
\displaystyle 2\lambda \iint F(x)\,\varphi(x,t)\,dx.
\end{array}
\end{equation}
\end{theorem}

\noindent Here, the universal constant $k(s,p,N)$ is the value of $\mathcal L_{s,p}(|y|^{-sp/(2-p)})$ at $|y|=1$.  By approximation we can use $\varphi$ belonging to $L^1(\ren)\cap L^\infty(\ren)$.

The associated evolution result says:

\begin{theorem}\label{thm.vss-vfr} Let $1<p<p_c$. There exists a constant $C_\infty(s,p,N) >0$ such that for any $T>0$
the function \begin{equation}\label{decay.vss.u}
U(x,t)=\,C_\infty  \,(T-t)^{1/(2-p)} { |x|^{-sp/(2-p)} }
\end{equation}
is a singular solution of the FPLEE at all times $0<t<T$ and points $x\ne 0$ with an integrable singularity at $r=0$.  We also have $C_\infty^{2-p}=(2-p)\,k$, with $k$ as before. The solution is not integrable at infinity (of the space variable) since \ $sp/(2-p)<N$ precisely for $p<p_c$. We may continue $U$ by $0$ for all $t\ge T$. $U$  is a weak solution of the FPLEE in the sense of \eqref{weak.sol.p}.
\end{theorem}

\noindent {\bf Remark.} Note that whenever  $qsp/(2-p)<N$, then \ $ F\in L_{loc}^q(\ren)$ with exponent $q$ in the range
$$
1\le q<q_*(s,p)=\frac{N(2-p)}{sp}.
$$
We have $q_*>1$ if $p<p_c$.

\medskip

\noindent \sl Proof of both results. \rm   We write the candidate solution for the FPLEE in the separate-variables form
$$
U(x,t)=C  \,(T-t)^{1/(2-p)}F(r), \  \quad  r= |x|,   \quad C>0.
$$
The equation for $F$ is then
$$
C^{p-1}{\mathcal L}_{s,p} F(r)= \frac{C}{2-p}\ F(r),
$$
which is just the elliptic equation \eqref{eq.eigen.vf} with  $(2-p)\lambda=C^{2-p}$\nc.

\noindent  We have to do some careful calculations with the candidate function $F=r^{-sp/(2-p)}$  for $p<p_c$. We first check that ${\mathcal L}_{s,p} F(r)$  is finite for every $r>0$, then by scaling, as we did before, we must have ${\mathcal L}_{s,p} F= kF$  \ for some  universal constant $k$ that is finite. This is just a consequence of scaling. Three possibilities arise: $k>0$, $k=0$ and $k<0$.

(I) In  case $k>0$ we  obtain a singular solution of the evolution equation of the form
$$
U(x,t)=C\,\frac{(T-t)^{1/(2-p)}}{|x|^{-sp/(2-p)}}
$$
for any $T>0$ by just adjusting the constant so that $C_\infty^{2-p}=(2-p)\,k$. This is the VSS that we want.
And we also solve the elliptic problem.

In order to proceed we need the following lemma.

\begin{lemma}
$F$ is a weak solution of the stationary equation, i.e.,
\begin{equation}\label{weaksol.vass.vf}
2k\int F(x)\varphi(x)\,dx = \iint   \frac{\Phi_p(F(x)-F(y))\,(\varphi(x)-\varphi(y))}{|x-y|^{N+sp}}\,dxdt.
\end{equation}
for every smooth  radial \nc function with compact support.
\end{lemma}

\noindent {\sl Proof. }   Since $\mathcal L_{s,p} F(x)=kF$ in the classical sense for $x\ne 0$, the only problem in proving \eqref{weaksol.vass.vf} is the singularity at  $x=y$ in the right-hand side. Now, given $x_0\in \ren$, $|x_0|=1$, for every small $\ve>0$ there is a $\delta_0<1 $ such that for $0<\delta<\delta_0$
$$
|\int_{B_\delta(x_0)} \frac{\Phi(F(x_0)-F(y))}{|x_0-y|^{N+sp}}\,dy |<\ve
$$
By scaling and using the power form of $F$, we get for every $x\ne 0$, not necessarily $|x|=1$,
$$
|\int_{B_\delta(x)} \frac{\Phi(F(x)-F(y))}{|x-y|^{N+sp}}\,dy |<\ve |x|^{-sp/(2-p)}=\ve F(x).
$$
On the other hand, we know that the full integral is $\mathcal L_{s,p} F(x)=kF(x)$, so that putting
$\Omega_\delta(x)=\ren\setminus B_\delta(x)$ we get
$$
|kF(x) -\int_{\Omega_\delta(x)} \frac{\Phi(F(x)-F(y))}{|x-y|^{N+sp}}\,dy |<\ve F(x)
$$
Integrating in $x$ over $\ren$ after multiplying by $\varphi(x)$, we have
$$
|k \int F(x)\varphi(x)\,dx -\iint_{\Omega_{\delta,1}} \varphi(x)\frac{\Phi(F(x)-F(y))}{|x-y|^{N+sp}}\,dxdy |
<\ve \int kF(x)|\varphi(x)|\,dx,
$$
where $\Omega_{\delta,1}=\{(x,y): |x-y|\le \delta|x|\} $. Likewise,
$$
|k \int F(y)\varphi(y)\,dx -\iint_{\Omega_{\delta,2}} \varphi(y)\frac{\Phi(F(y)-F(x))}{|x-y|^{N+sp}}\,dxdy |
<\ve \int kF(y)|\varphi(y)|\,dx
$$
in $\Omega_{\delta,2}=\{(x,y): |x-y|\le \delta|y|\} $. Adding both expressions, we get
$$
\big|2k \int F(y)\varphi(y)\,dx- \iint_{\Omega_{\delta,1}} \varphi(x)\frac{\Phi(F(x)-F(y))}{|x-y|^{N+sp}} -\iint_{\Omega_{\delta,2}} \varphi(y)\frac{\Phi(F(x)-F(y))}{|x-y|^{N+sp}}\,dxdy \big|
$$
$$
<2\ve k \int F(x)|\varphi(x)|\,dx.
$$
 Now observe that if we slight restrict the domain to $  \Omega_{\delta}=\Omega_{\delta,1} \cap \Omega_{\delta,2}$, we get
$$
\big|2k \int F\varphi\,dx- \iint_{\Omega_{\delta}} \frac{\Phi(F(x)-F(y))}{|x-y|^{N+sp}}(\varphi(x)-\varphi(y) \,dxdy \big|
<2\ve k \int F|\varphi|\,dx + J_1+J_2,
$$
where
$$
J_1=|\iint_{\Omega_{\delta,1}\setminus \Omega_{\delta}} \varphi(x)\frac{\Phi(F(x)-F(y))}{|x-y|^{N+sp}}\,dxdy|\,,
$$
and $J_2$ is its symmetric  version after changing $x$ and $y$. Next, we observe that for $(x,y)\in \Omega_{\delta,1}\setminus \Omega_{\delta}$
we have \ $\delta(1-\delta) |x|\le |x-y|\le \delta |x|$, \ and this integral can be computed in absolute value as follows:
$$
J_1\le C\int |\varphi(x)| dx\int_{A_\delta(x)} \frac{|\Phi(F(x)-F(y))|}{|x-y|^{N+sp}}\,dy
$$
where $A_\delta(x)=\{y: \delta|x|/2\le |x-y|\le \delta |x|\}$. Hence
$$
J_1\le C\int |\varphi(x)| |x|^{-\mu} \big(\int_{A_\delta(1)} \frac{|\Phi(F(x)-F(y))|}{|x-y|^{N+sp}}\,dy\big)\,dx\le
C_1\delta \int |\varphi(x)| |x|^{-\mu}\,dx.
$$
We conclude that
$$
\big|\iint_{\Omega_{\delta}} \frac{\Phi_p(F(x)-F(y))\,(\varphi(x)-\varphi(y))}{|x-y|^{N+sp}} \,dxdy -2k \int F\varphi\,dx\big|
<C_2\ve +C_3 \delta.
$$
Taking $\ve\to 0$, hence $\delta\to 0$ we get the result. \qed

\medskip

(II) We resume the proof. We want to eliminate the possibility that ${\mathcal L}_{s,p} F= 0$.  In order to prove that this cannot happen we recall an integral version of the maximum principe that uses the fact that $F$ is locally in $L^q$.  Given $c>0$  we use the fact that $U(x)\le c/2 $ for $|x|\ge R$ to compare $U$ with $c$ in the ball $B_R(0)$. We have ${\mathcal L}_{s,p} F= {\mathcal L}_{s,p} c= 0$. Therefore, using the definition of weak solution,
$$
0=\int {\mathcal L}_{s,p} F(x)\,\varphi(x)\,dx= \iint   \frac{\Phi_p(F(x)-F(y)\,(\varphi(x)-\varphi(y))}{|x-y|^{N+sp}}\,dxdt.
$$
Since $F$ is locally in $L^q$ for some $q>1$, for $\gamma$ small  $\varphi(x)=(F(x)-c)^\gamma_+$  is  a valid test function \nc (integrability of the integrand is guaranteed), and we have
$$
\iint \frac{\Phi_p(F(x)-F(y)\,((F(x)-c)^\gamma_+ -(F(y)-c)^\gamma_+))}{|x-y|^{N+sp}}\,dxdt=0.
$$
After analyzing the nonnegative integrand, we conclude that either $F$ is constant or  $F(x)\le c$ everywhere. Since $c>0$ is arbitrary, we conclude that $F=0$, a contradiction with the explicit formula for $F$.

(III)  It only remains to prove that we cannot have $k<0$. In that case we cannot used inverse-time self-similarity but plain self-similarity as in Section \eqref{sec.vss1}. We would obtain an increasing solution
$$
U(x,t)=C\,t^{1/(2-p)}\,|x|^{-sp/(2-p)}
$$
for any $T>0$ by just adjusting the constant $C^{2-p}=(2-p)\,|k|$.  This power growth is not possible for locally integrable functions. The argument is as follows:  we use vertical displacement downwards, and prove that $(U-\ve)_+$ must be bounded in $L^1(\ren)$ for all times. The Theorem is proved. \qed 

\medskip

\noindent {\bf Corollary 1.} Using $(F-\ve)_+ $ as a subsolution, we can prove that there are solutions that preserve
the singularity at a fixed point for some time. We conclude that there is no general smoothing effect $L^q\to L^\infty$ for
$1\le q<q_*$ in the very fast $p$ range.

\medskip

\noindent {\bf  Corollary 2.} The maximum principle allows to identify a class of solutions that extinguish in finite time,
i.e., those satisfying
$$
0\le u_0(x)\le C\, |x|^{-sp/(2-p)}
$$
for some $C>0$. This includes in particular all bounded data with compact support.

\noindent {\bf Remark.} A separated-variables solution with a fixed isolated singularity
exists  for the FPME \eqref{eq1.fpme} in the range  $m < m_c $, and  we get similar finite-time extinction profiles; the construction was performed in \cite{VolVa15}. We recall that for that equation the calculations are quite simple once the candidate has been identified.


\noindent\subsection{Nonexistence of fundamental solutions}

A basic tool in this topic is the scaling transformation \eqref{scal.trn1}
\begin{equation*}
{\mathcal T}_k u(x,t)= k^{N} u(kx, k^{N(p-2)+sp} t),\quad k>0\,.
\end{equation*}
For $p>p_c$ it allows to concentrate a given initial datum into a Dirac delta with the same mass as $k\to \infty$, while for $t>0$ ${\mathcal T}_k u$ looks for the value of the original solution far away in time, since $\gamma= N(p-2)+sp>0$. This was the main idea behind proof the existence of a fundamental solution in Theorem \ref{thm.exfs}.

Everything changes for lower values of $p$. Thus, in the case $p=p_c$, we get $\gamma=0$ and the transformation becomes
\begin{equation}\label{scal.trn1.pc}
({\mathcal T}_k u)(x,t)= k^{N} u(kx, t),\quad h>0\,.
\end{equation}
This transformation still preserves mass and in this case it also preserves time. So it is unable to look for large values of time and bring them back to $t\approx 1$. Even worse, in the fast diffusion case $1<p<p_c$ we have
$\gamma=N(p-2)+sp<0$ to that for $k>1$ ${\mathcal T}_k u $ looks at contracted time instead of looking to large times.

Our next result shows the failure of the approximation method to produce the existence of a fundamental solution.

\begin{theorem}\label{thm.nosss} Let $u$ be a semigroup solution of the FPLEE for exponent $1<p\le p_c$ with nonnegative, bounded and compactly supported initial data $u_0$ with mass $M>0$. Let $u_k$ the rescaled solution, so that $u_{0k}={\mathcal T}_k u_0$ converges $M$ times to the Dirac delta
Then for every $t>0$ we have
\begin{equation}\label{scal.as.dirac}
u_k(x,t)\to M\delta(x)\,.
\end{equation}
In other words, the limit solution preserves the initial point mass unchanged in time.
\end{theorem}

\noindent {\sl Proof.} In case $p=p_c$ we only need to note that
$$
u_k(\cdot,t)=S_t({\mathcal T}_k u_0)={\mathcal T}_k(S_t u_0).
$$

In the fast diffusion case $1<p<p_c$ we have
$$
u_k(\cdot,t)=S_t({\mathcal T}_k u_0)={\mathcal T}_k(S_{k^{-\gamma}t} u_0).
$$
But $S_{k^{-\gamma}t} u_0 \to u_0$  as $k\to\infty$ in all $L^q(\ren)$ so that the same conclusion holds with faster convergence. \qed

We will repeat a comment taken from the study of very fast diffusion for the PME in \cite{VazSmooth}. The term fast diffusion is accurate for small densities $u$ but for large densities rather the opposite happens, the equation becomes poorly diffusive and in particular, it is unable of minimally spreading a Dirac mass. While the explanation in \cite{VazSmooth} was based on a detailed physical analysis, it is only given here as a intuition for the proven fact that Dirac deltas do not spread at all.

It might be argued that a fundamental solution could come out from some other contrived construction. Let us point out that any method that leads to a self-similar solution is doomed to fail since the similarity exponents would be given by algebraic reasons by
formulas  \eqref{eq.sse1},  so that $\alpha=\beta=\infty$ for $p=p_c$, and they would become negative for $p<p_c$.

\medskip



\noindent\subsection{A note on the case $p=p_c$}\label{sec.pc}

The basic theory still applies as in previous section, and we have just shown a negative result about FS that includes $p=p_c$
with simple proof. Other topics for this limit case $p=p_c$ are more difficult. By analogy with the VSS solution profiles found for $p>p_c$ and $p<p_c$, the  candidate function becomes $F(x)=|x|^{-N}$ which is singular at $r=0$ with a barely non-integrable rate. This does not cause a problem for the existence of $\mathcal L_{s,p}(|x|^{-N})$ since the defining integral depends on \ $F^{p-1}$ and $p-1$ is less than 1.

We prove the following result that complements what is known for $p>p_c$, cf. Theorem \eqref{cor.eigen}, and $p<p_c$, cf. Theorem \eqref{cor.eigen.vf}.

\begin{theorem} Let $N\ge 1$, $0<s<1$. For the precise value $p=2N/(N+s)$ function $F(x)=|x|^{-N}$ solves the nonlinear elliptic equation
\begin{equation}\label{Fundsol.ell}
\mathcal L_{s,p}({-|x|^N})=0
\end{equation}
in the classical sense in $\ren\setminus \{0\}$.
\end{theorem}

\noindent \sl Proof. \rm  
(i) First of all, we argue as in Sections \ref{sec.barr1} and \ref{sec.barr2} that the integral is finite. Recall that for $x\approx y$ the integral is delicate and we need the results of \cite{KKL}, Section 3, or  \cite{DTGCV}, Lemma A.2.
We call $A=\left.\mathcal L_{s,p_c}(r^{-N})\right|_{r=1}$.

Next, we observe that for $p=p_c$ the  scaling transformation \eqref{scal.trn1} becomes for stationary solutions
\begin{equation}\label{scal.trn1,pc}
{\mathcal T}_k F(x)= k^{N} F(kx),\quad k>0\,,
\end{equation}
and obeys the rule  \ $ \mathcal L_{s,p}{\mathcal T}_k F(x)=k^{N}(\mathcal L_{s,p}F)(kx)$.
Since ${\mathcal T}_k $ leaves invariant the function $F(r)=r^{-N}$, we conclude that for every $r_0=|x|\ne 1$, and putting
$k=r_0^{-1}$,
we have
$$
\left. \mathcal L_{s,p_c}(r^{-N})\right|_{r=r_0}=\left. r_0^{-N}\mathcal L_{s,p_c}(r^{-N})\right\|_{r=1}
=A\,r_0^{-N}.
$$
Next step is to prove that $\mathcal L_{s,p}(|x|^{-N})$ vanishes at $r>0$ for $p=p_c$. We only need to check it at $r=1$, i.e., $A=0$. We get the result from  the continuity of the integral $\mathcal L_{s,p}(|x|^{-sp/(2-p)})$ at $|x|=1$ with respect the variation of the parameter $p$ around $p_c$ for fixed $s\in (0,1)$. Recall that  $sp/(2-p)=N$ for $p=p_c$. The continuity of this integral follows by continuous dependence of the integrand away from the problematic regions. To be precise, we only need to check the uniform integrability of the integral expression
$$
A_{p}=\int_{\ren}\frac{\Phi(|1-|y|^{-sp/(2-p)})}{|x-y|^{N+sp}}\,dy
$$
for $p\approx p_c$ and $|x|=1$, away from $y=0$ and away from infinity. Uniform integrability at infinity is easy.  For the uniform bound around the singularity at $x= y$ we again refer to \cite{KKL}.

We may now recall the already proved facts that $\mathcal L_{s,p}(r^{-sp/(2-p)})>0 $ for $1<p<p_c$ and is negative for $p_c<p<p_1$ to conclude that, necessarily, $\mathcal L_{s,p_c}(r^{-N)})=0$ for all $r>0$. \qed

We stop the investigation of this case at this moment since the very fast range leads to many novelties, worth a detailed study.
For the reader's sake we add a comparative comment on related (simpler) models.

 (i) Putting $s=1$, we are led consider the stationary  $p$-Laplacian equation with $1<p<2$ (fast range). Then we have $p=p_c=2N/(N+1)$, \cite{VazSmooth}. We use the candidate function $F(r)=r^{-N}$. By explicit computation we see that $\Delta_p({-r^N})=0$  for all $r> 0$, in agreement with our results. Moreover, the flux function $r^{N-1}(F')^{p_1} =c$, so that there is a Dirac delta, i.e.,
$$
-\Delta_p(|x|^{-N})=c\delta(x)
$$
in the sense of weak solutions.

(ii) We argue in a similar way for the PME with $N>2$: $u_t=\Delta u^m$ with $m>0$. This is again a nonlinear non-fractional equation and the critical exponent  $m=m_c=(N-2)/N$. We have the candidate $F=r^{-N}$ as limit of the profiles of the VSS for other exponents $m<1$. Then, $F^m=r^{-(N-2)}$ for $m=m_c$, hence $-\Delta (F^{m_c})$ is a delta. This is simply saying that $F$ is the fundamental solution for the Laplacian!

(iii) Our next example is the fractional version, FPME, where  we  have  $0<s<1$, and  $m=m_c=(N-2s)/N$, still a fast diffusion exponent ($0<m_c<1$). If we take $F=r^{-N}$ we have $F^m=r^{-(N-2s)}$, and it is clear that  $(-\Delta)^s (F^m)$ is zero for $r\ne 0$, see the Fourier analysis computations in  \cite{VPLP2020}.

Proving that there is a Dirac delta at the origin in our problem is an interesting issue, but we will not discuss it here since it leads to new technical developments\nc.

\section{Comments and extensions}

%

\noindent $\bullet$ The  exact tail behaviour  of the fundamental solution is very precise in the lower range $p_c<p<p_1$, but  it might be improved for $p_1\le p<2$. Numerical computations suggest a finer decay expression $F_1(x)\sim  C\, |x|^{N+sp}$. Of course, in this case $C_M$ depends on $M$. \nc

\noindent $\bullet$ The lack of agreement between the decay rates of the upper and lower bound in the critical case $p=p_1$ has prevented us from obtaining a Global Harnack Principle for that critical exponent. This is an interesting open problem.

\noindent $\bullet$ We have proved uniqueness of the self-similar fundamental solution. The uniqueness of the general fundamental solution is a delicate issue that we did not address here.

\noindent $\bullet$
 The question of uniqueness of the VSS and the corresponding nonlinear elliptic problem remains open. \nc

\noindent $\bullet$ Existence of solutions for measures as initial data should be investigated. See \cite{KMSire}. This is related to the question of initial traces.

\noindent $\bullet$ The question of rates of convergence for the asymptotic result \eqref{lim.ab.L1} of Theorem  \ref{thm.ab1} has not been considered. This important issue has been addressed for many other models of nonlinear diffusion. It is solved for many of them, but well known cases remain open.

\noindent $\bullet$ We did not consider in sufficient detail the case where $1<p<p_c$, where the fundamental solution does not exist. There is extinction in finite time  for many integrable solutions, see above and also \cite{BS2020}. As we have said, the study of this range needs further attention.

\noindent $\bullet$  In the existence theory we can consider wider classes of initial data, possibly growing at infinity. Optimal classes are known in the linear fractional equation (case $p=2$), \cite{BSV17}, and in the standard $p$-Laplacian equation (case $s=1$), cf. \cite{DiBeHerr89}. Of course, the asymptotic behaviour will not be the same as the one exhibited in this paper in terms of the fundamental solutions.

\noindent $\bullet$ Another interesting open issue is the presence of a right-hand side in the equation, maybe in the form of lower-order teems.  There are some works, see e.g. 
\cite{Teng2019} and its references.\nc

\noindent $\bullet$ We have considered a nonlinear equation of fractional type with nonlinearity \\ $\Phi(u)=|u|^{p-2}u$, and we have strongly used  in a number of technical steps the fact that $\Phi$ is a power, hence homogeneous. We wonder how much of the theory holds for more general monotone nonlinearities $\Phi$.


\vskip 1cm

\noindent {\textbf{\large \sc Acknowledgments.}} Author partially funded by Projects MTM2014-52240-P and   PGC2018-098440-B-I00 (Spain). Partially performed as an Honorary Professor at Univ. Complutense de Madrid.
 The author thanks F. del Teso for the numerical treatment that led to the self-similar profile displayed at the end of Section \ref{sec.fs}. He is also grateful to  A. Iannizzotto for information and discussions about his paper \cite{IannMS2016} that led to the continuity results we present here.

\medskip


{\small
\bibliographystyle{amsplain}

}

\medskip

\noindent {\sc Address:}

\noindent Juan Luis V\'azquez. Departamento de Matem\'{a}ticas, Universidad
Aut\'{o}noma de Madrid,\\ Campus de Cantoblanco, 28049 Madrid, Spain.  \\
e-mail address:~\texttt{juanluis.vazquez@uam.es}

\

\noindent {\bf Keywords: } Nonlinear parabolic equations, $p$-Laplacian operator, fractional operators, fundamental solutions,  asymptotic behaviour, very singular solutions.

\medskip

\noindent {\bf 2020 Mathematics Subject Classification.}
  	35K55,  	
   	35K65,   	
    35R11,   	
    35A08,   	
    35B40.   	

\vskip 1cm

\end{document}